\documentclass[preprint,12pt]{elsarticle}

\usepackage{hyperref}
\usepackage{graphicx}
\usepackage{subcaption}
\usepackage{amssymb}
\usepackage{amsmath}
\usepackage{multirow}
\usepackage[utf8]{inputenc}
\usepackage{cleveref}
\usepackage[section]{placeins}
\usepackage[ruled,vlined]{algorithm2e}
\usepackage{listings}
\usepackage{enumitem}
\usepackage[mathscr]{euscript}

\usepackage{xcolor}
\usepackage{xargs}
\usepackage[colorinlistoftodos,textsize=footnotesize]{todonotes}
\usepackage{lineno}
\usepackage{tikz}

\newcommandx{\unsure}[2][1=]{\todo[linecolor=red,backgroundcolor=red!25,bordercolor=red,#1]{#2}}
\newcommandx{\change}[2][1=]{\todo[linecolor=blue,backgroundcolor=blue!25,bordercolor=blue,#1]{#2}}
\newcommandx{\info}[2][1=]{\todo[linecolor=OliveGreen,backgroundcolor=OliveGreen!25,bordercolor=OliveGreen,#1]{#2}}

\newcommand{\ten}[1]{\ensuremath{\mathbf{#1}}}

\tikzset{
    cross/.pic = {
    \draw[rotate = 45] (-#1,0) -- (#1,0);
    \draw[rotate = 45] (0,-#1) -- (0, #1);
    }
}

\newcommand\fluids{
    \fill[blue!60](0.1,0.1) circle (1ex);
    \fill[blue!60](-0.05,0.9) circle (1ex);
    \fill[blue!60](0.1,1.95) circle (1ex);
    \fill[blue!60](0,3.0) circle (1ex);
    \fill[blue!60](-0.05,3.9) circle (1ex);
    \fill[blue!60](-0.1,4.8) circle (1ex);
    \fill[blue!60](-1.1,-0.1) circle (1ex);
    \fill[blue!60](-1.2,1.05) circle (1ex);
    \fill[blue!60](-0.9,2.05) circle (1ex);
    \fill[blue!60](-0.9,2.95) circle (1ex);
    \fill[blue!60](-1.1,4.1) circle (1ex);
    \fill[blue!60](-0.95,5.05) circle (1ex);
    \fill[blue!60](-2,0) circle (1ex);
    \fill[blue!60](-2,1.0) circle (1ex);
    \fill[blue!60](-1.95,2.0) circle (1ex);
    \fill[blue!60](-2,2.95) circle (1ex);
    \fill[blue!60](-1.9,4.0) circle (1ex);
    \fill[blue!60](-2,5.0) circle (1ex);
    \draw[dashed] (0.5,-0.5) -- (0.5,5.5) node[yshift=0.5cm] {Interface};
    \draw (-1.0, 5.5) node{Fluid};
    \draw (2.0, 5.5) node{Solid};
}

\newcommand\solids{
    \fill[red!60](1,0) circle (1ex);
    \fill[red!60](1,1) circle (1ex);
    \fill[red!60](1,2) circle (1ex);
    \fill[red!60](1,3) circle (1ex);
    \fill[red!60](1,4.0) circle (1ex);
    \fill[red!60](1,5.0) circle (1ex);
    \fill[red!60](2,0) circle (1ex);
    \fill[red!60](2,1.) circle (1ex);
    \fill[red!60](2,2.0) circle (1ex);
    \fill[red!60](2,3) circle (1ex);
    \fill[red!60](2,4.0) circle (1ex);
    \fill[red!60](2,5.0) circle (1ex);
    \fill[red!60](3,0) circle (1ex);
    \fill[red!60](3,1.0) circle (1ex);
    \fill[red!60](3,2.0) circle (1ex);
    \fill[red!60](3,3) circle (1ex);
    \fill[red!60](3,4.0) circle (1ex);
    \fill[red!60](3,5.0) circle (1ex);
}

\newcommand\io{
    \fill[blue!90](0.9549218423,-0.04221) circle(1ex);
    \fill[blue!90](0.9563673372,0.966944) circle(1ex);
    \fill[blue!90](0.9668138534,2.045488) circle(1ex);
    \fill[blue!90](1.011271014,2.9962833) circle(1ex);
    \fill[blue!90](0.9729329842,4.040094) circle(1ex);
    \fill[blue!90](1.015977573,5.0373940) circle(1ex);
    \fill[blue!90](2.023183498,-0.035654) circle(1ex);
    \fill[blue!90](1.972560509,1.0005706) circle(1ex);
    \fill[blue!90](1.959527978,1.9692822) circle(1ex);
    \fill[blue!90](1.994235407,3.0133835) circle(1ex);
    \fill[blue!90](2.011956196,4.0432622) circle(1ex);
    \fill[blue!90](1.989535452,5.0236606) circle(1ex);
    \fill[blue!90](2.97682717,0.03333736) circle(1ex);
    \fill[blue!90](2.95544884,1.00277262) circle(1ex);
    \fill[blue!90](3.038032744,2.0337744) circle(1ex);
    \fill[blue!90](3.024208056,3.0235824) circle(1ex);
    \fill[blue!90](2.97974001,3.96700024) circle(1ex);
    \fill[blue!90](2.952165244,5.0348294) circle(1ex);

    \fill[blue!60](0.01797076887,0.02027) circle(1ex);
    \fill[blue!60](0.04368048786,0.98296) circle(1ex);
    \fill[blue!60](-0.04902756808,2.0484) circle(1ex);
    \fill[blue!60](0.03860085251,2.95145) circle(1ex);
    \fill[blue!60](0.03829308234,3.98986) circle(1ex);
    \fill[blue!60](0.03069009269,5.00578) circle(1ex);
    \fill[blue!60](-1.027177474,0.013286) circle(1ex);
    \fill[blue!60](-0.975182599,0.985942) circle(1ex);
    \fill[blue!60](-0.9950074012,1.95265) circle(1ex);
    \fill[blue!60](-0.9670742206,2.96435) circle(1ex);
    \fill[blue!60](-0.9759024185,4.04093) circle(1ex);
    \fill[blue!60](-1.020562935,4.950827) circle(1ex);
    \fill[blue!60](-1.952384278,0.011329) circle(1ex);
    \fill[blue!60](-2.035966346,1.041900) circle(1ex);
    \fill[blue!60](-1.971135398,2.044477) circle(1ex);
    \fill[blue!60](-2.031527974,2.973417) circle(1ex);
    \fill[blue!60](-2.031252392,4.044029) circle(1ex);
    \fill[blue!60](-1.999422362,4.984838) circle(1ex);
    \fill[blue!60](-2.974749829,-0.02605) circle(1ex);
    \fill[blue!60](-2.974244968,0.963327) circle(1ex);
    \fill[blue!60](-2.989710956,1.982291) circle(1ex);
    \fill[blue!60](-2.971306524,3.026063) circle(1ex);
    \fill[blue!60](-3.032865433,4.025190) circle(1ex);
    \fill[blue!60](-3.013507329,4.969604) circle(1ex);
    \fill[blue!60](-3.958393846,0.037311) circle(1ex);
    \fill[blue!60](-4.018227797,0.965342) circle(1ex);
    \fill[blue!60](-3.980661133,2.031720) circle(1ex);
    \fill[blue!60](-4.020376981,2.986839) circle(1ex);
    \fill[blue!60](-4.038549523,3.951917) circle(1ex);
    \fill[blue!60](-4.023768211,4.996410) circle(1ex);

    \fill[blue!30](-5,0) circle (1ex);
    \fill[blue!30](-5,1.0) circle (1ex);
    \fill[blue!30](-5,2.0) circle (1ex);
    \fill[blue!30](-5,3) circle (1ex);
    \fill[blue!30](-5,4.0) circle (1ex);
    \fill[blue!30](-5,5.0) circle (1ex);
    \fill[blue!30](-6,0) circle (1ex);
    \fill[blue!30](-6,1.0) circle (1ex);
    \fill[blue!30](-6,2.0) circle (1ex);
    \fill[blue!30](-6,3) circle (1ex);
    \fill[blue!30](-6,4.0) circle (1ex);
    \fill[blue!30](-6,5.0) circle (1ex);
    \fill[blue!30](-7,0) circle (1ex);
    \fill[blue!30](-7,1.0) circle (1ex);
    \fill[blue!30](-7,2.0) circle (1ex);
    \fill[blue!30](-7,3) circle (1ex);
    \fill[blue!30](-7,4.0) circle (1ex);
    \fill[blue!30](-7,5.0) circle (1ex);

\draw[dashed] (-7.5, 5.5) -- node[yshift=0.5cm]{Fluid}(3.5, 5.5)
node[yshift=0.5cm, xshift=-1.5cm]{Outlet};
    \draw[dashed] (-7.5, -0.5) -- (3.5, -0.5);
    \draw[dashed] (-7.5, -0.5) -- node[yshift=3.5cm, xshift=1.5cm]{Inlet}(-7.5, 5.5);
    \draw[dashed] (-4.5, -0.5) -- (-4.5, 5.5);
    \draw[dashed] (0.5, -0.5) -- (0.5, 5.5);
    \draw[dashed] (3.5, -0.5) -- (3.5, 5.5);

}

\journal{}

\begin{document}

\begin{frontmatter}

  \title{How to train your solver: Verification of boundary
  conditions for smoothed particle hydrodynamics}

  \author[IITB]{Pawan Negi\corref{cor1}}
  \ead{pawan.n@aero.iitb.ac.in}
  \author[IITB]{Prabhu Ramachandran}
  \ead{prabhu@aero.iitb.ac.in}
\address[IITB]{Department of Aerospace Engineering, Indian Institute of
  Technology Bombay, Powai, Mumbai 400076}

\cortext[cor1]{Corresponding author}

\begin{abstract}
  The weakly compressible smoothed particle hydrodynamics (WCSPH) method has
  been employed to simulate various physical phenomena involving fluids and
  solids. Various methods have been proposed to implement the solid wall,
  inlet/outlet, and other boundary conditions. However, error estimation and
  the formal rates of convergence for these methods have not been discussed
  or examined carefully. In this paper, we use the method of manufactured
  solution (MMS) to verify the convergence properties of a variety of
  commonly employed of various solid, inlet, and outlet boundary
  implementations. In order to perform this study, we propose various
  manufactured solutions for different domains. On the basis of the
  convergence offered by these methods, we systematically propose a
  convergent WCSPH scheme along with suitable methods for implementing the
  boundary conditions. We also demonstrate the accuracy of the proposed
  scheme by using it to solve the flow past a circular cylinder. Along with
  other recent developments in the use of adaptive resolution, this paves the
  way for accurate and efficient simulation of incompressible or
  weakly-compressible fluid flows using the SPH method.

\end{abstract}


\end{frontmatter}


\section{Introduction}
\label{sec:intro}

The Smoothed particle hydrodynamics (SPH) method is widely used to solve
fluid dynamics problems\cite{violeau2012fluid}. One of the widely used
variants of the SPH method is the weakly compressible SPH (WCSPH). In
WCSPH, the pressure is obtained using an artificial equation of state
\cite{wcsph-state-of-the-art-2010}. Many WCSPH schemes viz.\ transport
velocity formulation (TVF) \cite{Adami2013}, entropically damped artificial
compressibility (EDAC) SPH \cite{edac-sph:cf:2019}, $\delta^{+}$-SPH
\cite{sun2019consistent}, and dual-time SPH
\cite{ramachandran_dual-time_2021} have been proposed in the last decade.
Recently, \citet{negi2021numerical} introduced several SPH schemes and
formally showed the convergence to be second order. However, they only
considered domains with periodic boundaries. For these schemes to be useful
in practical simulation it is essential to also have second order
convergent boundary condition implementations.

In the context of SPH, accurate implementation of boundary conditions is a
grand challenge problem~\cite{vacondio_grand_2020}. Many authors (See the
review by \citet{Violeau16}) have proposed various methods to implement solid
boundary conditions in SPH. Some authors
\cite{takeda1994a,maciaTheoreticalAnalysisNoSlip2011,randles1996smoothed,Adami2013}
use few layers of fixed solid (boundary) particles and use different
methods to extrapolate the properties from fluid to the solid particles.
\citet{colagross2003a} proposed the creation of solid particles by
reflecting the fluid particles about the solid-fluid interface and
retaining the properties. \citet{marrone-deltasph:cmame:2011} used fixed
ghost particles to implement the boundary condition wherein a reflection of
these ghost particles about the solid-fluid interface are used to evaluate
properties on the ghost particles. Recently, \citet{fourtakas2019} proposed
to use a dynamically generated local stencil for the particles near the
boundary. Some
authors~\cite{hashemi2012,ferrand2013,monaghan2009a,sph:fsf:monaghan-jcp94}
propose to use a single layer of particles to represent the solid boundary
and use methods to correctly evaluate the forces. Others like
\citet{ferrand2013} proposed semi-analytical methods to compensate for the
loss of kernel support near the boundary.

For the case of open boundaries, the use of a weakly compressible
formulation poses unique problems since pressure waves travel with an
artificial speed of sound. These waves must pass through the open
boundaries without reflecting into the domain. \citet{federico2012} propose
a do-nothing kind of outlet where the fluid particles are converted from
fluid to outlet particles while retaining their properties.
\citet{tafuni_io_2018} propose to mirror the inlet and outlet particles to
calculate properties and its gradient. The inlet and outlet properties are
set using a Taylor series expansion from the respective ghost particle to
the inlet/outlet particle. Recently, \citet{negi_nrbc_2020} proposed a
modified version of the method proposed by \citet{lastiwka_nrbc_2009} where
a time-averaged value is passed to the inlet/outlet and properties derived
from the characteristics of the flow are added using a Shepard
interpolation~\cite{shepard_1968}. Many authors
\cite{negi_nrbc_2020,maciaTheoreticalAnalysisNoSlip2011,valizadeh2015}
compared various boundary condition implementations qualitatively without
performing a convergence study. Furthermore, in the context of SPH, the
errors are often shown qualitatively by comparing the results of the
simulation \cite{Adami2013,edac-sph:cf:2019,sun2019consistent}.

In this paper, we verify the convergence of various boundary condition
implementations that have been proposed and identify second-order convergent
methods that are suitable for use for the most common solid wall, as well as
inlet and outlet boundary conditions. In order to identify a convergent
boundary implementation, we use the method of manufactured solutions (MMS)
introduced by \citet{negiHowTrainYour2021a} for SPH. In MMS, a manufactured
solution (MS) is created such that the boundary condition is satisfied at the
boundary interface of interest. Since the MS is not a solution of the weakly
compressible Navier-Stokes (NS) equation, we obtain a residue on substituting
the MS in the NS equations. This residue is added as a source term to the NS
equation in the scheme and it is expected that the solver will recover the MS
as the solution. In order to test the boundary condition, one requires a
convergent scheme in the bulk of the fluid for which we use the second-order
convergent scheme proposed by \citet{negi2021numerical} which is further
verified using MMS in \cite{negiHowTrainYour2021a}.

We believe that this is the first time that second order convergent
boundary condition implementations have been identified and verified in the
context of WCSPH. However, the implementation of these methods to solve a
real-life problem is non-trivial. In view of that, we propose a complete
algorithm and demonstrate the accuracy by solving a simple flow past a
circular cylinder problem. This shows that the resulting SPH scheme can be
applied to a variety of problems, paving the way for SPH to be an effective
alternative to traditional finite volume based codes especially when
combined with some recent advancements of adaptive resolution SPH
schemes~\cite{muta_efficient_2022,haftu_parallel_2022}.

In the next section, we discuss the SPH method and the second-order
discretizations. We then discuss various solid boundary condition
implementations, viz.\ pressure Neumann, slip and no-slip in \cref{sec:bc},
and open boundary conditions, viz.\ inlet and outlet in \cref{sec:obc} for
both pressure and velocity. In \cref{sec:mms_bc}, we discuss the MMS in
general and its construction for specific boundary condition. We compare
all the boundary condition implementation in the \cref{sec:results}
followed by conclusion in \cref{sec:conclusions}. All the results in this
manuscript are reproducible, and the source code for the simulations can be
found at \url{https://gitlab.com/pypr/mms_sph_bc}.

\section{The SPH method}
\label{sec:sph}

\citet{chorin_numerical_1967} proposed the weakly-compressible method to
solve fluid flow problems. In the weakly-compressible method, we solve the
governing equations given by
\begin{equation}
  \begin{split}
    \frac{d \varrho}{dt} &= -\varrho \nabla \cdot \ten{u}\\
    \frac{d \ten{u}}{dt } &= -\frac{\nabla p}{\varrho} + \nu \nabla^2 \ten{u},\\
  \end{split}
  \label{eq:wcsph}
\end{equation}
where $\varrho$ , $\ten{u}$, and $p$ are the density, velocity, and
pressure of the fluid, and $\nu$ is the dynamic viscosity of the fluid.
Additionally, we use an artificial equation of state (EOS) to link the
pressure with the density. In SPH, we use the EOS given by
\begin{equation}
  p = c_o^2 (\varrho - \varrho_o),
  \label{eq:eos}
\end{equation}
where $\varrho_o$ and $c_o$ are the reference density and artificial speed
of sound, respectively. In this paper, we use the L-IPST-C scheme proposed
by \citet{negi2021numerical} to discretize the governing equations. The
continuity equation is discretized as
\begin{equation}
  \frac{d \varrho_i}{ dt} = -\varrho_i \sum_j (\ten{u}_j - \ten{u}_i) \cdot
  \tilde{\nabla} W_{ij} \omega_j,
  \label{eq:cont_desc}
\end{equation}
where $\tilde{\nabla} W_{ij}$ is the kernel gradient corrected using the
correction proposed by \citet{bonet_lok:cmame:1999}, and $\omega_j =
\frac{1}{\sum_j W_{ij}}$ is the numerical volume. The function $W_{ij} =
W(\ten{x}_i - \ten{x}_j, h_{ij})$ is a smoothing kernel used in SPH, where
$\ten{x}_i$ is the position of the destination particle, $\ten{x}_j$ is the
position of the source particle, and $h_{ij} = 0.5 (h_i + h_j)$ is the
average smoothing radius\footnote{In this paper, we consider only constant
smoothing radius in the domain.}. The momentum equation is discretized as
\begin{equation}
  \frac{d \ten{u}_i}{dt} = \sum_j \left( \frac{ (p_i - p_j)}{\varrho_i}
  \tilde{\nabla} W_{ij} \omega_j + \nu (\nabla \ten{u}_j - \nabla \ten{u}_i)
  \cdot \tilde{\nabla} W_{ij} \omega_j \right),
  \label{eq:mom_desc}
\end{equation}
where $\nabla \ten{u}_i = \sum_j (\ten{u}_j - \ten{u}_i) \otimes
\tilde{\nabla} W_{ij} \omega_j$. The particles are shifted using the
iterative particle shifting technique proposed by \citet{huang_kernel_2019}
after every 10 timesteps. The particle properties are updated after
shifting using first-order Taylor series approximation given by
\begin{equation}
  \phi(\tilde{\ten{x}}) = \phi(\ten{x}) + (\tilde{\ten{x}} - \ten{x}) \nabla \phi(\ten{x}),
\end{equation}
where $\tilde{\ten{x}}$ is the position after shifting, $\ten{x}$ is the
position before shifting, and $\phi$ is any fluid property. We use a
second-order Runge-Kutta time integration scheme to integrate the
continuity and momentum equations. Since we are interested in the spatial
convergence of the scheme, we use a constant timestep corresponding to the
highest resolution, i.e. $500 \times 500$ for all the simulations given by
\begin{equation}
  \Delta t = \frac{h}{c_o + U} = \frac{1.2 \times 1/500}{20 + 1} = 0.00012\ sec.
\end{equation}
We assume a maximum velocity $U=1m/s$ and corresponding speed of sound $20
m/s$ for all our simulations.

In order to apply boundary conditions, many authors have proposed different
methods to implement Neumann pressure, slip, and no-slip boundary
conditions in SPH. The main objective is to extrapolate velocity and
pressure from the fluid particle to the ghost particle representing solid
such that the desired condition is satisfied. In the next section, we
discuss various boundary conditions in brief.

\section{Solid boundary conditions}
\label{sec:bc}

In this section, we discuss various boundary condition implementations widely
used in the SPH literature. We classify these implementations on the basis of
the requirement of the secondary particle arrays. In SPH, two types of
particles are used to implement boundary conditions, viz.\ ghost, and virtual
particles. The ghost particles carry the extrapolated properties from the
fluid and influence the fluid particles. Whereas the virtual particles are
used to evaluate some intermediate value of a property from fluid and do not
affect the actual flow.

Many authors~\cite{ferrand2013,hashemi2012,marongiu2007} have proposed
different approaches where a single layer of ghost particles are used. In
order to assess the effect of the number of layers on the accuracy of the
second order accurate gradient and Laplacian discretization, we perform a
simple numerical test. In this test, we consider a finite 2D domain of size
$1m \times 1m$. We discretize the domain using particles at different
resolutions, and initialize various properties using
\begin{equation}
  \begin{split}
  u(x, y) &= \sin(4 \pi (x +y)),\\
  v(x, y) &= \cos(4 \pi (x +y)), \\
  p(x, y) &= sin(4 \pi x) + sin(4 \pi y).
  \end{split}
  \label{eq:check}
\end{equation}

\begin{figure}[htbp]
  \centering
  \includegraphics[width=0.48\linewidth]{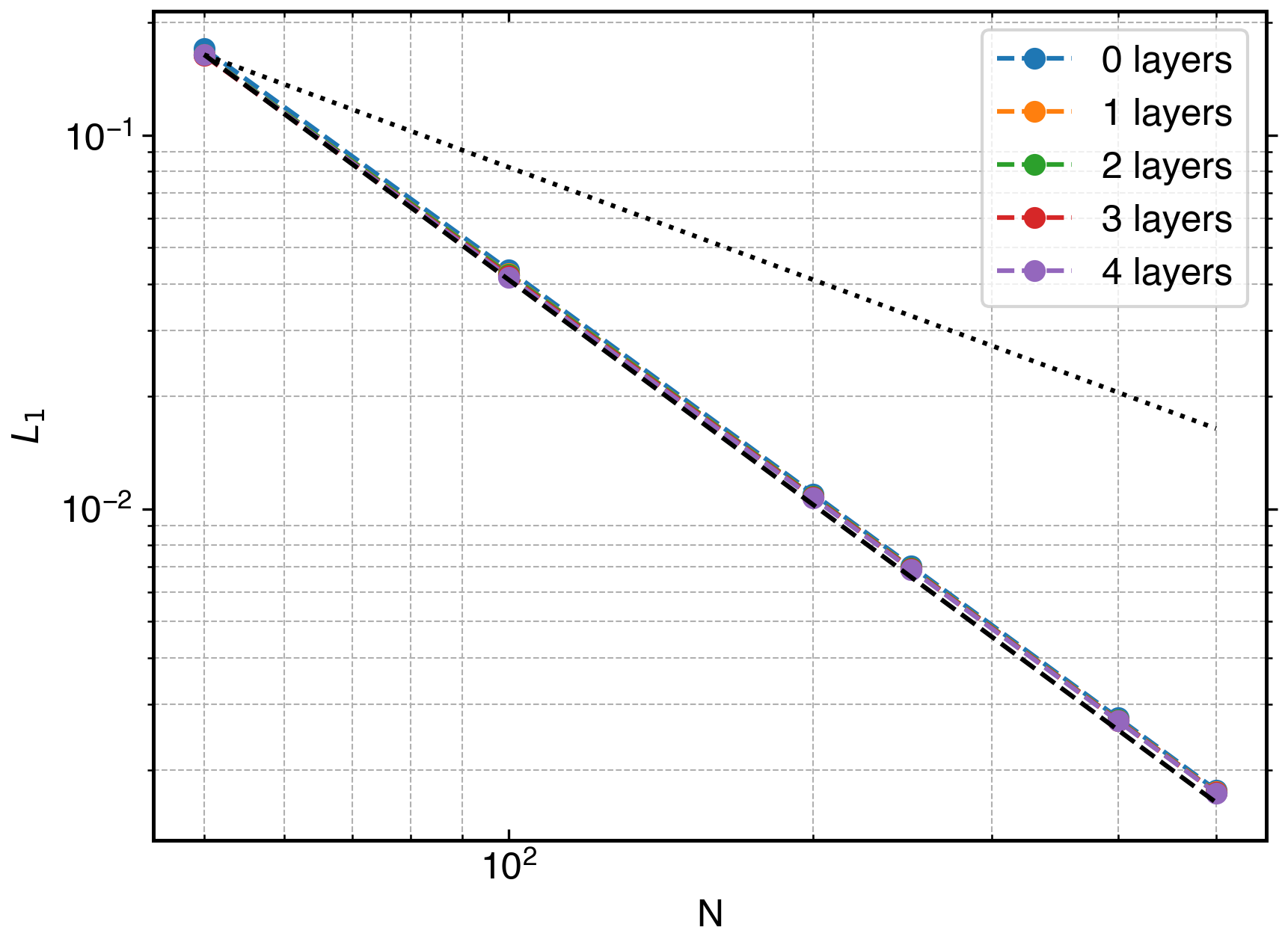}
  \includegraphics[width=0.48\linewidth]{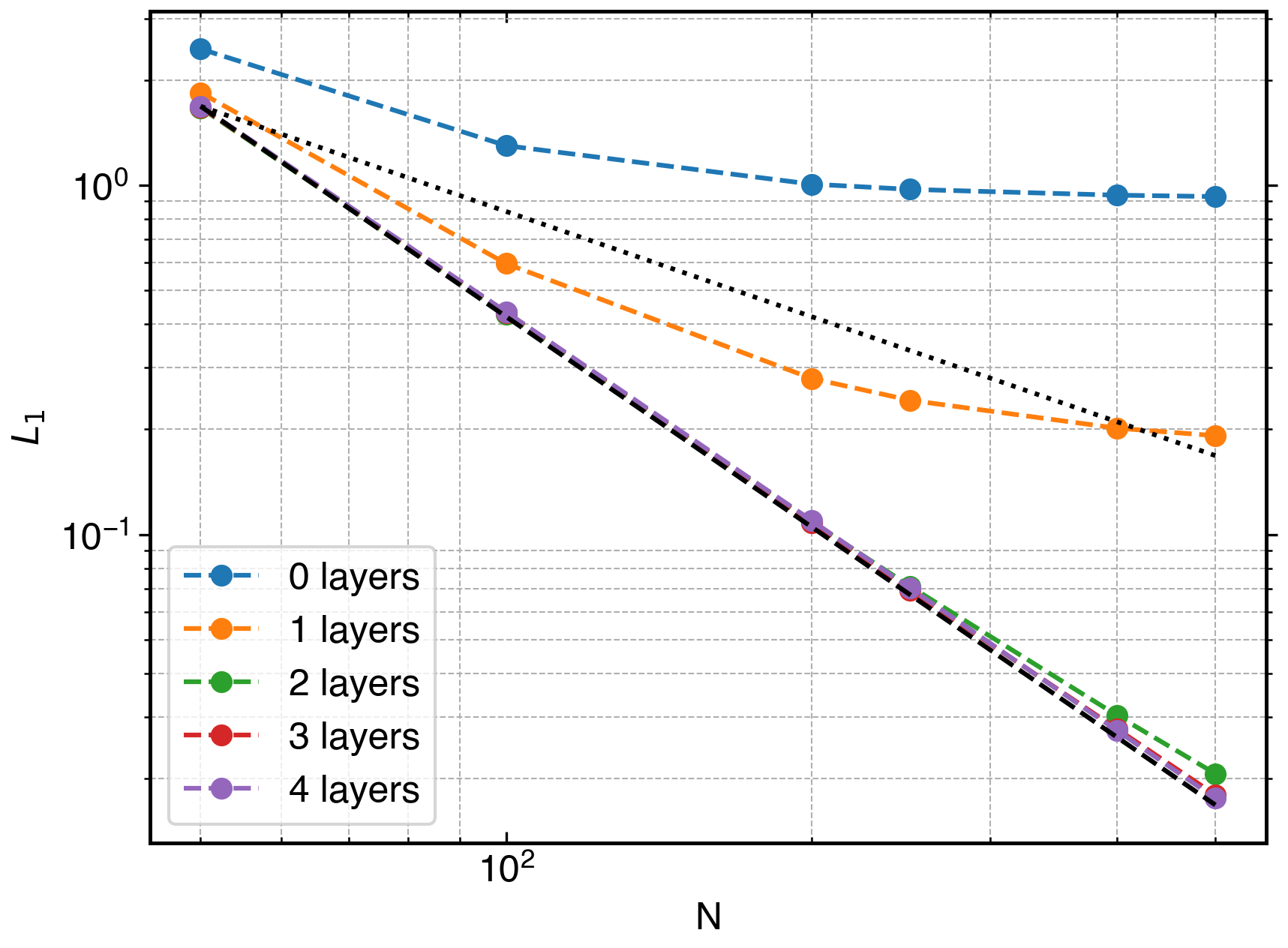}
\caption{$L_{1}$ error in pressure gradient (left) and Laplacian (right)
approximation with change in the skipped number of layers.}
  \label{fig:l1_layer}
\end{figure}

We compute the SPH approximation of the pressure gradient and Laplacian of
velocity on the particles using the discretization in \cref{eq:mom_desc}.
Since the particles near the boundary do not have complete support, they
will have higher errors compared to the inner particles. We compute the
$L_1$ error in the approximation using
\begin{equation}
  L_1(N_p) = \frac{\sum_j^{N_p} |f_j(\ten{x}) - \tilde{f}(\ten{x})|}{N_p},
\end{equation}
where $f(x)$ is the computed value and $\tilde{f}(\ten{x})$ is the actual
value, $N_p$ are the number of particles in the domain $(1 -  2n \Delta x)
\times (1 - 2n \Delta x)$, where $n$ is the number of layers skipped from
each side. In \cref{fig:l1_layer}, we plot the $L_1$ error for the
particles skipping $n$-layers of particles from the outermost boundary for
pressure gradient and Laplacian approximations. We note that both the
discretization used are second-order accurate in a periodic
domain~\cite{negi2021numerical}. However, we observe that the pressure
gradient is second-order accurate even when zero layers of particles are
surrounding the domain of interest. Whereas, for the Laplacian
approximation, we require at least 2 layers of particles. Since all the
second-order accurate formulations for Laplacian approximation require
gradient on all the neighboring particles \cite{negi2021numerical}
therefore the requirement for the number of neighboring particles is double
to that in the case of the gradient approximation for the second-order
accurate convergence.

This numerical test shows that the boundary implementation which uses a
single layer of particles is bound to have error in Laplacian approximation
resulting in an inaccurate solution irrespective of having an accurate
boundary condition implementation. Furthermore, since all the second order
viscosity formulations require the evaluation of gradient on the boundary
particle~\cite{negi2021numerical}, we do not pursue the implementation of
boundaries based on the local point symmetry method in
\cite{ferrari2009,fourtakas2015,fourtakas2019} for the no-slip boundary
condition. The various solid boundary implementations considered are as
follows:

\subsection{Using a single layer of ghost particles on the boundary surface}

\subsubsection{With virtual particles}

\begin{figure}
  \centering
  \begin{tikzpicture}
    \fluids
    \fill[red!60](0.5,0) circle (1ex);
    \fill[red!60](0.5,1) circle (1ex);
    \fill[red!60](0.5,2) circle (1ex);
    \fill[red!60](0.5,3) circle (1ex);
    \fill[red!60](0.5,4) circle (1ex);
    \fill[red!60](0.5,5) circle (1ex);
    \fill[black](0.5,3) circle (0.2ex);
    \draw[->] (0.5,3) -- node[anchor=center, yshift=0.3cm] {kh} (2.0,4) ;
    \draw[dashed] (0.5, 3) circle (10ex);

    \fill[black](0.5,3) circle (0.5ex);
    \fill[black](-0.5,3) circle (0.5ex);
    \fill[black](-1.5,3) circle (0.5ex);
    \fill[black](-2.5,3) circle (0.5ex);
    \fill[black](-3.5,3) circle (0.5ex);
  \end{tikzpicture}
\caption{Arrangement of ghost and virtual particles for given fluid
particle denoted by blue circles. The dashed line is the boundary surface,
red particles are ghost (solid) particles. The black dots are generated to
evaluate gradient using finite difference at the red boundary points.}
  \label{fig:bc_mo}
\end{figure}
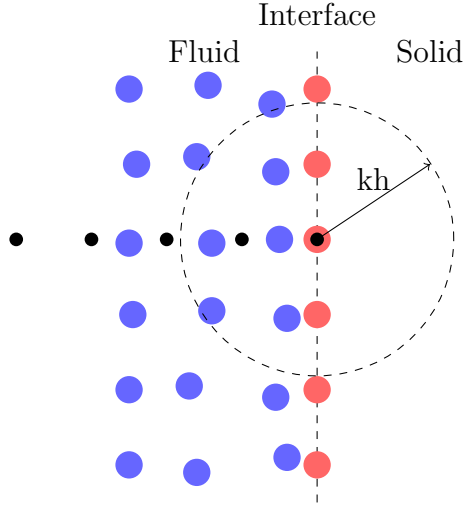

In these methods, only one layer of particles are used on the boundary
surface. \citet{marongiu2007} proposed a characteristics-based evolution
equation for pressure update at these boundary particles are given by
\begin{equation}
  \frac{d \varrho}{ d t} = c_o \frac{\partial \varrho}{ \partial n} -
  \varrho \frac{\partial u_n}{ \partial n} - \frac{\varrho \ten{g} \cdot
  \ten{n}}{c_o}
  \label{eq:char_evol}
\end{equation}
where $\ten{n}$ is the normal of the boundary surface pointing into the
fluid. In order to evaluate the gradient at the boundary point, five-point
finite difference approximation is used. These five points are generated
along the normal of the boundary particle at a spacing equal to the average
particle spacing represented by black point in \cref{fig:bc_mo} for a
single particle. The values of the properties at these black points are
evaluated using Shepard interpolation~\cite{shepard_1968}.

\subsubsection{Without virtual particles}

\citet{hashemi2012} proposed methods to implement pressure and
no-slip boundary conditions using one layer of boundary particles on the
surface. However, they do not use any extra set of elements to derive the
values on the boundary elements denoted by red particles in the
\cref{fig:bc_mo}. In order to satisfy the no-slip boundary condition the
velocity of the red particles is kept the same as the velocity of the solid
boundary. Whereas, for pressure boundary implementation, momentum balance
is performed along the normal of the boundary surface given by
\begin{equation}
  \frac{\nabla p}{\varrho} \cdot \ten{n} = -\frac{d \ten{u}}{d t} \cdot
  \ten{n} + \frac{ \mu \nabla^2 \ten{u}}{\varrho} \cdot \ten{n} + \ten{g} \cdot
  \ten{n},
  \label{eq:hash}
\end{equation}
where $\ten{n}$ is the normal of the boundary surface pointing into the
fluid. The \cref{eq:hash} is discretized using the second-order consistent
approximation to obtain the pressure at $i$\textsuperscript{th} solid
boundary particles are given by
\begin{equation}
  p_i = \frac{\left( \frac{p_j}{\varrho_i} \tilde{\nabla} W_{ij} \omega_j
  \right) \cdot \ten{n}_i - \left< -\frac{d \ten{u}}{d t} \cdot \ten{n} +
  \frac{ \mu \nabla^2 \ten{u}}{\varrho} \cdot \ten{n} + \ten{g} \cdot
  \ten{n} \right>_i}{ \left( \frac{1}{\varrho_i} \tilde{\nabla} W_{ij}
  \omega_j \right) \cdot \ten{n}_i}.
  \label{eq:hash_desc}
\end{equation}

\subsection{Using multiple layers of ghost particles outside boundary surface}

\subsubsection{With virtual particles}

\begin{figure}
  \centering
  \begin{tikzpicture}
    \fluids
    \solids
    \draw (0, 0) pic[rotate = 30, thick, red] {cross=5pt};
    \draw (0,1) pic[rotate = 30, thick, red] {cross=5pt};
    \draw (0,2) pic[rotate = 30, thick, red] {cross=5pt};
    \draw (0,3) pic[rotate = 30, thick, red] {cross=5pt};
    \draw (0,4.0) pic[rotate = 30, thick, red] {cross=5pt};
    \draw (0,5.0) pic[rotate = 30, thick, red] {cross=5pt};
    \draw (-1,0) pic[rotate = 30, thick, red] {cross=5pt};
    \draw (-1,1.) pic[rotate = 30, thick, red] {cross=5pt};
    \draw (-1,2.0) pic[rotate = 30, thick, red] {cross=5pt};
    \draw (-1,3) pic[rotate = 30, thick, red] {cross=5pt};
    \draw (-1,4.0) pic[rotate = 30, thick, red] {cross=5pt};
    \draw (-1,5.0) pic[rotate = 30, thick, red] {cross=5pt};
    \draw (-2,0) pic[rotate = 30, thick, red] {cross=5pt};
    \draw (-2,1.0) pic[rotate = 30, thick, red] {cross=5pt};
    \draw (-2,2.0) pic[rotate = 30, thick, red] {cross=5pt};
    \draw (-2,3) pic[rotate = 30, thick, red] {cross=5pt};
    \draw (-2,4.0) pic[rotate = 30, thick, red] {cross=5pt};
    \draw (-2,5.0) pic[rotate = 30, thick, red] {cross=5pt};
  \end{tikzpicture}
\caption{Arrangement of ghost and virtual particle for given blue fluid
particles and solid boundary shown by dashed line. The red circle on the
right are the ghost particles, and the red crosses are created by
reflecting the ghost about the interface.}
  \label{fig:marr}
\end{figure}
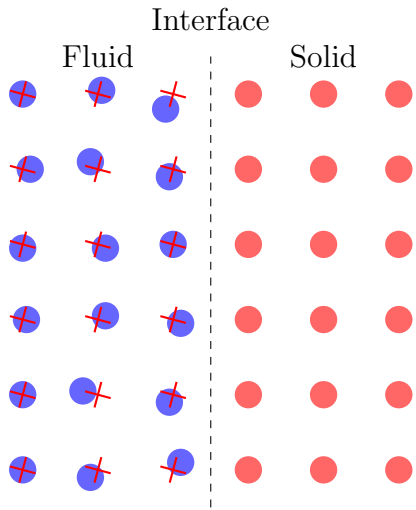

\citet{marrone-deltasph:cmame:2011} proposed a method where fixed virtual
particles are generated by reflecting the ghost particles about the
interface. The created particles are illustrated in the \cref{fig:marr},
where red crosses are the virtual particles and red particles on the right
represent the ghost particles. The properties on the ghost particles are
set as $\varrho_g=\varrho_v$, $p_g = p_v$, and the velocity is set
according to the slip or no-slip condition required, where $*_g$ represent
the property value on ghost particle and $*_v$ represent the property value
on the corresponding virtual particle. In case of slip, the velocity normal
to the wall is reversed, whereas in case of no-slip, the velocity is set
negative of the value on the corresponding virtual particle. The properties
on the virtual particles are evaluated using
\begin{equation}
  f_i = \sum_j f_j \tilde{W}_{ij} \omega_j,
\end{equation}
where $f_*$ is the desired property, $\omega_j$ is the truncated volume of
the fluid particles and $\tilde{W}_{ij}$ is the kernel corrected using the
method by \citet{liu_restoring_2006}. The sum $j$ is taken over all the
fluid particles in the support of the kernel,

\subsubsection{Without virtual particles}

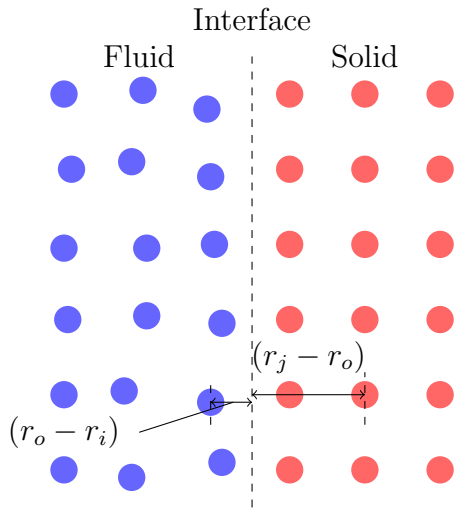
\begin{figure}
  \centering
  \begin{tikzpicture}
    \fluids
    \solids
    \draw[dashed] (-0.05,0.6) -- (-0.05,1.2) ;
    \draw[<->] (-0.05,0.9) -- (0.5,0.9);
    \draw[-] node[xshift=-2cm , yshift=0.5cm]{$(r_o-r_i)$} (-1, 0.5) -- (0.25, 0.9);
    \draw[dashed] (2,0.6) -- (2,1.3) ;
    \draw[<->] (0.5,1) -- node[yshift=0.5cm]{$(r_j-r_o)$} (2,1);
  \end{tikzpicture}
\caption{Arrangement of particles for the method proposed by
\citet{takeda1994a}. The ghost particles are shown in red, and the fluid are
shown in blue color particles. The value on the red particles are
interpolated using the nearest fluid particle along the normal of the
boundary surface.}
  \label{fig:takeda}
\end{figure}

In SPH literature, most of the methods for boundary condition
implementation use multiple layers of ghost particle such that the kernel
has full support. \citet{takeda1994a} proposed to extrapolate the
properties using a linear interpolation depending upon the distance of the
ghost particle from the nearest fluid particle as shown in
\cref{fig:takeda}. The extrapolated velocity is given by
\begin{equation}
  \ten{u}_j = -\ten{u}_i \frac{r_j - r_o}{r_o - r_i},
  \label{eq:takeda}
\end{equation}
where $r_j-r_o$ and $r_o-r_i$ are the distances along the normal from the
boundary for $j$\textsuperscript{th} ghost and $i$\textsuperscript{th}
fluid particle, respectively. We select the nearest fluid particle for a
particular ghost particle along the normal of the boundary.

\begin{figure}
  \centering
  \begin{tikzpicture}
    \fluids
    \solids
    \fill[blue!30](0.1,0.1) circle (1ex);
    \fill[blue!30](-0.05,0.9) circle (1ex);
    \fill[blue!30](0.1,1.95) circle (1ex);
    \fill[blue!30](0,3.0) circle (1ex);
    \fill[blue!30](-0.05,3.9) circle (1ex);
    \fill[blue!30](-0.1,4.8) circle (1ex);

    \fill[black](0.,3) circle (0.2ex);
    \draw[->] (0.,3) -- node[anchor=center, yshift=0.3cm] {kh} (1.5,4) ;
    \draw[dashed] (0.0, 3) circle (10ex);

    \fill[black](1.,2) circle (0.2ex);
    \draw[->] (1.,2) -- node[anchor=center, yshift=0.3cm] {kh} (2.5,3) ;
    \draw[dashed, red] (1, 2) circle (10ex);
  \end{tikzpicture}
\caption{Arrangement of particles for the method proposed by
\citet{randles1996smoothed}, and \citet{Adami2012} for the fluid particles
in blue and boundary interface shown by dashed line. The red circles on the
left represents the ghost solid particles. The property of light blue fluid
particle is manipulated to satisfy boundary condition.}
  \label{fig:rand}
\end{figure}
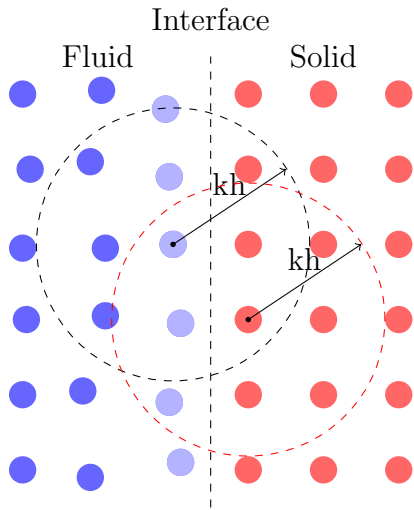

\citet{randles1996smoothed} proposed a boundary condition for pressure. The
prescribed value is assigned to the ghost particles in red and the values
on the light blue fluid particles as shown in \cref{fig:rand} are set such
that the desired condition is satisfied. The value on light blue fluid
particles is given by
\begin{equation}
  p_i = p_{bc} + \frac{\sum_{j \in I}  (p_j-p_{bc}) W_{ij} \omega_j}{1 - \sum_{j \in
  B} W_{ij} \omega_j},
\end{equation}
where $I$ is the set of blue particles, $B$ is the set of red particle, and
pressure is evaluated at each $i$\textsuperscript{th} light blue particle
shown in \cref{fig:rand}.

\citet{Adami2012} proposed the method where Shepard interpolation is
employed to extrapolate properties from fluid to ghost particles. The
property at a ghost particle is given by
\begin{equation}
  f_i = \frac{\sum_j f_j W_{ij}}{\sum_j W_{ij}},
\end{equation}
where sum $j$ is over all the fluid particles in the support of the kernel
as shown by red dashed circle in \cref{fig:rand}, at the
$i$\textsuperscript{th} ghost particle. Furthermore,
\citet{esmailisikarudi2016} proposed to perform a first order accurate
extrapolation to evaluate properties on ghost particles.

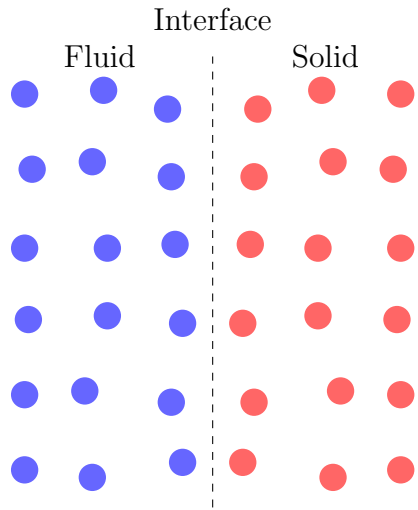
\begin{figure}
  \centering
  \begin{tikzpicture}
    \fluids

    \fill[red!60](0.9,0.1) circle (1ex);
    \fill[red!60](1.05,0.9) circle (1ex);
    \fill[red!60](.9,1.95) circle (1ex);
    \fill[red!60](1,3.0) circle (1ex);
    \fill[red!60](1.05,3.9) circle (1ex);
    \fill[red!60](1.1,4.8) circle (1ex);
    \fill[red!60](2.1,-0.1) circle (1ex);
    \fill[red!60](2.2,1.05) circle (1ex);
    \fill[red!60](1.9,2.05) circle (1ex);
    \fill[red!60](1.9,2.95) circle (1ex);
    \fill[red!60](2.1,4.1) circle (1ex);
    \fill[red!60](1.95,5.05) circle (1ex);
    \fill[red!60](3,0) circle (1ex);
    \fill[red!60](3,1.0) circle (1ex);
    \fill[red!60](2.95,2.0) circle (1ex);
    \fill[red!60](3,2.95) circle (1ex);
    \fill[red!60](2.9,4.0) circle (1ex);
    \fill[red!60](3,5.0) circle (1ex);
  \end{tikzpicture}
\caption{Arrangement of particles for the method proposed by
\citet{colagross2003a} for the fluid particles in blue, and boundary
represent by the dashed line. The ghost particles in red representing solid are
created by the mirror reflection about the interface.}
\label{fig:cola}
\end{figure}

\citet{colagross2003a} proposed to mirror the fluid particles near the
solid interface, about the interface to generate ghost particles as shown
in \cref{fig:cola}. These ghost particles carry the velocity of opposite
sign to implement no penetration. The value of pressure and density is kept
the same. In order to implement this method, we create new particles for
solids from the fluid particles before each time step.

\section{Open boundary conditions}
\label{sec:obc}

Open boundary conditions are required to simulate a wind-tunnel kind of
simulation in SPH. It consist of an inlet from where the particles are
added to the domain, and an outlet from where the particles exits the
domain as shown in \cref{fig:io}. The particles added to the domain should
have prescribed inlet velocity and should not introduce any artifacts in
the flow whereas the outlet should remove particles from the flow without
affecting the flow. Since we use a weakly-compressible scheme, pressure
waves are inevitable. Therefore, the inlet/outlet must have non-reflecting
property along with adherence to the boundary condition. In \cref{fig:io},
we show a schematic arrangement of the inlet, fluid, and outlet domains.
The dashed line in between different regions are the boundary interfaces.
The inlet/fluid particles are converted to fluid/outlet particles once they
cross the corresponding interface. Many
authors~\cite{federico2012,lastiwka_nrbc_2009,tafuni_io_2018,negi2021numerical}
have proposed different methods to implement inlet/outlet boundary
condition in SPH. Open boundary condition implementations considered in
this paper are as follows:

\begin{figure}[htbp]
  \centering
  \begin{tikzpicture}
    \io
  \end{tikzpicture}
  \caption{Schematic of the arrangement of fluid with the inlet and outlet
  particles. The top and bottom are supported by solid particles not shown
  in the figure.}
  \label{fig:io}
\end{figure}
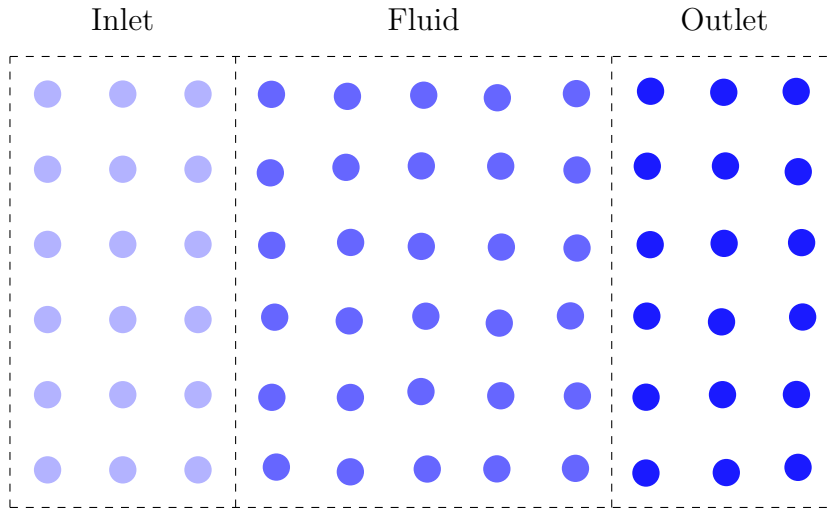

\subsection{Do-nothing}

\citet{federico2012} proposed the method to implement outlet boundary where
the properties of the fluid particles crossing the fluid-outlet interface
are fixed in time. Therefore, the particle, after crossing the interface,
advects with a constant velocity.

\subsection{Mirror}

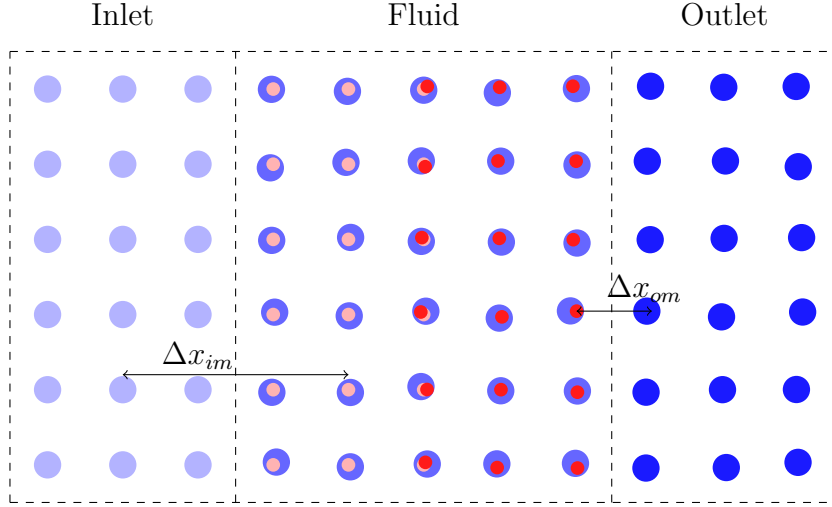
\begin{figure}[htbp]
  \centering
  \begin{tikzpicture}
    \io
    \fill[red!30](-4,0) circle (0.5ex);
    \fill[red!30](-4,1.0) circle (0.5ex);
    \fill[red!30](-4,2.0) circle (0.5ex);
    \fill[red!30](-4,3) circle (0.5ex);
    \fill[red!30](-4,4.0) circle (0.5ex);
    \fill[red!30](-4,5.0) circle (0.5ex);
    \fill[red!30](-3,0) circle (0.5ex);
    \fill[red!30](-3,1.0) circle (0.5ex);
    \fill[red!30](-3,2.0) circle (0.5ex);
    \fill[red!30](-3,3) circle (0.5ex);
    \fill[red!30](-3,4.0) circle (0.5ex);
    \fill[red!30](-3,5.0) circle (0.5ex);
    \fill[red!30](-2,0) circle (0.5ex);
    \fill[red!30](-2,1.0) circle (0.5ex);
    \fill[red!30](-2,2.0) circle (0.5ex);
    \fill[red!30](-2,3) circle (0.5ex);
    \fill[red!30](-2,4.0) circle (0.5ex);
    \fill[red!30](-2,5.0) circle (0.5ex);

    \fill[red!90](0.0450781577,-0.04221) circle(0.5ex);
    \fill[red!90](0.0436326628,0.966944) circle(0.5ex);
    \fill[red!90](0.0331861466,2.045488) circle(0.5ex);
    \fill[red!90](-0.011271014,2.9962833) circle(0.5ex);
    \fill[red!90](0.0270670158,4.040094) circle(0.5ex);
    \fill[red!90](-0.015977573,5.037394) circle(0.5ex);
    \fill[red!90](-1.023183498,-0.035654) circle(0.5ex);
    \fill[red!90](-0.972560509,1.0005706) circle(0.5ex);
    \fill[red!90](-0.959527978,1.9692822) circle(0.5ex);
    \fill[red!90](-0.994235407,3.0133835) circle(0.5ex);
    \fill[red!90](-1.011956196,4.0432622) circle(0.5ex);
    \fill[red!90](-0.989535452,5.0236606) circle(0.5ex);
    \fill[red!90](-1.97682717,0.03333736) circle(0.5ex);
    \fill[red!90](-1.95544884,1.00277262) circle(0.5ex);
    \fill[red!90](-2.038032744,2.0337744) circle(0.5ex);
    \fill[red!90](-2.024208056,3.0235824) circle(0.5ex);
    \fill[red!90](-1.97974001,3.96700024) circle(0.5ex);
    \fill[red!90](-1.952165244,5.0348294) circle(0.5ex);

\draw[<->] (-3,1.2) -- (-6,1.2) node[xshift=1cm, yshift=0.25cm]{$\Delta
x_{im}$};
    \draw[<->] (0.0331861466,2.045488) -- (1.0331861466,2.045488) node[xshift=-0.1cm, yshift=0.3cm]{$\Delta
    x_{om}$};

  \end{tikzpicture}
  \caption{Schematic of the arrangement of fluid with the inlet and outlet
particles and their corresponding virtual particles are shown in light red and
red, respectively.}
  \label{fig:mirror}
\end{figure}

\citet{tafuni2018} proposed a method applicable to both inlets and outlets.
The properties from the fluid are stored on the virtual particles generated
by reflecting the inlet/outlet particles about the interface as shown in
\cref{fig:mirror}. The desired property and its gradient are evaluated for
each mirror particle using a first-order consistent approximation. From
the mirror particle, the property of the corresponding outlet/inlet
particle is evaluated using
\begin{equation}
  f_o = f_i - \Delta x_{om} \nabla f_i,
  \label{eq:mirr_tay_exp}
\end{equation}
where $f_i$ and $\nabla f_i$ are the properties on the mirror particle and
$x_{om}$ is the distance between the outlet and mirror particle. In
addition to this method, we test the convergence of the method when the
Taylor series correction is not applied such that the
\cref{eq:mirr_tay_exp} is simplified to
\begin{equation}
  f_o = f_i.
  \label{eq:mirr_hyd}
\end{equation}
We refer to this method as simple-mirror.

\subsection{Characteristics}

\citet{lastiwka_nrbc_2009} proposed a non-reflecting inlet and outlet
boundary implementation using the method of characteristics. The
characteristics of the flow are given by
\begin{equation}
  \begin{split}
    J_1 &= -c_o^2 (\rho - \rho_{ref}) + (p - p_{ref}),\\
    J_2 &= \rho c_o (u - u_{ref}) + (p - p_{ref}),\\
    J_3 &= -\rho c_o (u - u_{ref}) + (p - p_{ref})i,\\
  \end{split}
  \label{eq:charac}
\end{equation}
where $\rho_{ref}$, $u_{ref}$, and $p_{ref}$ are the reference density,
velocity, and pressure. These characteristics are extrapolated to the
inlet/outlet particles, and the properties are calculated using the
extrapolated characteristics are given by
\begin{equation}
  \begin{split}
  \rho &=  \rho_{ref} + \frac{1}{c_o^2} \left(-J_1 + \frac{1}{2}J_2 + \frac{1}{2} J_2 \right),\\
  u & = u_{ref} + \frac{1}{2 \rho c_o } \left(J_2 - J_3\right),\\
  p & = p_{ref} + \frac{1}{2} \left(J_2 + J_3\right).\\
  \end{split}
\end{equation}
In order to implement the inlet boundary $J_1$ and $J_2$ should be set zero
and $J_3$ must be extrapolated from the fluid to the inlet particles.
Whereas, for outlet $J_1$ and $J_2$ are extrapolated from the fluid to the
outlet, and $J_3$ is set to zero. Recently, \citet{negi_nrbc_2020},
proposed to obtain the reference value of the properties by taking a time
average of properties when the acoustic intensity of the particle is below
some prescribed value. This reference value is determined based on the
inlet velocity. We refer to this method as hybrid method.

In the next section, we discuss the procedure one should follow for the
construction of manufactured solution in the method of manufactured
solutions to test boundary condition implementation in the context of SPH.

\section{Method of manufactured solutions}
\label{sec:mms_bc}

It is of utmost importance that the discretization and approximations
employed in a solver are of the desired accuracy. Both verification and
validation play an essential role to determined the correctness of a
solver. In validation, we test whether the intended physics modeled by the
differential equation converges as required. Whereas, in verification, we
determine whether the discretization of various terms in the governing
equation are coded correctly to reflect the intended rate of convergence.
The MMS is a verification technique widely used for the verification of
finite volume
method~\cite{waltzManufacturedSolutionsThreedimensional2014,choudharyCodeVerificationBoundary2016}
and finite element method~\cite{gfrererCodeVerificationExamples2018} codes.
The MMS can also be used to determine specific terms causing a lower order
of convergence or if the code has a mistake. Recently,
\citet{negiHowTrainYour2021a} proposed procedures to apply MMS in the
context of SPH. They show various methods to find coding mistakes as well
as procedures to obtain manufactured solution (MS) to verify boundary
condition implementation at straight boundaries.

In the MMS, we substitute an MS to the governing partial differential
equations to obtain a residue. For example, for an MS of the form
$\tilde{\ten{u}}(x, y, z, t)$, and governing equation $\mathscr{L}\ten{u} =
g$, where $\mathscr{L}$ is an operator and $g$ is a source term, we obtain
\begin{equation}
  \mathscr{R} = \mathscr{L}\tilde{\ten{u}} - g
\end{equation}
as the residue. We use this residue as a source term in the solver to
obtain the MS from the solver. In SPH, the source term can be easily
applied to a particle by just substituting the value of the coordinates in
the source function. We solve the equations using the solver on a domain
with fluid and solid particles. Since we use a second-order convergent
scheme to solve the fluid flow, the errors are generated due to the
boundary condition implementations only. Therefore, the error at various
resolutions shows the rate of convergence of the boundary condition
implementation.

In this paper, we construct various MS such that it satisfies the boundary
condition on different domain shapes. We take the following steps in order
to construct an MS for the boundary $c(x, y)=0$:
\begin{enumerate}
\item Construct a base MS such that it satisfies the general requirement of
MMS \cite{negiHowTrainYour2021a} as follows.
  \begin{itemize}
\item The MS must be differentiable up to the highest order present in the
governing equations.
\item The MS must satisfy the required physics of the governing equation,
for example, if the solver requires the density to be continuous, the MS
must have continuous density.
\item It must not prevent the successful completion of the solver.
\item It must be bounded in the domain of interest.
  \end{itemize}
\item In the context of SPH, we additionally require the velocity field to
be non-solenoidal \cite{negiHowTrainYour2021a}.
\item Modify the property of interest such that the boundary condition
at $c(x, y)=0$ is satisfied. For example, for the no-slip boundary only
velocity need to be modified.
\item We note that one should ensure that the MS of the property of
interest is non-zero on the boundary before the boundary condition is
satisfied. For example, if we need $\ten{u} \cdot \ten{n}$ at the boundary
to be zero, we must make sure that $\ten{u} \neq 0$ at the boundary.
\end{enumerate}
In appendix \ref{apn:mms}, we discuss the construction of MS for all the
methods discussed in this paper. In the next section, we use these MS
(defined in appendix \ref{apn:mms}) to test various boundary condition
implementation discussed in \cref{sec:bc} and \cref{sec:obc}.

\section{Results}
\label{sec:results}

In this section, we verify the convergence of various boundary condition
implementations discussed in \cref{sec:bc} and \cref{sec:obc} using MMS. We
first show the convergence of various solid boundary condition
implementations followed by the inflow and outflow boundary
implementations. For all the test cases, we simulate 100 timesteps for
resolutions in the range $100 \times 100$ to $500 \times 500$, and evaluate
the $L_{1}$ error in the domain given by
\begin{equation}
  L_{1} = \sum_j \sum_i \frac{|(f(x_i, y_i, t_j) - f_o(x_i, y_i, t_j))|}{N \times N_t}
\end{equation}
where $f$ is the property of interest and $f_o$ is the property value
determined using the MS, $N$ is the number of particles in the domain, $N_t$
are number of instances in time\footnote{we saved 5 time solutions for a
particular testcase.}. We fix the timestep corresponding to the highest
resolutions in order to isolate the error due to boundary condition
implementation.

We use the PySPH \cite{pysph2020} framework to implement all the methods.
All the results presented in this paper are reproducible and can be easily
generated using the automation framework \textit{automan} \cite{pr:automan:2018}. In
the interest of reproducibility, we provide the entire source code at
\url{https://gitlab.com/pypr/mms_sph_bc}.

\subsection{Comparison of solid boundary condition implementations}

\begin{figure}[htbp]
  \centering
  \includegraphics[width=\linewidth]{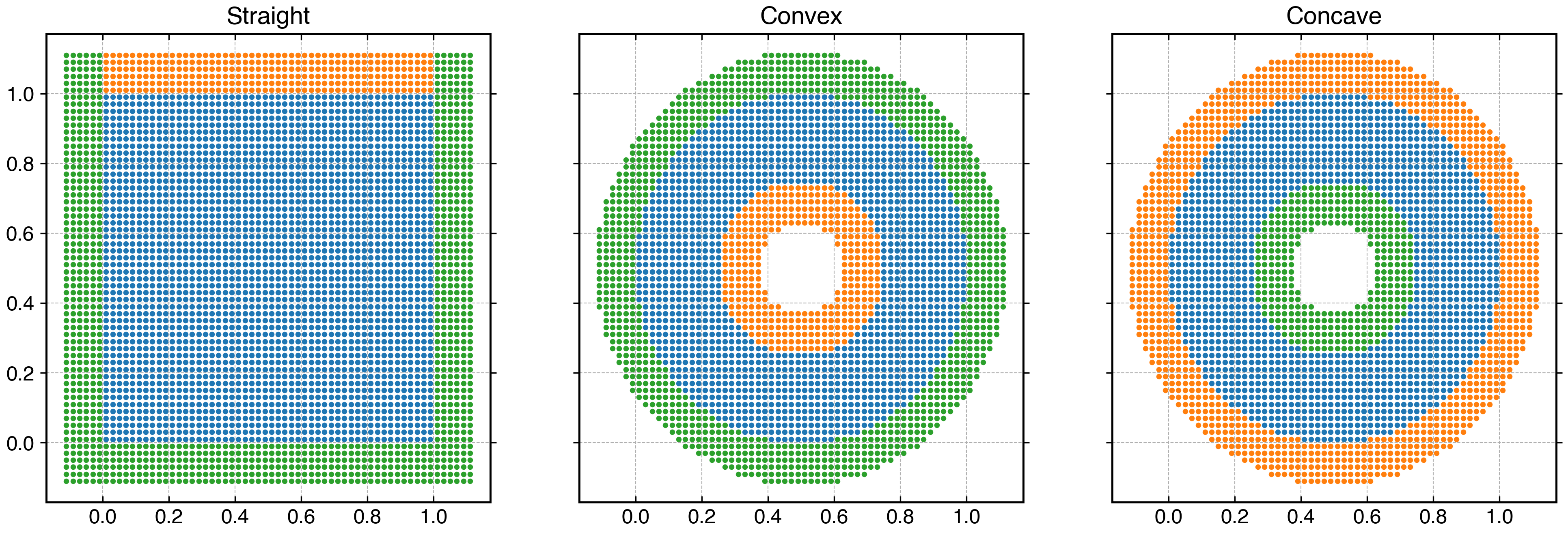}
  \caption{Types of domains considered to test the convergence of solid
  boundary implementation. The fluid particles are represented by the blue
  color, the particles in green represent ghost particles on which
  properties are set using MS, and the particles in orange are used to test
  the convergence of boundary implementation of interest.}
  \label{fig:domains}
\end{figure}

In this section, we verify all the solid boundary methods discussed in
\cref{sec:bc}. The solid boundary can be straight or curved. In this paper,
we do not consider a non-smooth geometric features, like a corner. For
corners, one requires a discontinuous MS, and at discontinuity the higher
order terms fail to show the actual order of convergence. However, the
method showing second-order convergence for smooth boundary will perform
better than other methods. We consider three types of boundary shapes viz.\
straight, convex, and concave as shown in \cref{fig:domains}. The fluid
particles are represented by the blue color, and the green color particles
represent the ghost particles for which we set the properties using the MS.
The orange colored particles are used to verify a particular method. The
domain referred to as `straight', the ghost particles in orange have a
constant normal. The domain referred to as `convex', the boundary is a
convex surface, whereas the domain referred to as `concave', the boundary
is a concave surface. We note that both convex and concave domains are
identical, however the boundary surfaces of interest are different.

\begin{figure}[htbp]
  \centering
  \includegraphics[width=0.7\linewidth]{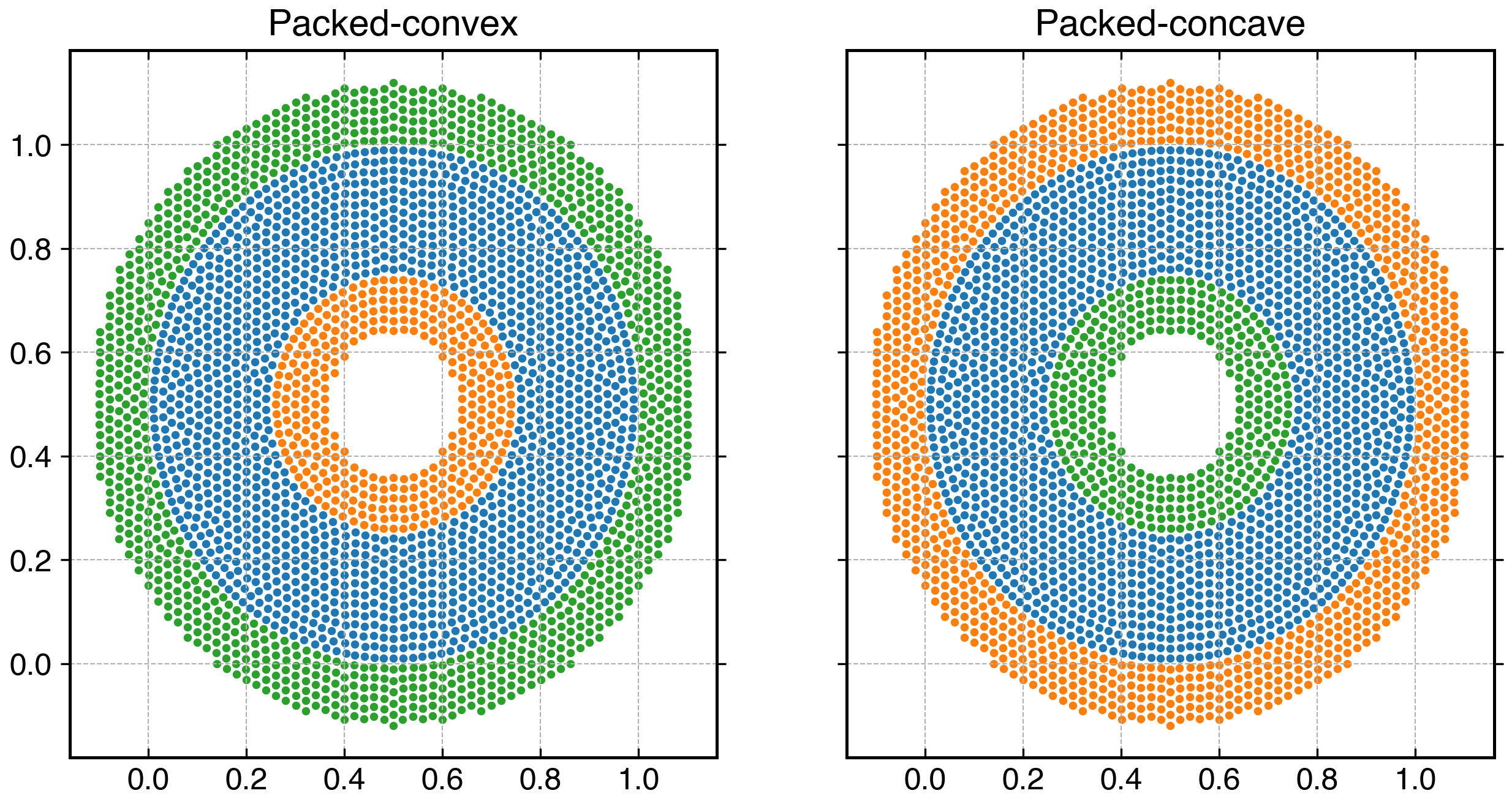}
\caption{The packed version of convex and concave domains. The fluid
particles are represented by blue color, the particles in green represents
ghost particles on which properties are set using MS, and the particles in
orange are used to test the convergence of desired boundary condition
implementation.}
  \label{fig:domains_pack}
\end{figure}

We observe that the convex and concave domains have staircase pattern at
the boundary due to the use of the Cartesian arrangement of particles to
represent the boundary. We use the method proposed by
\citet{negi2019improved}, to pack the particles near both, the inner and
outer surface. In \cref{fig:domains_pack}, we plot the packed version of
the convex and concave domains, and we refer to these as `packed-convex'
and `packed-concave' respectively.

\subsubsection{Neumann pressure boundary condition}

\begin{figure}[htbp]
  \centering
  \includegraphics[width=0.8\linewidth]{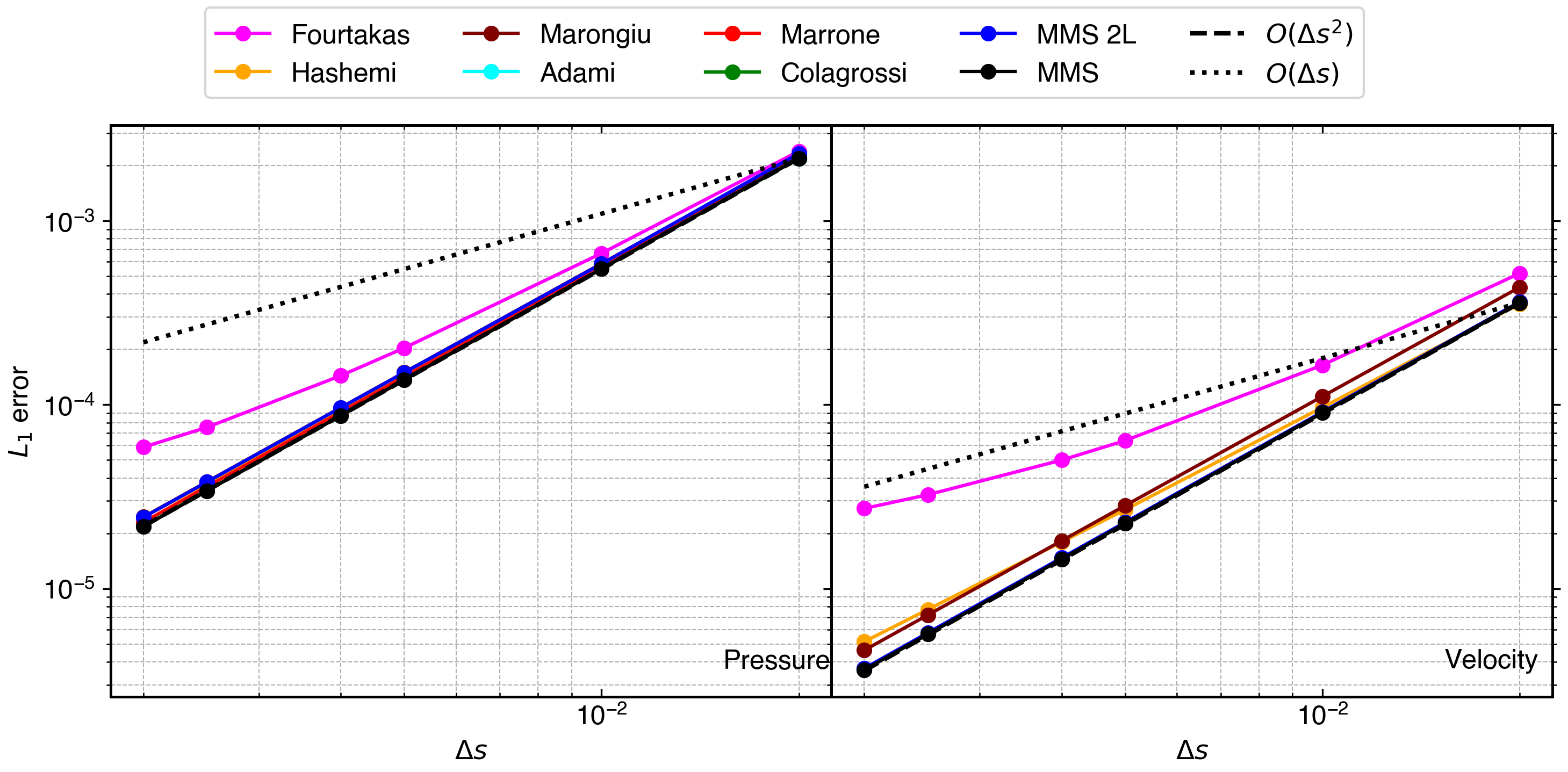}
\caption{$L_{1}$ error in pressure and velocity after 100 time steps for
different Neumann pressure boundary implementations in the straight domain
in \cref{fig:domains}.}
  \label{fig:d1_p_bc}
\end{figure}

In this section, we apply various methods described in \cref{sec:bc} to
apply the Neumann pressure boundary condition on the orange particles shown in
the domains in \cref{fig:domains}, and \cref{fig:domains_pack}. We first
use the MS in \cref{eq:mms_pres_num_d1}. The MS ensures that $\nabla p
\cdot \ten{n} = 0$ at the boundary of interest in the straight domain. We
refer to a particular method using the corresponding first author names.
The `MMS' and `MMS-2L' are the cases where properties on solid are updated
using the MS, and six layers and two layers of ghost particles are used to
represent solid, respectively. In \cref{fig:d1_p_bc}, we plot the $L_{1}$
error in the pressure and velocity after 100 timesteps.

The `MMS-2L' plot shows that the L-IPST-C scheme converges even when two
layers of solid particles are employed for a straight boundary. In the case of
the straight domain, all the methods considered in this paper are second-order
convergent except the method proposed by \citet{fourtakas2019}, where a
virtual stencil is used to complete the support of the particles. However,
the rate of convergence is similar to as reported in \cite{fourtakas2019}.

\begin{figure}[htbp]
  \centering
  \includegraphics[width=0.8\linewidth]{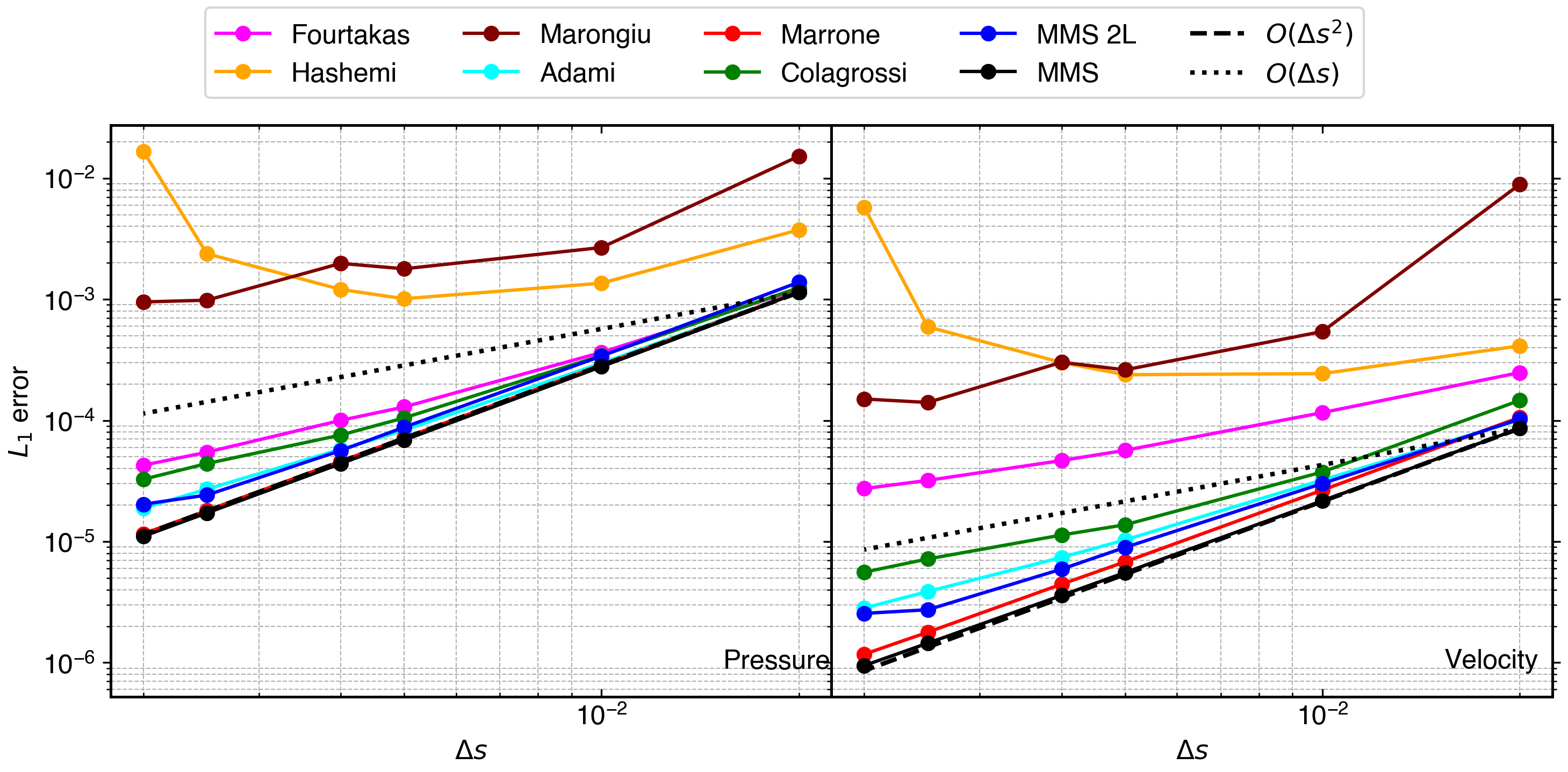}
\caption{$L_{1}$ error in pressure and velocity after 100 time steps for
different Neumann pressure boundary implementations in the convex domain
as shown in \cref{fig:domains}.}
  \label{fig:d2_p_bc}
\end{figure}

\begin{figure}[htbp]
  \centering
  \includegraphics[width=0.8\linewidth]{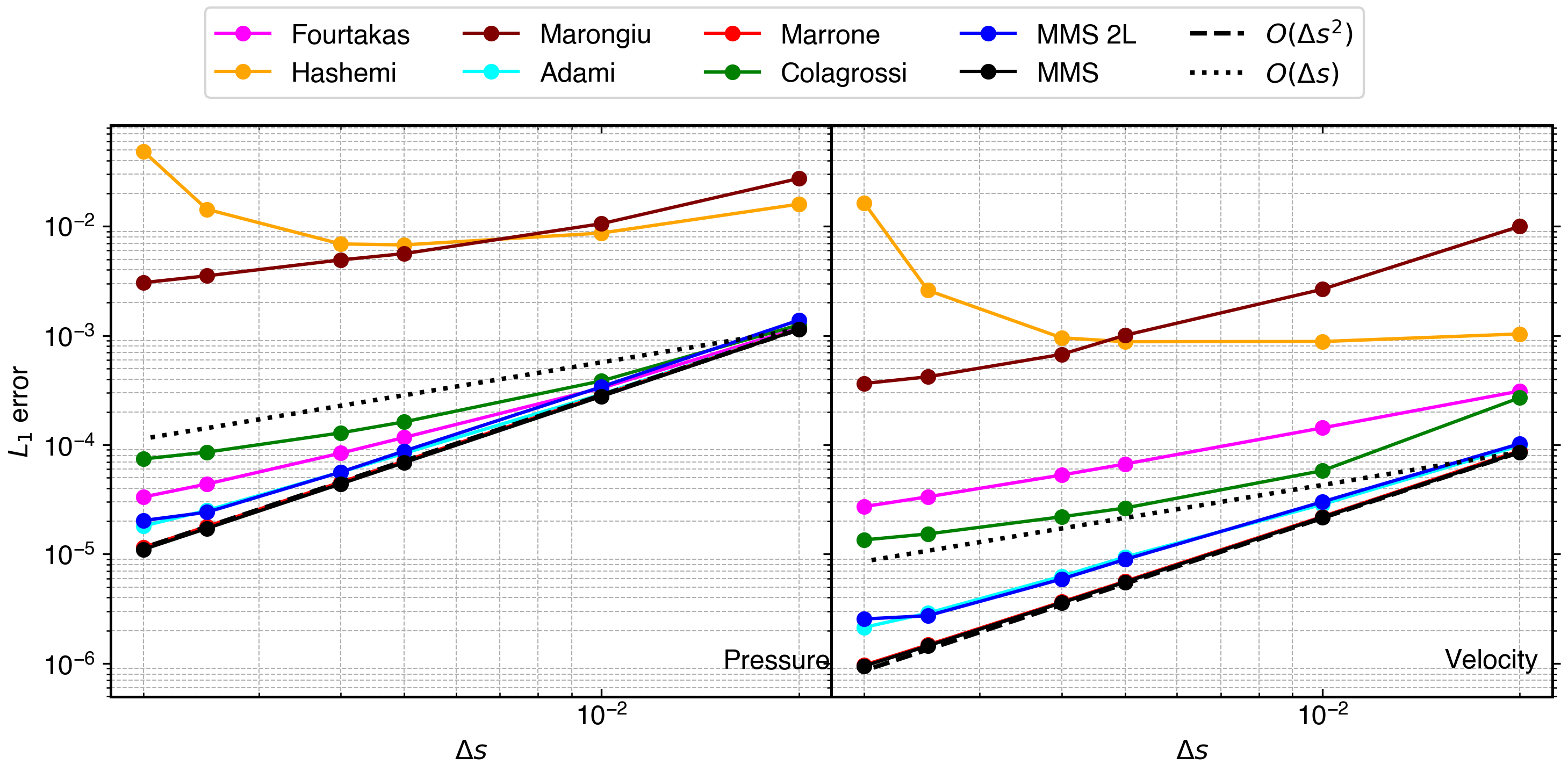}
\caption{$L_{1}$ error in pressure and velocity after 100 time steps for
different Neumann pressure boundary implementations in the concave domain
as shown in \cref{fig:domains}.}
  \label{fig:d3_p_bc}
\end{figure}

We also test all the methods on the convex and concave domains. In order to
verify, we use the MS in \cref{eq:mms_pres_num_d5} which satisfies the
boundary condition for the surface of interest in both the domains. In
\cref{fig:d2_p_bc}, and \cref{fig:d3_p_bc}, we plot the $L_{1}$ error for
pressure and velocity for the convex and concave domains, respectively.

We observe that the method proposed by \citet{hashemi_modified_2012} diverges
since only a single layer of particles are used to represent solid, which is
insufficient even for corrected gradient computation. However, the method
proposed by \citet{marongiu2007} shows first-order convergence. Moreover, the
rate of convergence is non-monotonous for a convex boundary. The reason behind
a lower order of convergence is the use of a single layer of particles and
zeroth order interpolation on the virtual particles, which are further used in
the fifth order finite difference interpolation. The method proposed by
\citet{fourtakas2019} shows first order convergence in pressure and 1.5 in
velocity as expected. The convergence of the method proposed by
\citet{Adami2012}, and `MMS-2L' are very close to 1.5. The `MMS-2L' has a
slight decrease in convergence compared to `MMS', which shows that the minimum
number of layers required for an accurate Neumann boundary is higher for a
curved surface compared to a straight boundary. Clearly, the method proposed
by \citet{marrone-deltasph:cmame:2011} is second order convergent.

\begin{figure}[htbp]
  \centering
  \includegraphics[width=0.8\linewidth]{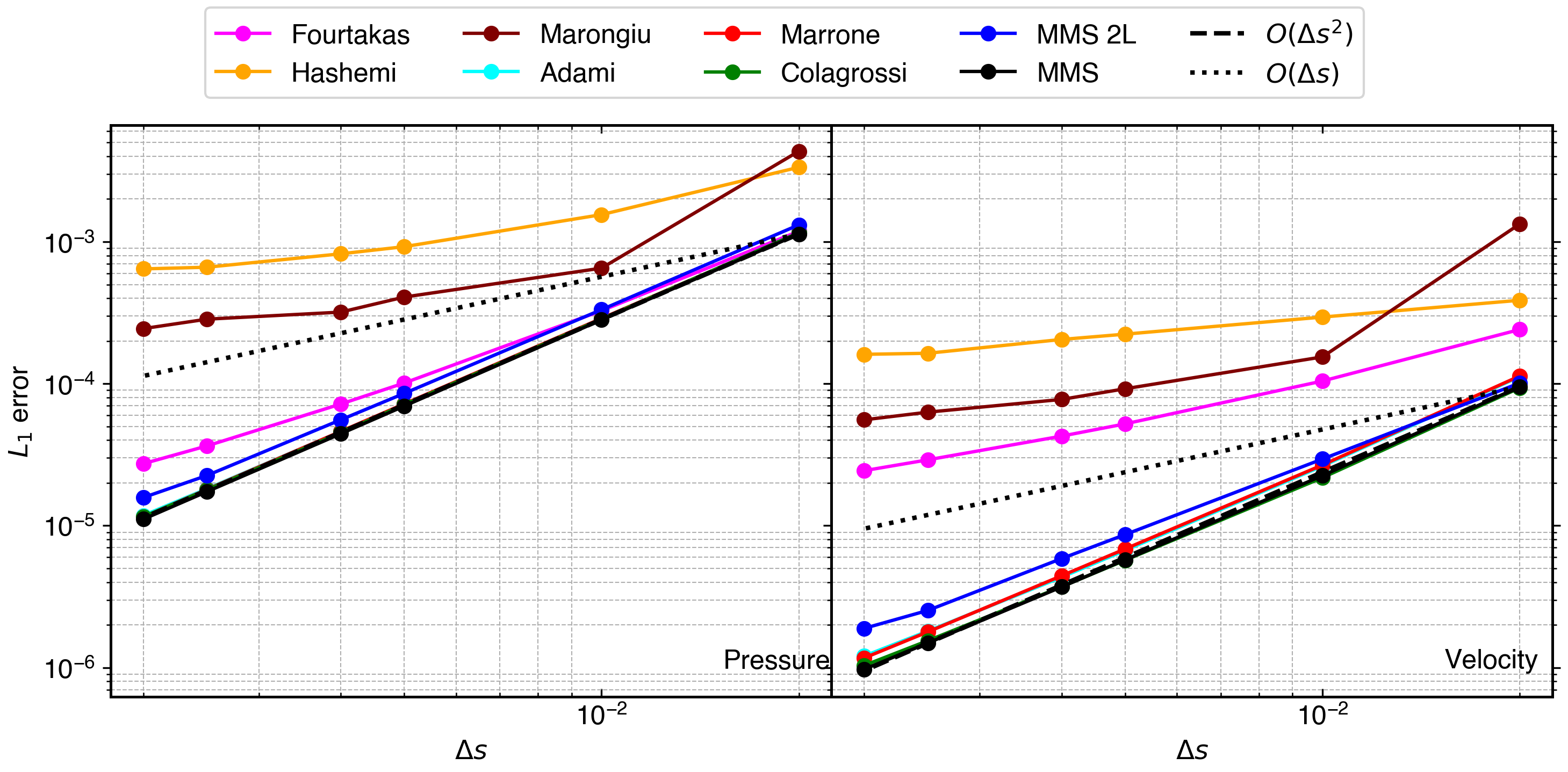}
\caption{$L_{1}$ error in pressure and velocity after 100 time steps for
different Neumann pressure boundary implementations in the packed-convex
domain as shown in \cref{fig:domains_pack}.}
  \label{fig:d4_p_bc}
\end{figure}

\begin{figure}[htbp]
  \centering
  \includegraphics[width=0.8\linewidth]{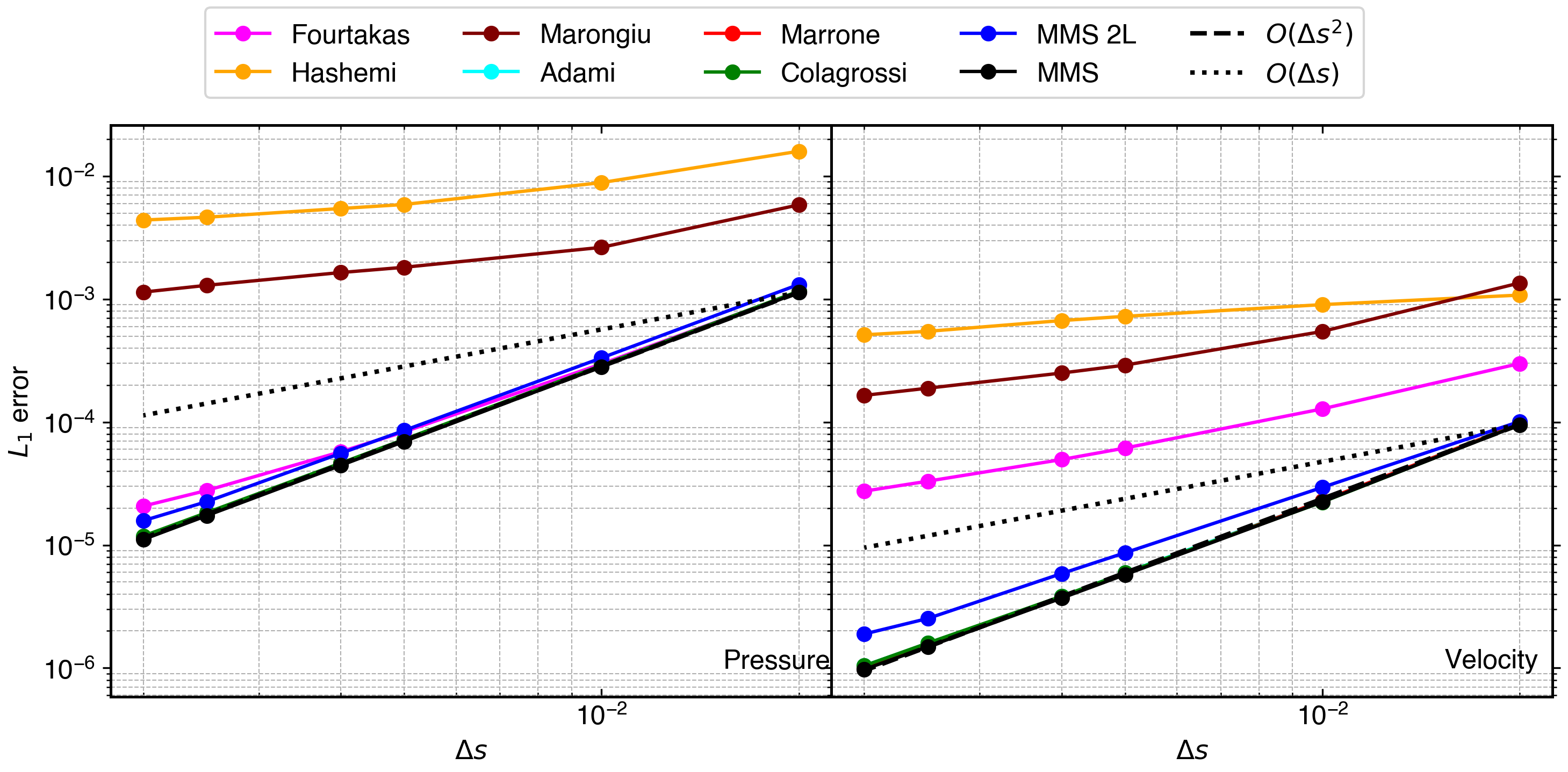}
\caption{$L_{1}$ error in pressure and velocity after 100 time steps for
different Neumann pressure boundary implementations in the packed-concave
domain as shown in \cref{fig:domains_pack}.}
  \label{fig:d5_p_bc}
\end{figure}

In order to remove the effect of jagged edges on the convergence of the
boundary condition implementations, we performed the numerical experiment on
the packed domain viz.\ packed-convex and packed concave as shown in
\cref{fig:domains_pack}. Since the boundary surfaces are the same, we use the
same MS in \cref{eq:mms_pres_num_d5}. In \cref{fig:d4_p_bc}, and
\cref{fig:d5_p_bc}, we plot the $L_1$ error for pressure and velocity for
both the domains.

We observe that all the methods show a better rate of convergence. We note
that unlike earlier in the Hashemi method errors do not increase. Furthermore,
Marongiu method shows an almost constant rate of convergence for a convex
domain. The rate of convergence increases for all the methods compared to an
unpacked domain. The convergence of the method proposed by
\citet{colagross2003a} increases by a large amount since the particles after
mirroring have good distribution. The method by
\citet{marrone-deltasph:cmame:2011} and \citet{colagross2003a} overlaps and
shows second-order convergence. This test also demonstrates the effectiveness
of packing for curved surfaces.

\subsubsection{Slip boundary condition}

\begin{figure}[htbp]
  \centering
  \includegraphics[width=0.8\linewidth]{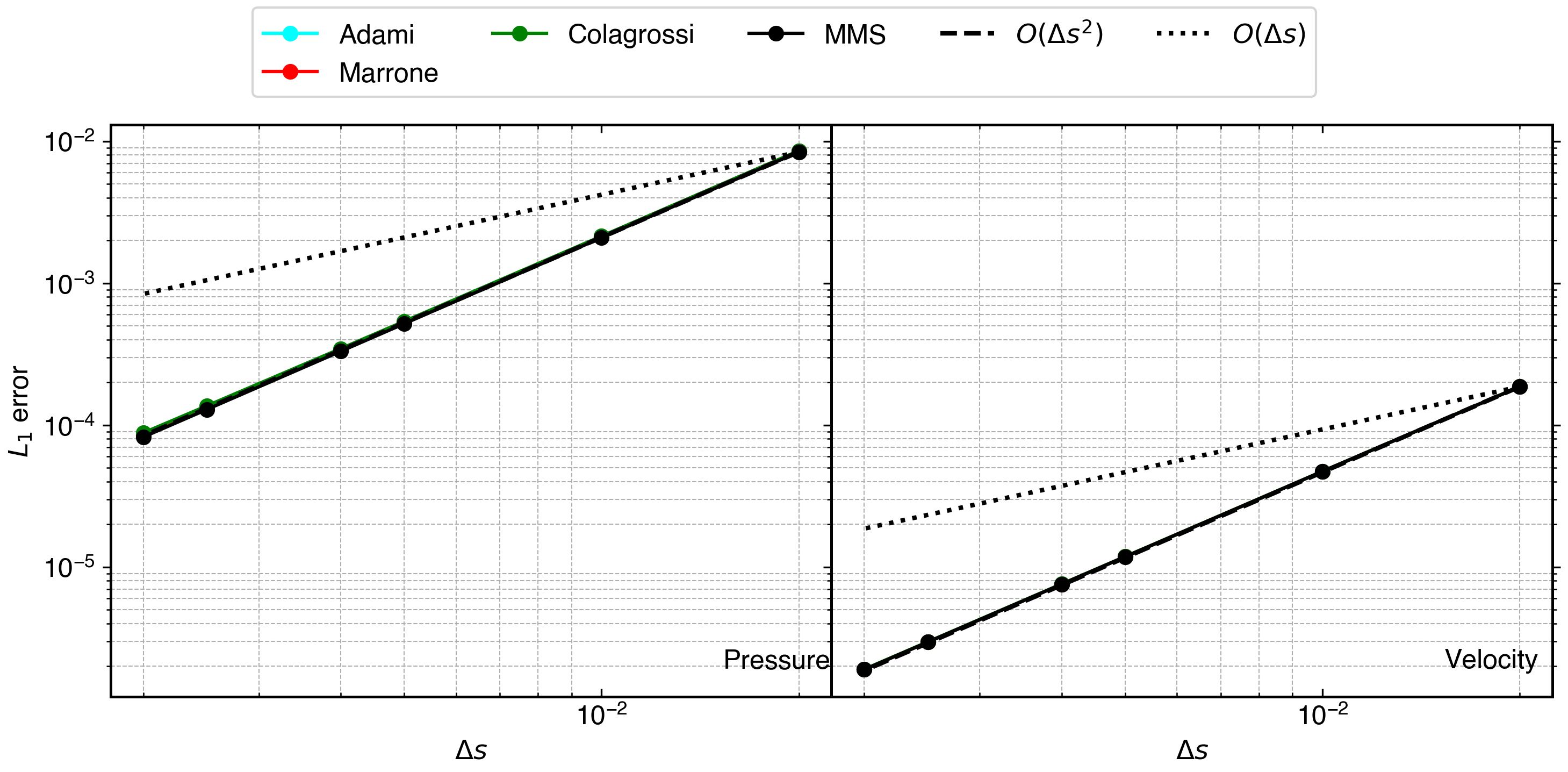}
\caption{$L_{1}$ error in pressure and velocity after 100 time steps for
different slip boundary implementations in straight domain as shown in
\cref{fig:domains}.}
  \label{fig:d1_us_bc}
\end{figure}

In this section, we test various slip boundary condition implementations
discussed in \cref{sec:bc}. In order to test these methods, we use all the
different domains considered in the previous results. For the straight
domain, we use the MS in \cref{eq:mms_slip_d1}. In order to construct this
MS we ensure that $\ten{u} \cdot \ten{n} =0$ at the boundary. In
\cref{fig:d1_us_bc}, we plot the $L_{1}$ error in pressure and velocity
after 100 timesteps. Clearly, all the methods show second-order
convergence. In general, the slip boundary condition is not a realistic
boundary condition, and it is usually used to remove the effect of walls
not affecting the flow. However, to complete the discussion, we test these
methods in other domains.

\begin{figure}[htbp]
  \centering
  \includegraphics[width=0.8\linewidth]{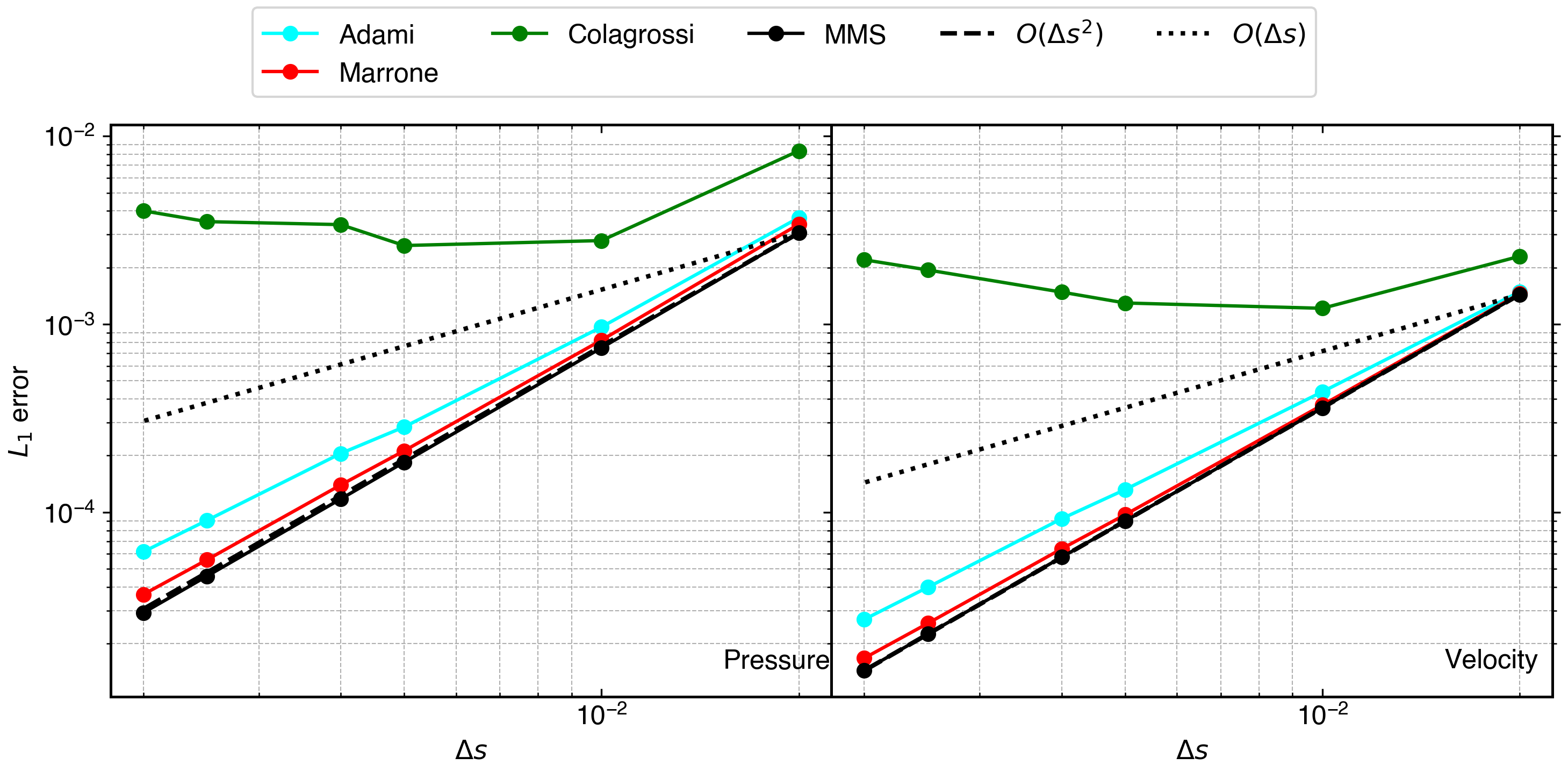}
\caption{$L_{1}$ error in pressure and velocity after 100 time steps for
different slip boundary implementations in the convex domain as shown in
\cref{fig:domains}.}
  \label{fig:d2_us_bc}
\end{figure}

\begin{figure}[htbp]
  \centering
  \includegraphics[width=0.8\linewidth]{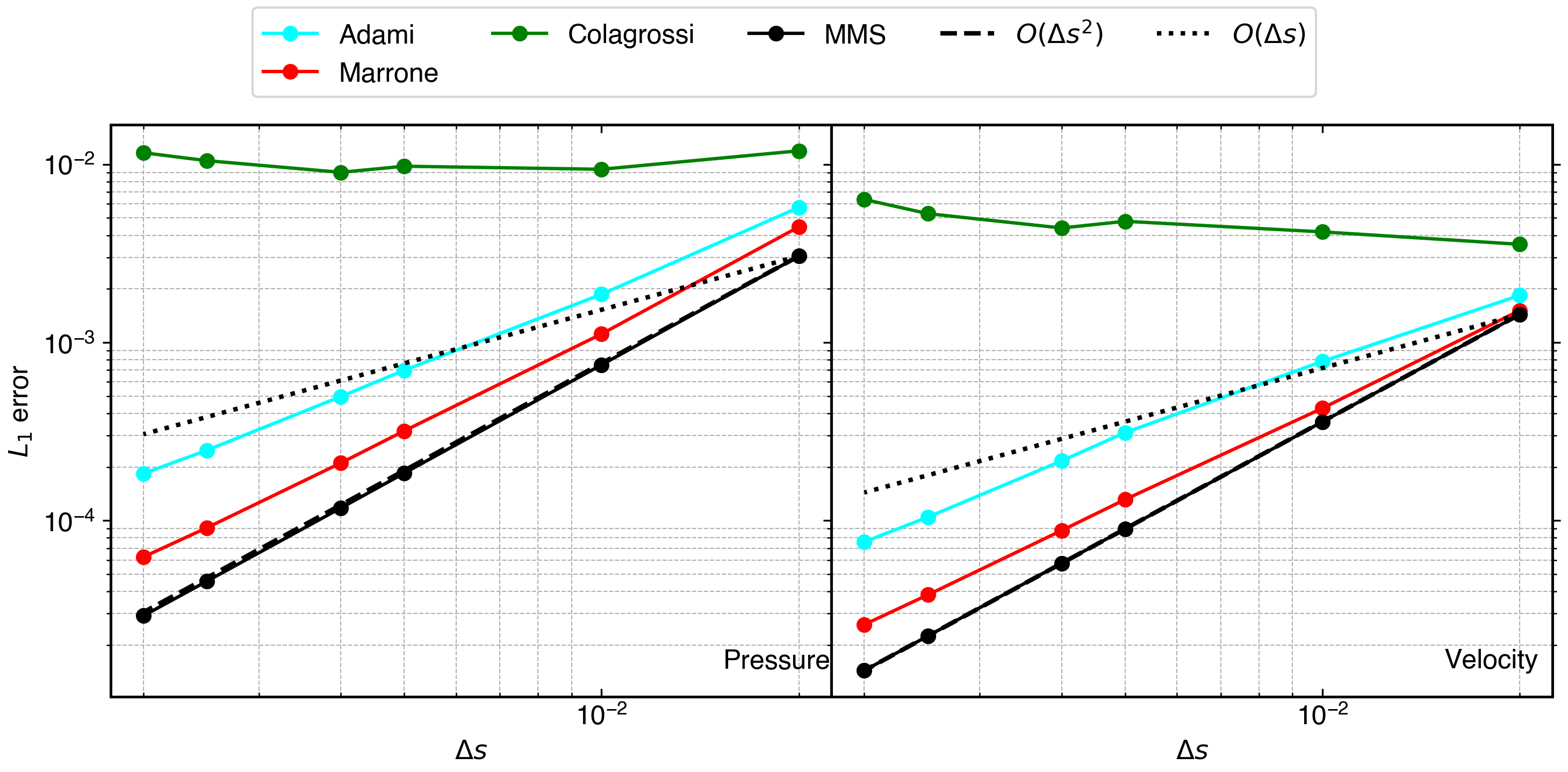}
\caption{$L_{1}$ error in pressure and velocity after 100 time steps for
different slip boundary implementations in the concave domain as shown in
\cref{fig:domains}.}
  \label{fig:d3_us_bc}
\end{figure}

We construct the MS for convex and concave domains in \cref{eq:mms_slip_d5}
that satisfies $\ten{u} \cdot \ten{n}=0$ at respective boundary surfaces of
interest. In \cref{fig:d2_us_bc} and \cref{fig:d3_us_bc}, we plot the $L_1$
error for pressure and velocity in convex and concave domains after 100
timesteps, respectively. Clearly, the method proposed by
\citet{colagross2003a} diverges for higher resolutions. The method proposed
by \citet{Adami2012} shows convergence rate close to 1.6 whereas the method
by \citet{marrone-deltasph:cmame:2011} is very close to second-order
convergence.

\begin{figure}[htbp]
  \centering
  \includegraphics[width=0.8\linewidth]{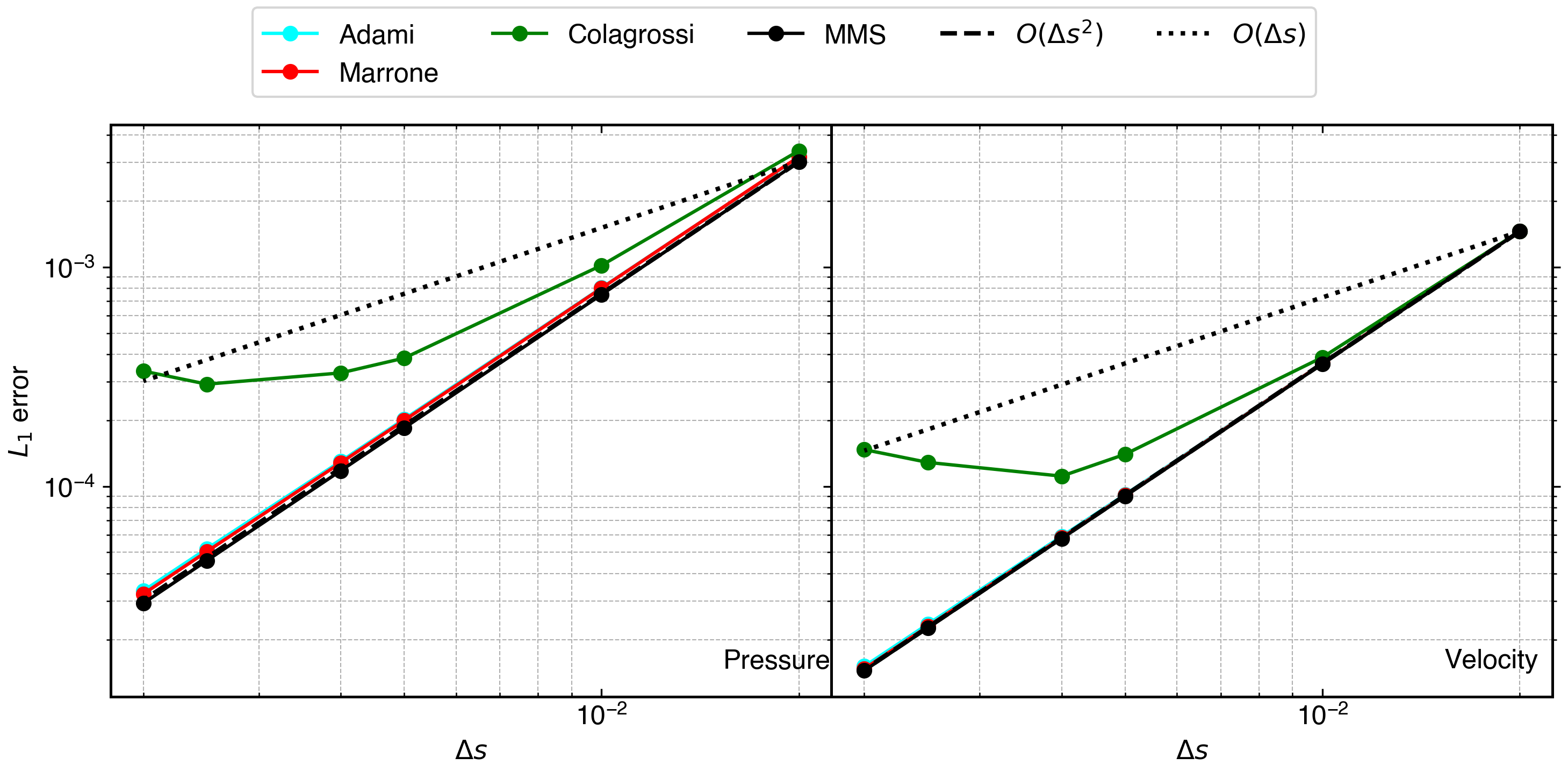}
\caption{$L_{1}$ error in pressure and velocity after 100 time steps for
different slip boundary implementations in the packed-convex domain as shown in
\cref{fig:domains_pack}.}
  \label{fig:d4_us_bc}
\end{figure}

\begin{figure}[htbp]
  \centering
  \includegraphics[width=0.8\linewidth]{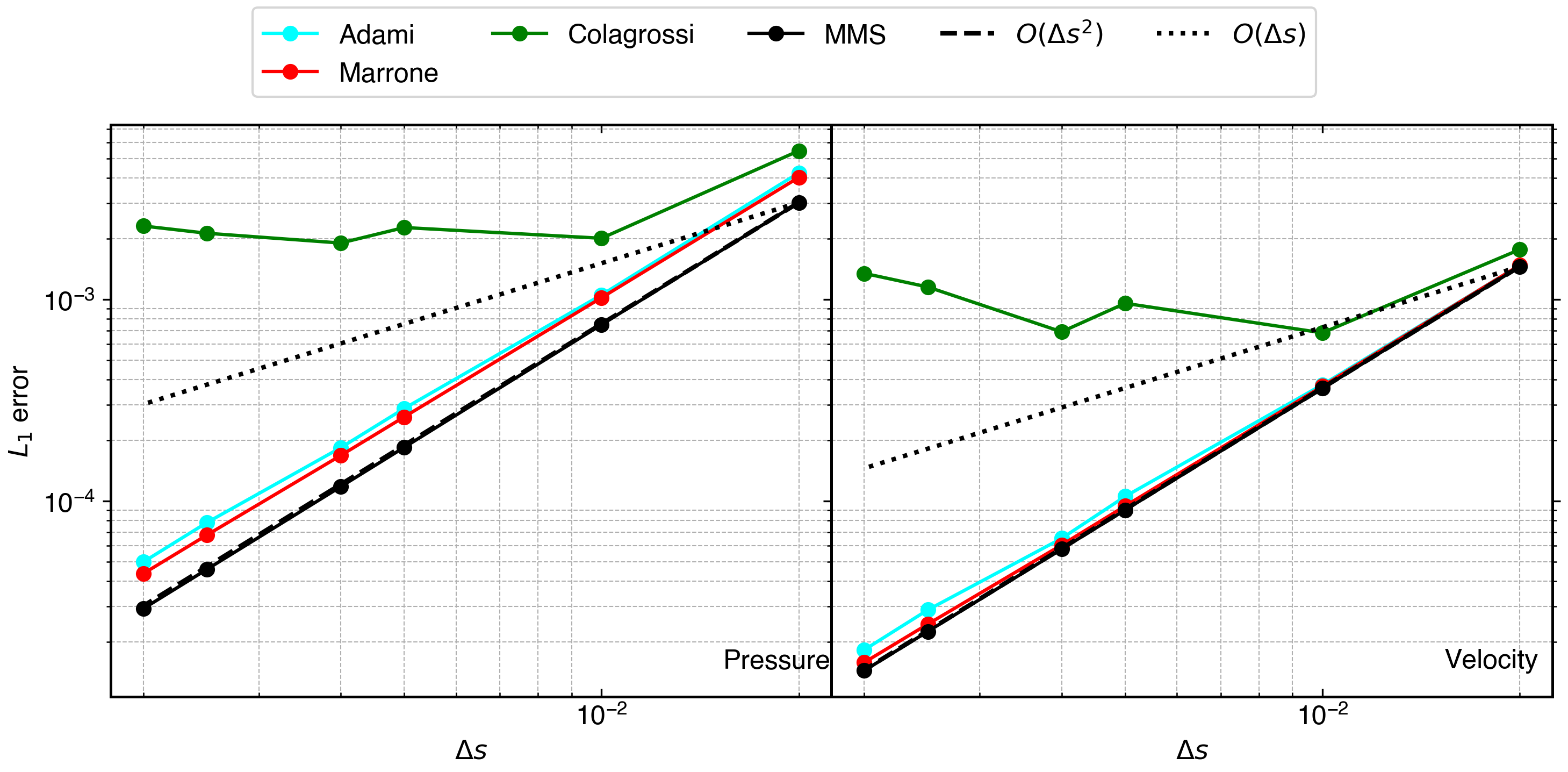}
\caption{$L_{1}$ error in pressure and velocity after 100 time steps for
different slip boundary implementations in the packed-concave domain as shown in
\cref{fig:domains_pack}.}
  \label{fig:d5_us_bc}
\end{figure}

We use the same MS for the packed-convex and packed-concave domains. In
\cref{fig:d4_us_bc} and \cref{fig:d5_us_bc}, we plot the $L_1$ error for
pressure and velocity for packed-convex and packed-concave domains after
100 timesteps, respectively. As expected, The order of convergence is
improved. In the packed domain, the convergence for the method by
\citet{Adami2012}, and \citet{marrone-deltasph:cmame:2011} shows second
order convergence. The method proposed by \citet{colagross2003a} converges
for lower resolutions in the case of the packed-convex domain but diverges in
the case of packed-concave domain. This shows that mirroring the fluid
particles for a curved surface does not result in a convergent boundary
condition implementation.

\subsubsection{No-slip boundary condition}

\begin{figure}[htbp]
  \centering
  \includegraphics[width=0.8\linewidth]{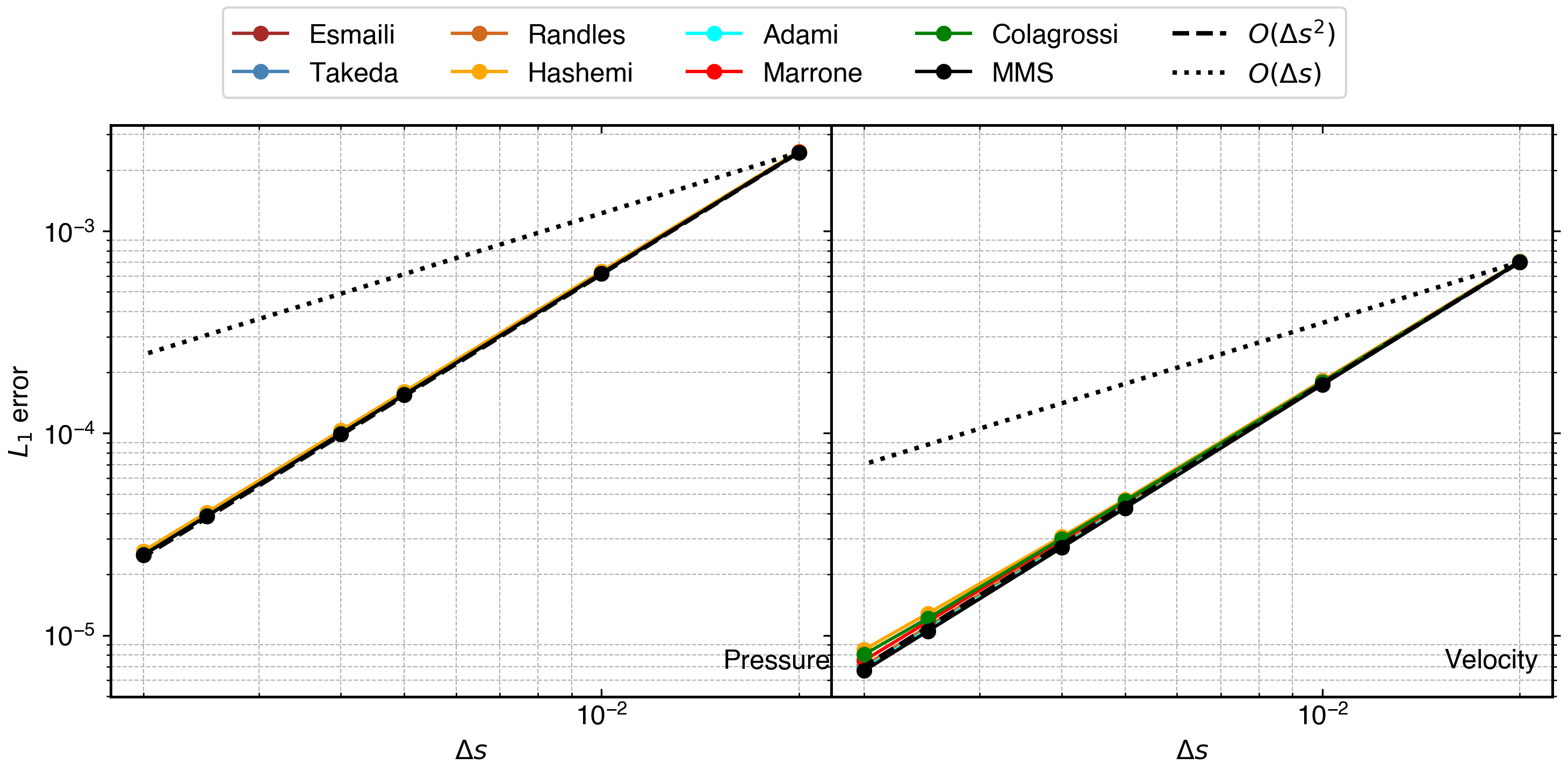}
\caption{$L_{1}$ error in pressure and velocity after 100 time steps for
different no-slip boundary implementations in the straight domain as shown in \cref{fig:domains}.}
  \label{fig:d1_uns_bc}
\end{figure}

In this section, we test different no-slip boundary implementations
discussed in \cref{sec:bc} using the domains used in the previous section. In
all the no-slip boundary condition implementations, we apply no-penetration
along with the no-slip boundary. In order to construct an MS for no-slip
boundary condition, we satisfy $\ten{u} = 0$ at the boundary. For the
straight domain, we use the MS in \cref{eq:mms_noslip_d1}. In
\cref{fig:d1_uns_bc}, we plot the $L_1$ error for pressure and velocity in
the domain after 100 timesteps. Clearly, all the methods show a convergence
rate very close to second order.

\begin{figure}[htbp]
  \centering
  \includegraphics[width=0.8\linewidth]{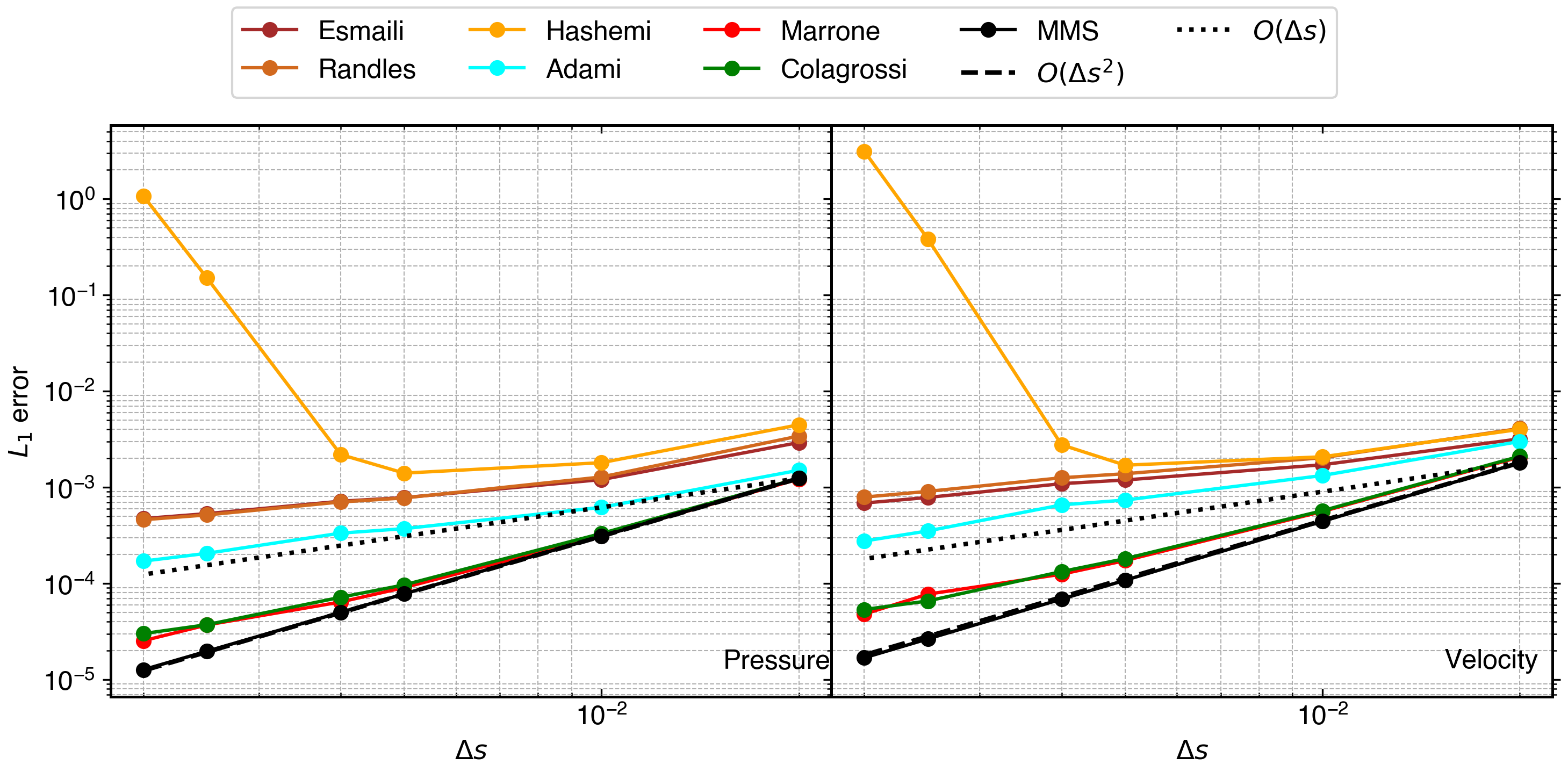}
\caption{$L_{1}$ error in pressure and velocity after 100 time steps for
different no-slip boundary implementations in the convex domain as shown in
\cref{fig:domains}.}
  \label{fig:d2_uns_bc}
\end{figure}

\begin{figure}[htbp]
  \centering
  \includegraphics[width=0.8\linewidth]{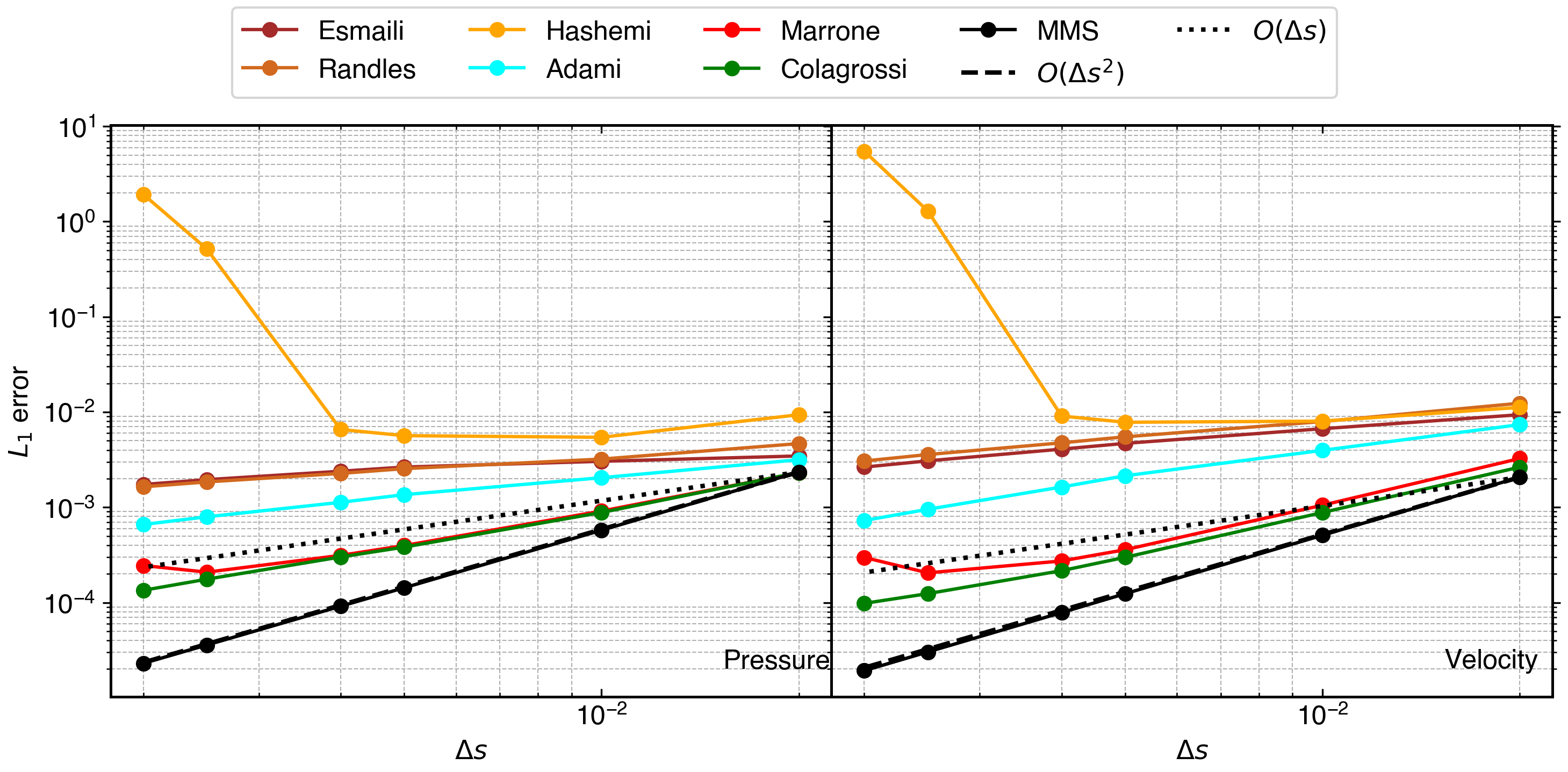}
\caption{$L_{1}$ error in pressure and velocity after 100 time steps for
different no-slip boundary implementations in the concave domain as shown in \cref{fig:domains}.}
  \label{fig:d3_uns_bc}
\end{figure}

Generally, we find that objects on which we intend to apply no-slip
boundary are curved. Therefore, we simulate all the methods on the domains
having curved surfaces. For the convex domain, we use the MS in
\cref{eq:mms_noslip_d5} whereas for concave domain, we use
\cref{eq:mms_noslip_d6}. In \cref{fig:d2_uns_bc} and \cref{fig:d3_uns_bc},
we plot the $L_{1}$ error in pressure and velocity for convex and concave
domains, respectively. Clearly, the method proposed by
\citet{hashemi_modified_2012} diverges for a curved surface. The errors in
the solutions are more in the concave domain compared to the convex domain.
The method by \citet{randles1996smoothed}, \citet{esmailisikarudi2016},
\citet{Adami2012} shows first-order convergence. Whereas methods by
\citet{colagross2003a}, and \citet{marrone-deltasph:cmame:2011} shows close
to 1.5. Some methods like by \citet{takeda1994a} cannot be applied on the
jagged boundary as some particles may lie on the surface, which may result
zero in the denominator of \cref{eq:takeda}.

\begin{figure}[htbp]
  \centering
  \includegraphics[width=0.8\linewidth]{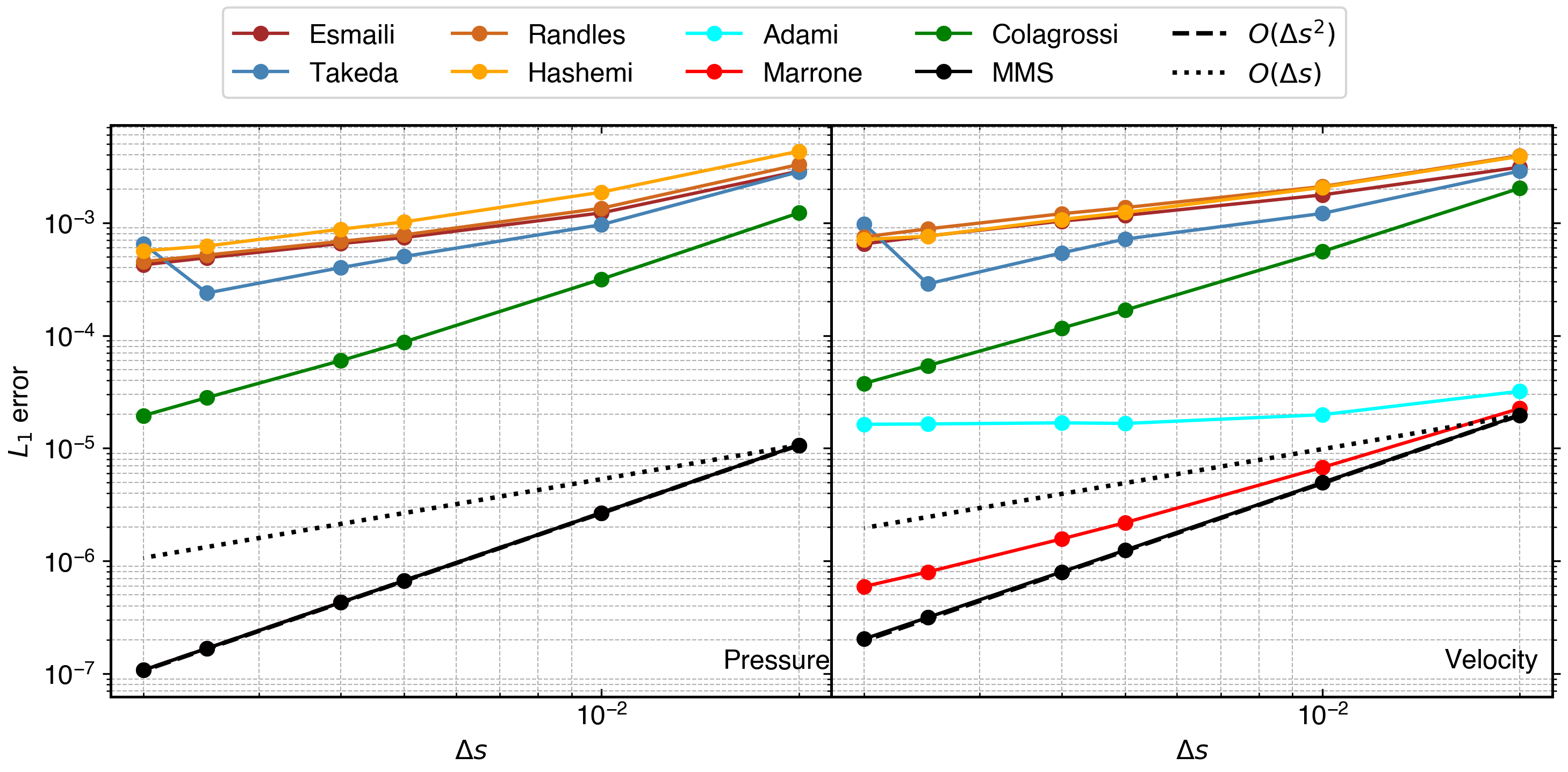}
\caption{$L_{1}$ error in pressure and velocity after 100 time steps for
different no-slip boundary implementations in the packed-convex domain as shown in
\cref{fig:domains_pack}. Note that Marrone results overlaps the results of MMS.}
  \label{fig:d4_uns_bc}
\end{figure}

\begin{figure}[htbp]
  \centering
  \includegraphics[width=0.8\linewidth]{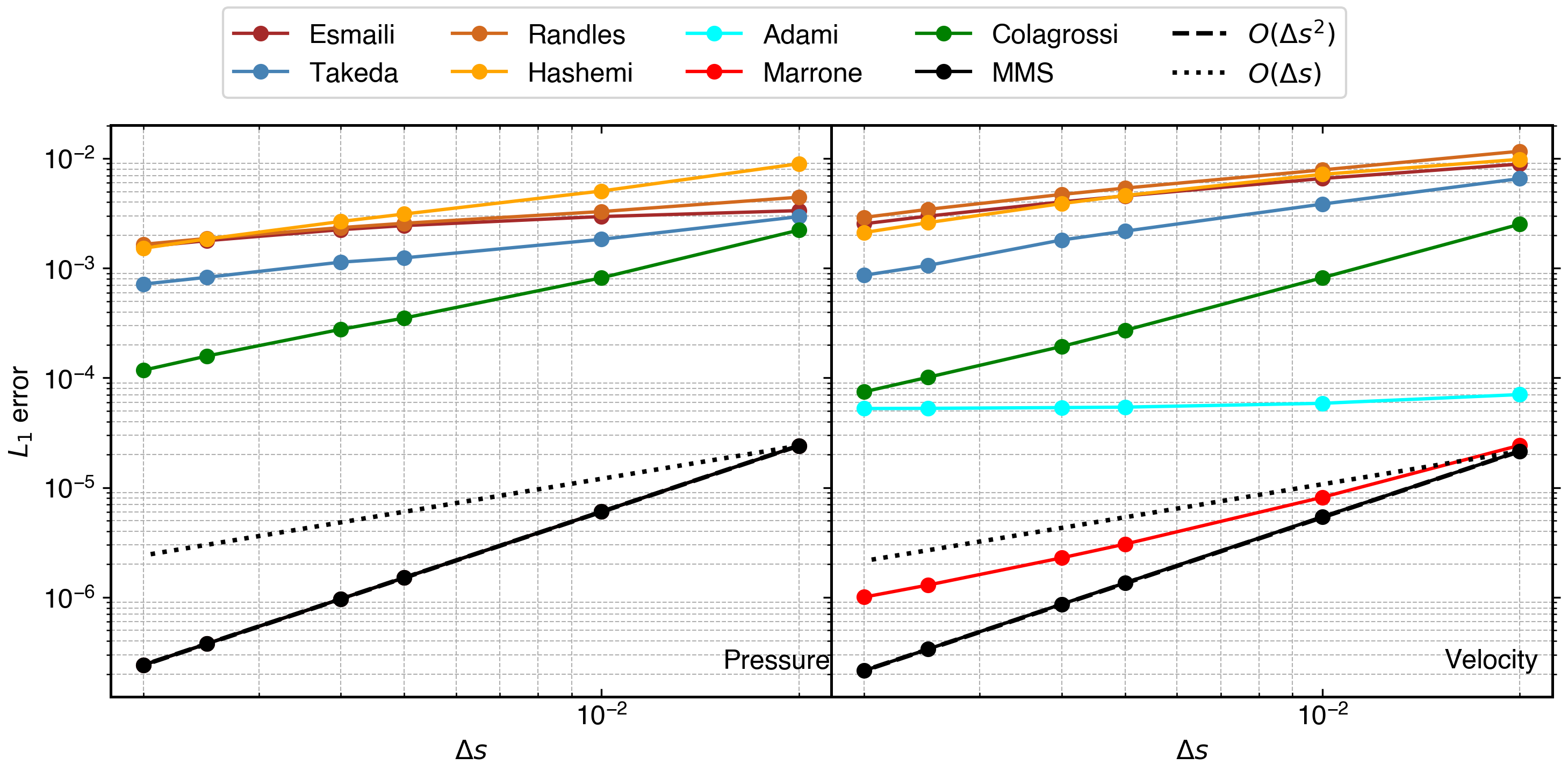}
\caption{$L_{1}$ error in pressure and velocity after 100 time steps for
different no-slip boundary implementations in the packed-concave domain as shown in
\cref{fig:domains_pack}. Note that Marrone results overlaps the results of MMS.}
  \label{fig:d5_uns_bc}
\end{figure}

Similar to other tests, we use the packed version of convex and concave
domains to test the convergence of the methods on packed domains. We use
the MS in \cref{eq:mms_noslip_d5}, and \cref{eq:mms_noslip_d6} for the
packed-convex and packed-concave domains, respectively. In the figure
\cref{fig:d4_uns_bc}, and \cref{fig:d5_uns_bc}, we plot the $L_1$ error in
pressure and velocity for packed-convex and packed-concave domains,
respectively. As expected, the convergence is improved. The method by
\citet{Adami2012} does not show any convergence due to zero order
interpolation used on the ghost particles. Method by \citet{takeda1994a},
\citet{hashemi_modified_2012}, \citet{randles1996smoothed},
\citet{esmailisikarudi2016} shows close to first order convergence. Clearly,
the method by \citet{marrone-deltasph:cmame:2011} shows convergence close
to the method when MS is used on the ghost particles. Further, the method
of \citet{colagross2003a} also shows good convergence, however the error
compared to \citet{marrone-deltasph:cmame:2011} method is 2 order of
magnitude higher. In the case of the packed-concave domain in \cref{fig:d5_uns_bc},
the order of convergence shown by all methods is lower compared to
packed-convex domain results.

\begin{table*}
  \centering
  \begin{tabular}{c|l|l|l}
    Method & Neumann Pressure & Slip & No-Slip\\
    \hline
    Adami~\cite{Adami2012}&($1.17\times10^{-5}$)2.00&($1.82\times10^{-5}$)1.95&($5.26\times10^{-5}$)0.09\\
    Colagrossi~\cite{colagross2003a}&($1.18\times10^{-5}$)2.00&($1.34\times10^{-3}$)0.19&($7.47\times10^{-5}$)1.52\\
    Esmaili~\cite{esmailisikarudi2016}&-&-&($2.54\times10^{-3}$)0.55\\
    Fourtakas~\cite{fourtakas2019}&($2.07\times10^{-5}$)1.66&-&-\\
    Hashemi~\cite{hashemi2012}&($4.39\times10^{-3}$)0.47&-&($2.11\times10^{-3}$)0.72\\
    Marongiu~\cite{marongiu2007}&($1.14\times10^{-3}$)0.63&-&-\\
    Marrone~\cite{marrone-deltasph:cmame:2011}&($1.15\times10^{-5}$)2.00&($1.58\times10^{-5}$)1.97&($1.00\times10^{-6}$)1.35\\
    Randles~\cite{randles1996smoothed}&-&-&($2.91\times10^{-3}$)0.62\\
    Takeda~\cite{takeda1994a}&-&-&($8.63\times10^{-4}$)0.86\\
  \end{tabular}
\caption{Table showing the summary of the error (in brackets) at the
resolution $500 \times 500$ and order of convergence of various boundary
condition methods in the packed-concave domain.}
  \label{tab:sum_bc_solid}
\end{table*}

In order to summarize the results, since the straight and convex domain
shows better results compared to concave domain, we consider the results
for a concave domain only. Furthermore, we compile results for a packed
domain only since the packed domains are preferred over the unpacked ones.
In the case of the no-slip and slip boundary, we focus on the convergence
of velocity, and in the case of the Neumann pressure, we focus only on the
convergence of pressure. In \cref{tab:sum_bc_solid}, we tabulate the error
at the highest resolution, i.e. $\Delta x = 1/500$, and the approximate
order of convergence for all the boundary conditions and methods. Clearly,
in the case of the Neumann pressure boundary condition, Adami, Colagrossi,
and Marrone converges well. In the case of the slip boundary conditon only
the Adami and Marrone methods work. Whereas in the case of the no-slip
boundary only Colagrossi and Marrone method show reasonable convergence.
Clearly, Marrone method is able to reach lowest error as well as show
convergence for all the types of boundary conditions.

\subsection{Comparison of open boundary condition implementations}

\begin{figure}[htbp]
  \centering
  \includegraphics[width=\linewidth]{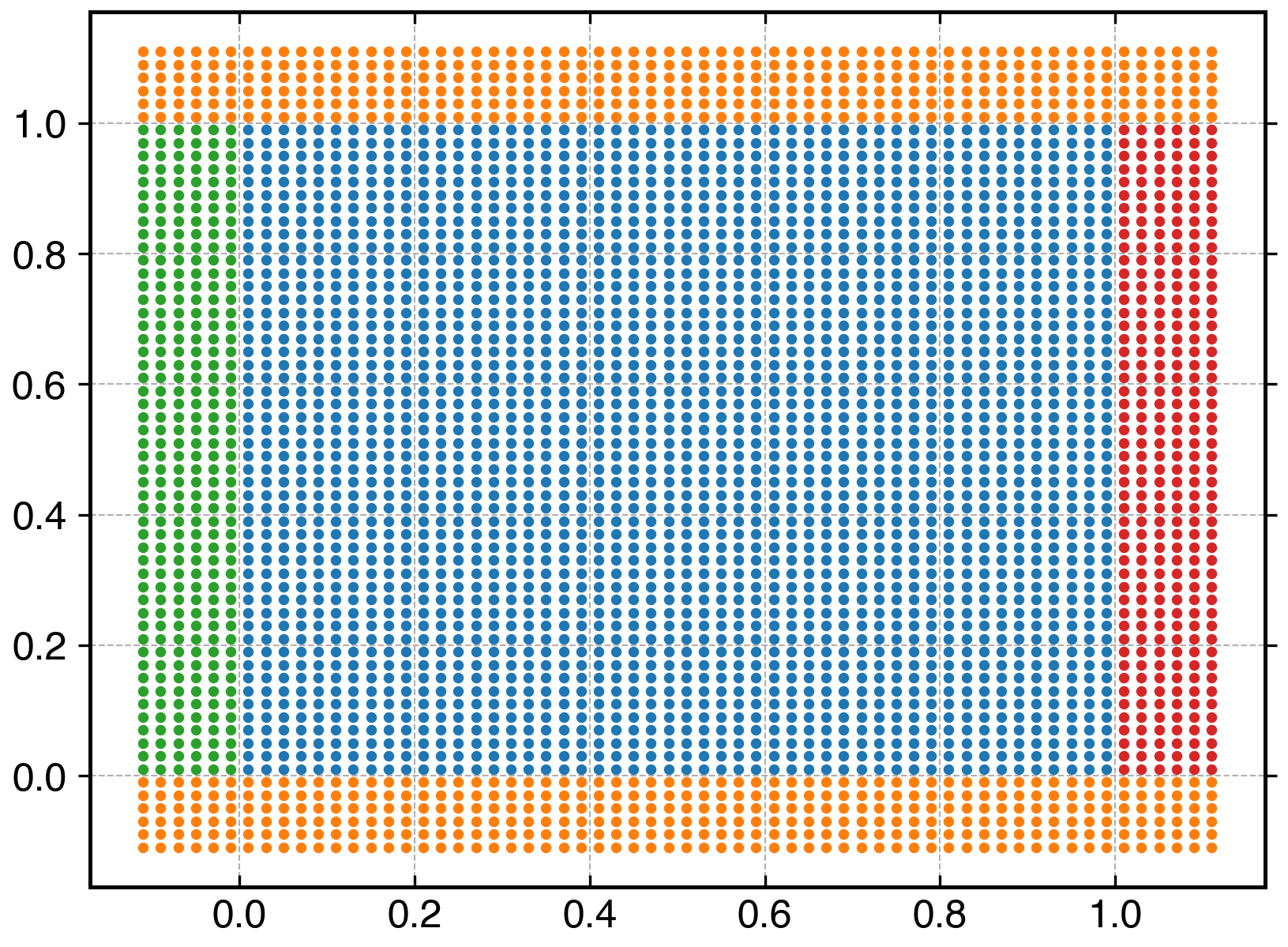}
\caption{The domain used for the verification of inlet and outlet boundary
implementation. The blue particles represent the fluid, orange particle
represent the wall, green particles are inflow particles, and red particles are the
outflow particles.}
  \label{fig:domain_io}
\end{figure}

In this section, we test various inlet and outlet boundary condition
implementations discussed in \cref{sec:obc}. In order to test the boundary
condition implementation, we require the inlet and outlet boundary to
continuously add and remove particles from the domain, respectively.
Furthermore, the inlet and outlet condition requires that the flow is only
along the normal at the boundary. In order to satisfy these conditions, we
use a $1m \times 1m$ domain, with an inlet and outlet on the left and right,
respectively. In \cref{fig:domain_io}, we show the domain with ghost
particles representing the inflow (in green), outflow (in red), and wall (in
orange). In the following sections, we discuss inlet boundary condition
implementations first, followed by outlet boundaries.

\subsubsection{Inlet boundary}

\begin{figure}[htbp]
  \centering
  \includegraphics[width=0.8\linewidth]{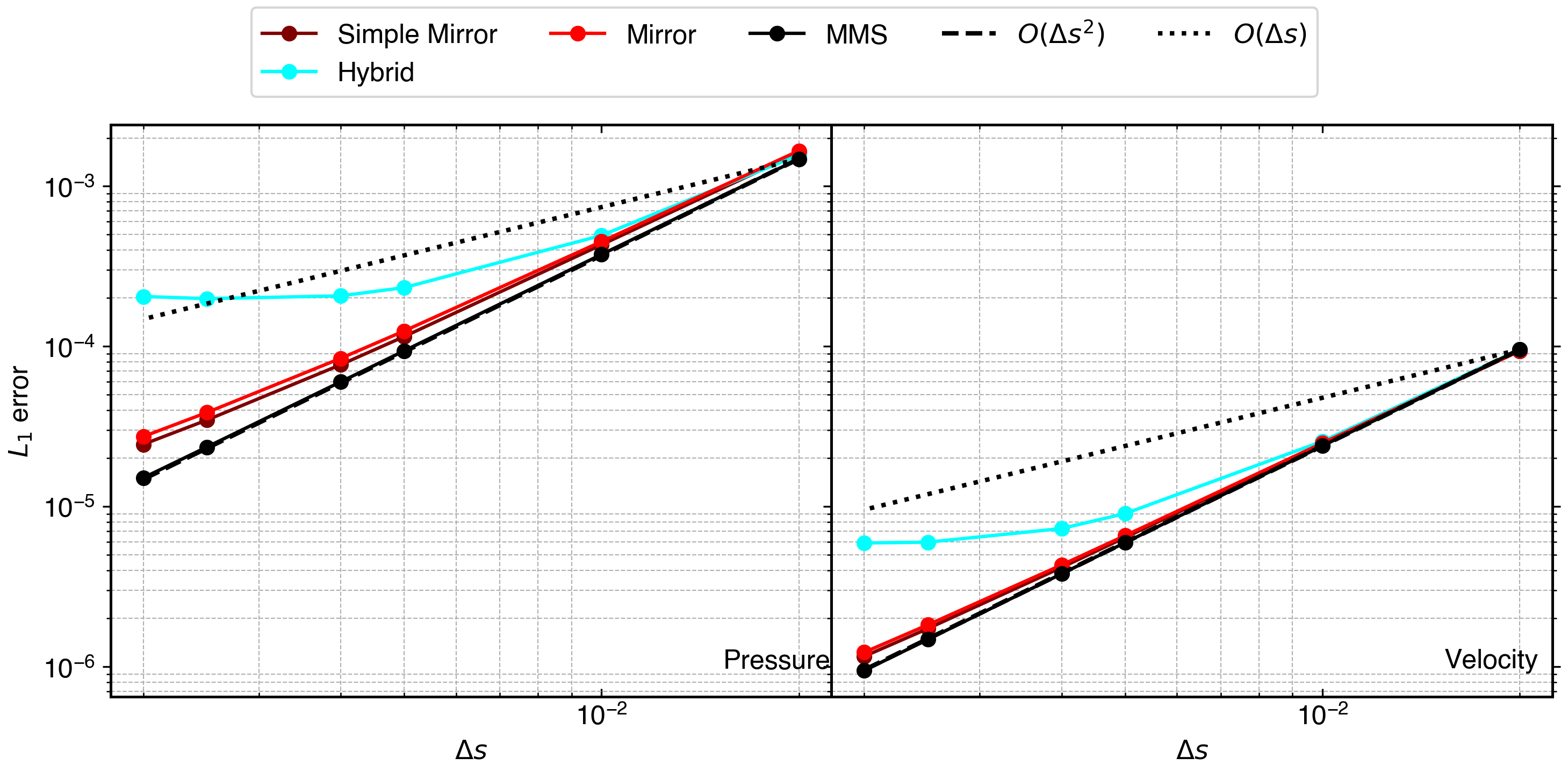}
\caption{$L_{1}$ error in pressure and velocity after 100 time steps for
different inlet velocity boundary implementations in the domain as shown in
\cref{fig:domain_io}.}
  \label{fig:u_in}
\end{figure}

In order to test the inlet velocity boundary condition, we use the MS in
\cref{eq:mms_io_vel}. In \cref{fig:u_in}, we plot the $L_1$ error in
pressure and velocity after 100 timesteps for all the velocity inflow
boundary implementations. We test all the methods discussed in \cref{sec:obc}
viz.\ mirror, simple-mirror, and hybrid. We observe that both mirror and
simple-mirror perform well for a velocity inlet boundary condition.
Whereas the hybrid is bounded by the limiting error in both pressure and
velocity.

\begin{figure}[htbp]
  \centering
  \includegraphics[width=0.8\linewidth]{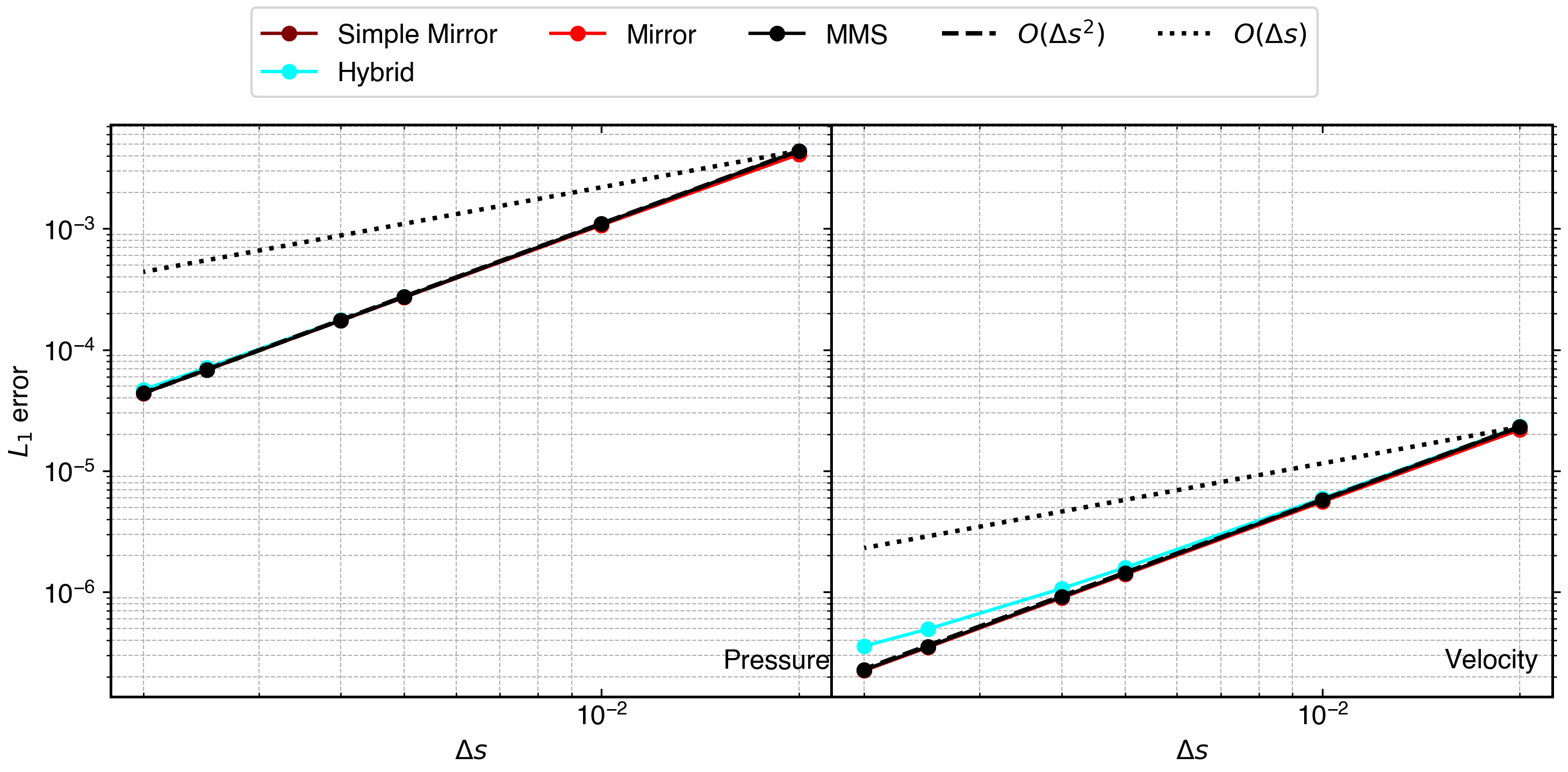}
\caption{$L_{1}$ error in pressure and velocity after 100 time steps for
different inlet pressure boundary implementations in the domain as shown in
\cref{fig:domain_io}.}
  \label{fig:p_in}
\end{figure}

In order to test the pressure inflow boundary implementation, we use the MS
in \cref{eq:mms_io_pres}. In \cref{fig:p_in}, we plot the $L_1$ error in
pressure and velocity after 100 timesteps for all the velocity inflow
boundary implementations. Clearly, the boundary implementation for a
pressure inflow boundary is second-order accurate for all the methods. In
the case of the hybrid method, a slight deviation in the convergence can be seen.

\begin{figure}[htbp]
  \centering
  \includegraphics[width=0.8\linewidth]{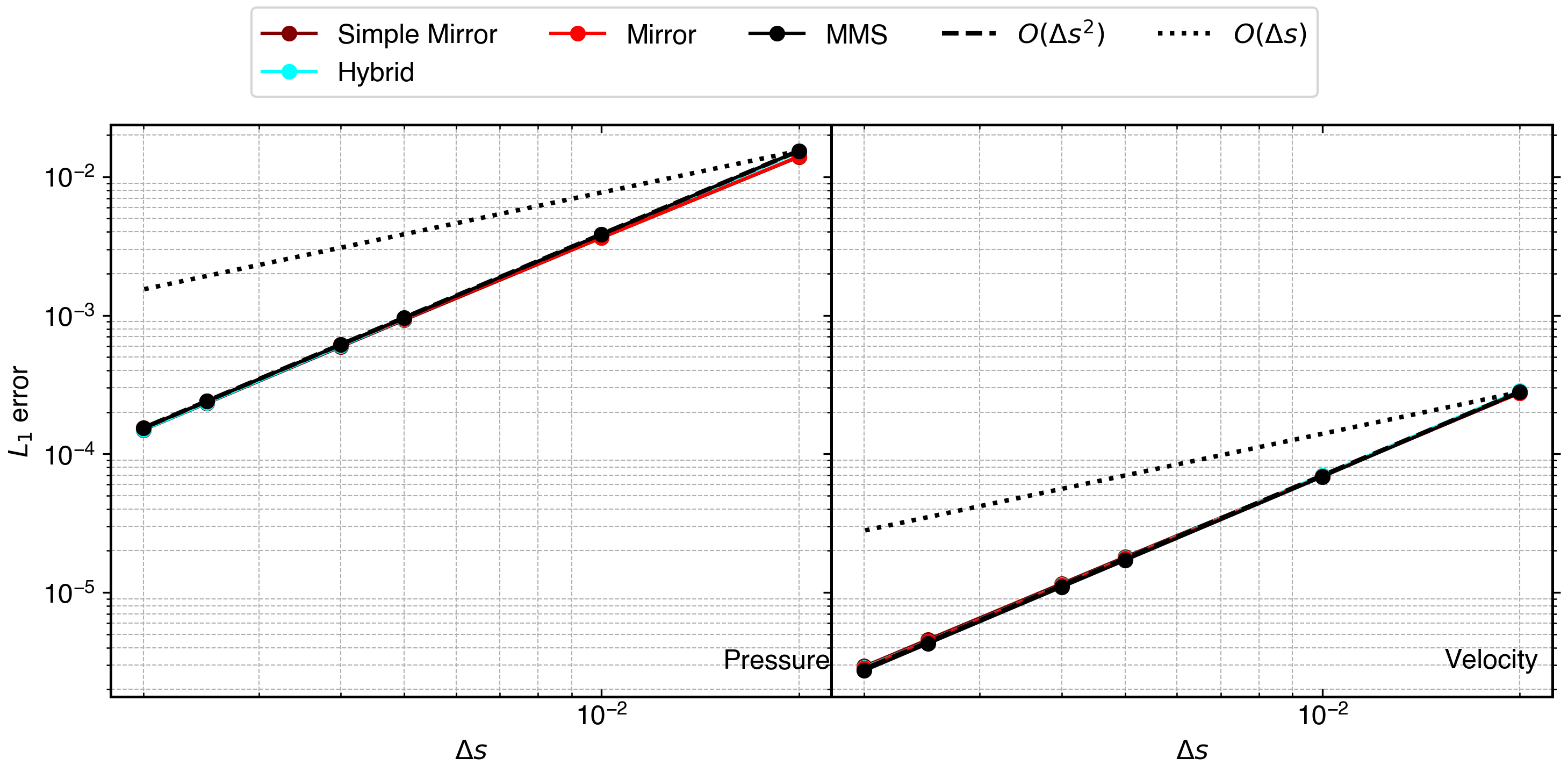}
\caption{$L_{1}$ error in pressure and velocity after 500 time steps for
different inlet velocity wave going upstream boundary implementations as shown in
the domain in \cref{fig:domain_io}.}
  \label{fig:v_wave_in}
\end{figure}

\begin{figure}[htbp]
  \centering
  \includegraphics[width=0.8\linewidth]{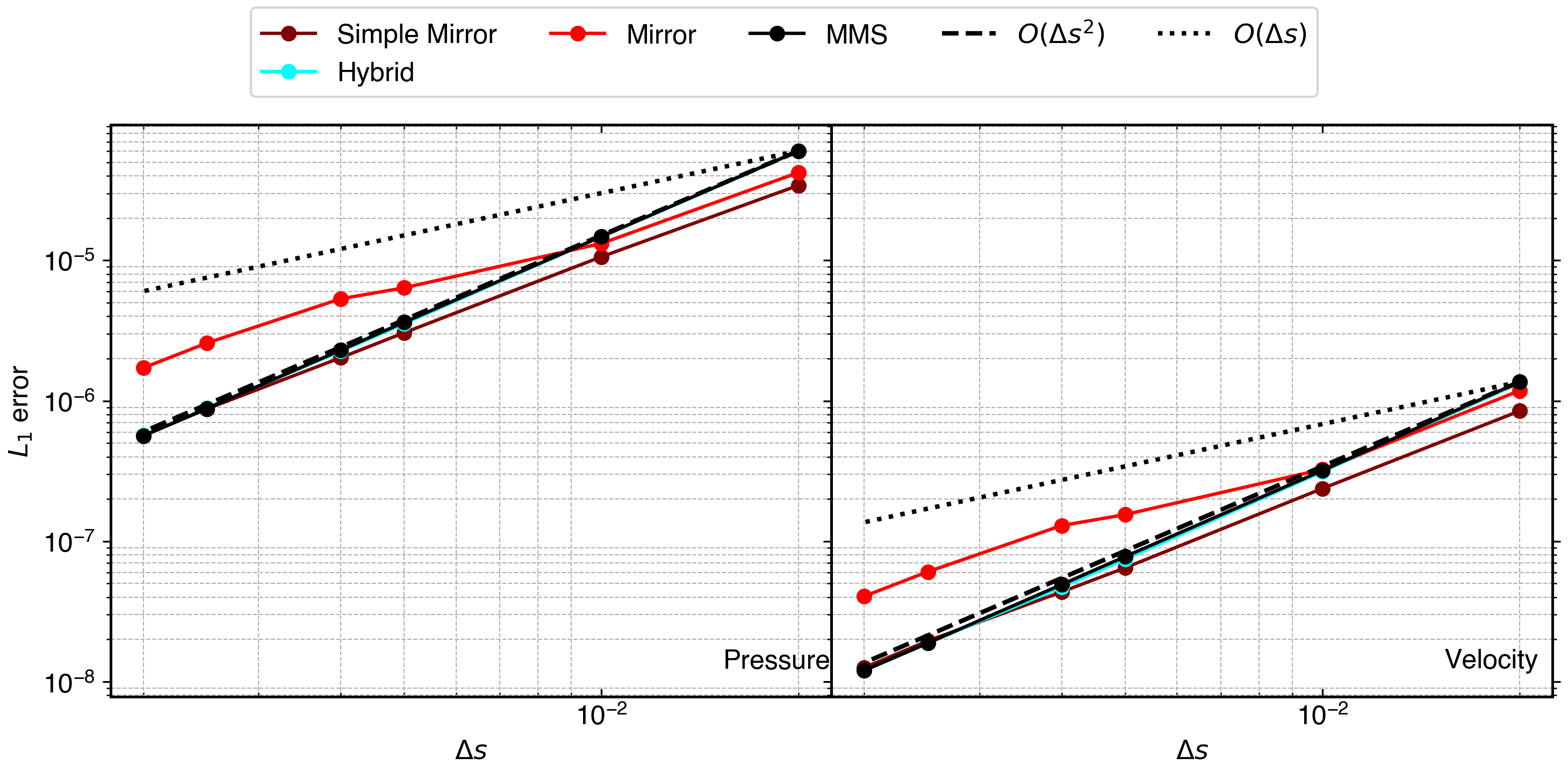}
\caption{$L_{1}$ error in pressure and velocity after 500 time steps for
different inlet pressure wave going upstream boundary implementations as shown in
the domain in \cref{fig:domain_io}.}
  \label{fig:p_wave_in}
\end{figure}

Since in WCSPH, due to weakly compressible assumption, waves travel with a
speed of artificial velocity of sound. We use the MS simulating a wave
passing out of the inlet, these kind of waves are encountered when a jump
start is performed on a wind tunnel kind of simulation. We simulate the
problem for 500 iterations in order to allow the wave to completely pass
through the inlet/outlet. We use the MS in \cref{eq:mms_in_vel} for the
inlet velocity wave. In \cref{fig:v_wave_in}, we plot the $L_1$ error
in pressure and velocity for all the methods. Clearly, the hybrid method
also shows second-order convergence along with other methods.

In order to simulate a pressure wave going out of the inlet, we use the MS
in \cref{eq:mms_in_pres}. In \cref{fig:p_wave_in}, we plot the $L_1$ error
in pressure and velocity for all the methods. The mirror method shows a
slight increase in error for higher resolutions. The simple-mirror remains
at the same level of error compared to the hybrid method, which shows
second-order convergence.

\subsubsection{Outlet boundary}

\begin{figure}[htbp]
  \centering
  \includegraphics[width=0.8\linewidth]{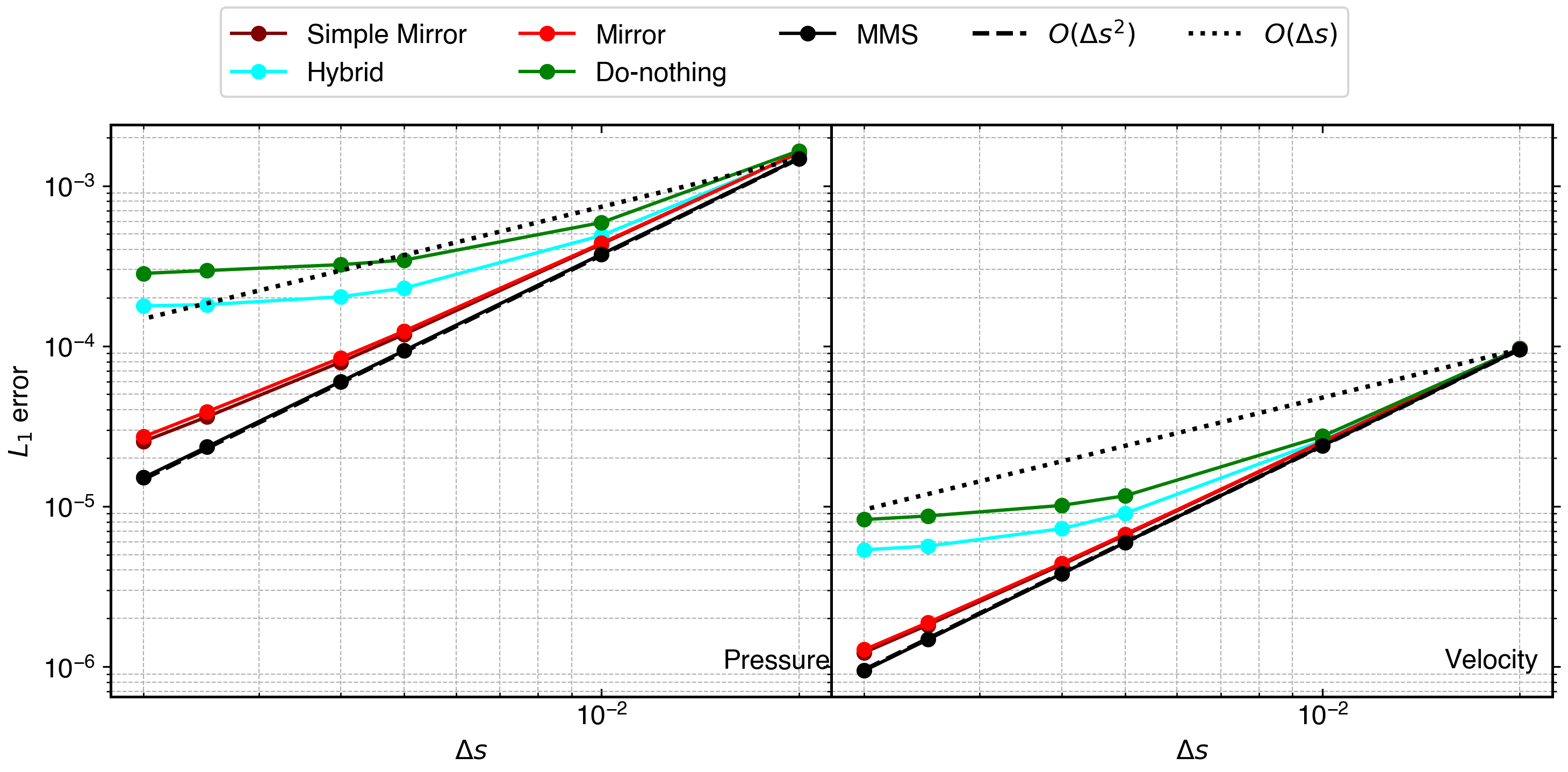}
\caption{$L_{1}$ error in pressure and velocity after 100 time steps for
different outlet velocity boundary implementations in the domain as shown in
\cref{fig:domain_io}.}
  \label{fig:u_out}
\end{figure}

The outflow is different compared to the inlet as we usually do not have
any information about the ghost particles in these regions. In order to test
the outflow velocity boundary condition, we use the MS in
\cref{eq:mms_io_vel}. In \cref{fig:u_out}, we plot the $L_1$ error in
pressure and velocity after 100 timesteps for all the velocity outflow
boundary implementation. The do-nothing and the hybrid boundary are both
bounded by a limiting error which is proportional to the speed of sound. As
before, both mirror and simple-mirror show second-order convergence for
velocity outlet boundary condition.

\begin{figure}[htbp]
  \centering
  \includegraphics[width=0.8\linewidth]{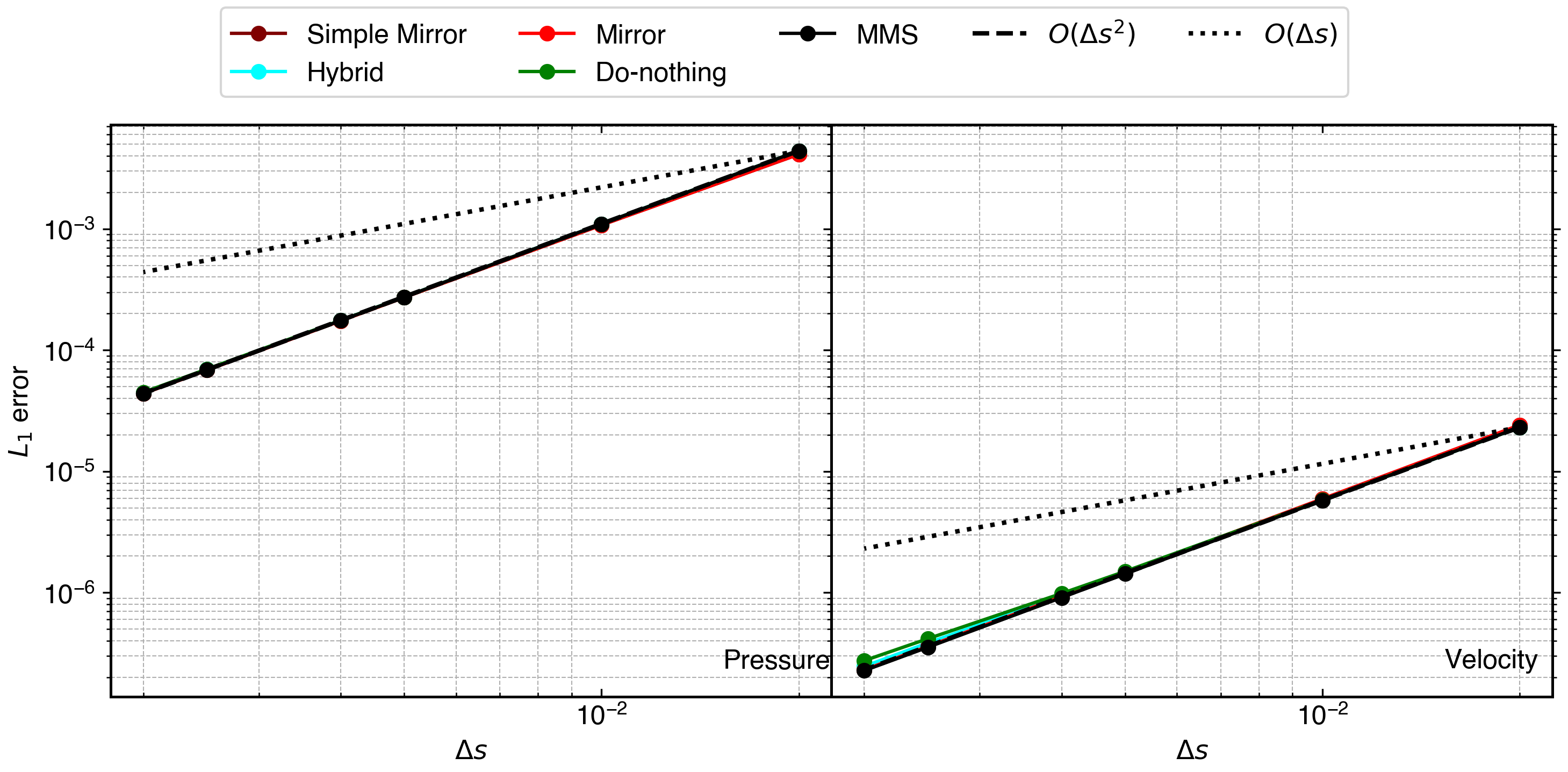}
\caption{$L_{1}$ error in pressure and velocity after 100 time steps for
different outlet pressure boundary implementations in the domain as shown in
\cref{fig:domain_io}.}
  \label{fig:p_out}
\end{figure}

In order to test the pressure outflow boundary implementation, we use the MS
in \cref{eq:mms_io_pres}. In \cref{fig:p_out}, we plot the $L_1$ error in
pressure and velocity after 100 timesteps for all the velocity outflow
boundary implementation. Clearly, all the methods show second-order
convergence.

\begin{figure}[htbp]
  \centering
  \includegraphics[width=0.8\linewidth]{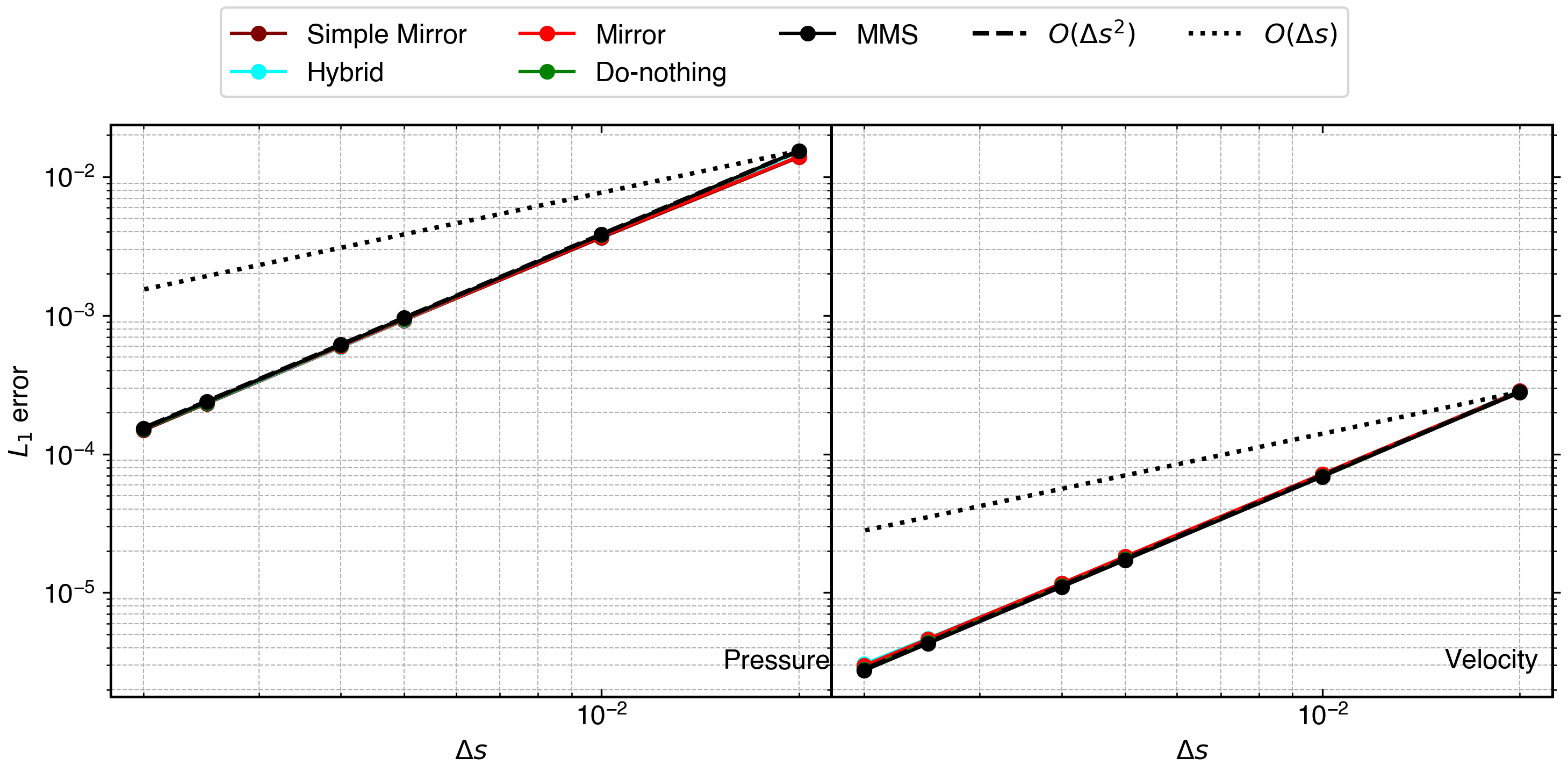}
\caption{$L_{1}$ error in pressure and velocity after 500 time steps for
different outlet velocity wave going downstream boundary implementations in
the domain as shown in \cref{fig:domain_io}.}
  \label{fig:v_wave_out}
\end{figure}

\begin{figure}[htbp]
  \centering
  \includegraphics[width=0.8\linewidth]{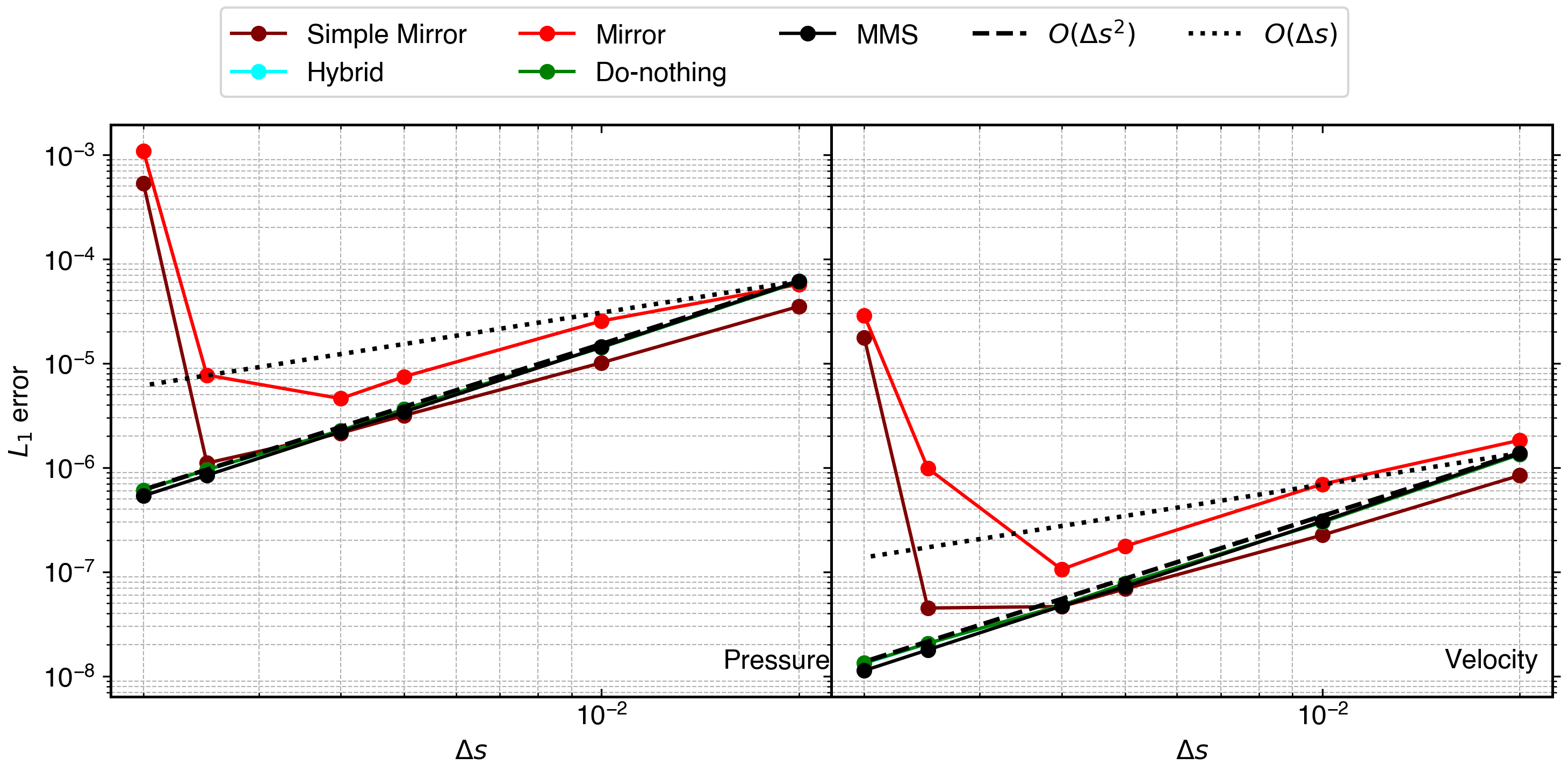}
\caption{$L_{1}$ error in pressure and velocity after 500 time steps for
different outlet pressure boundary wave going downstream implementations in
the domain as shown in \cref{fig:domain_io}.}
  \label{fig:p_wave_out}
\end{figure}

To investigate the behavior of the outlet boundary implementation under the
influence of a passing wave, we use the MS given in \cref{eq:mms_out_pres}
and \cref{eq:mms_out_vel} for outlet pressure and outlet velocity
boundary implementations, respectively. In \cref{fig:v_wave_out}, and
\cref{fig:p_wave_out}, we plot the $L_1$ error for pressure and velocity
after 500 timesteps for both the MS. In the case of the velocity wave, all the
methods show second order convergence. However, in the case of the pressure
wave in \cref{fig:p_wave_out}, the mirror and simple-mirror method
diverges. This shows that the mirror and simple-mirror methods are not
truly non-reflecting and are unable to pass a pressure wave. These results
support the finding of \citet{negi2021numerical}, where a short domain
with the mirror outlet was found to be unstable.

\begin{table}
  \centering
  \begin{tabular}{c|l|l|l|l}
    Method & Velocity in & Pressure in \\
    \hline
    Hybrid\cite{negi_nrbc_2020}&$(2.87 \times 10^{-6})$ 2.00&$(5.66 \times 10^{-7})$ 2.02\\
    Mirror\cite{tafuni_io_2018}&$(2.84 \times 10^{-6})$ 1.97&$(1.71 \times 10^{-6})$ 1.38\\
    Simple Mirror&$(2.93 \times 10^{-6})$ 1.98&$(5.64 \times 10^{-7})$ 1.81\\
  \end{tabular}
  \caption{Summary of results for the wave traveling upstream out of the
  inlet for all the methods. Error at highest resolution is shown in brackets.}
  \label{tab:comp_obc_in}
\end{table}

\begin{table}
  \centering
  \begin{tabular}{c|l|l|l|l}
    Method &Velocity out & Pressure out \\
    \hline
    Do-nothing\cite{federico2012}&$(2.80 \times 10^{-6})$ 2.00&$(6.08 \times 10^{-7})$ 2.00\\
    Hybrid\cite{negi_nrbc_2020}&$(3.05 \times 10^{-6})$ 1.99&$(6.05 \times 10^{-7})$ 2.00\\
    Mirror\cite{tafuni_io_2018}&$(2.97 \times 10^{-6})$ 1.97&$(1.01 \times 10^{-3})$ -3.63\\
    Simple Mirror&$(2.86 \times 10^{-6})$ 1.97&$(5.34 \times 10^{-4})$ -4.21\\
  \end{tabular}
  \caption{Summary of results for the wave traveling downstream out of the
  outlet for all the methods. Error at highest resolution is shown in brackets.}
  \label{tab:comp_obc_out}
\end{table}

In order to summarize the results for the open boundary conditions, we
consider only the results for the traveling wave since, in WCSPH, it is
important that the waves that are generated must be allowed to pass through
inlet/outlet without affecting the flow. In case of a velocity wave, we
focus on errors in velocity, whereas in the case of a pressure wave, we
focus on errors in pressure. In \cref{tab:comp_obc_in} and
\cref{tab:comp_obc_out}, we tabulate the error for the highest resolution
and the approximate order of convergence for all the methods simulating
traveling wave MS. Clearly, the mirror method show significant decrease in
order of convergence in the case of the pressure wave moving upstream. In
case of the wave traveling downstream, both mirror and simple-mirror
diverge. The hybrid method is applicable and converge for both the
scenarios.

\subsection{Performance Comparison}
\label{sec:perf}

In this section, we compare the performance of three solid boundary
condition implementations, viz. Marrone, Colagrossi, and Adami. For this
testcase, we use a 10 core, dual socket Intel(R) Xeon(R) CPU E5-2650 v3
processor CPU. In the context of complexity, the Adami method requires only
one loop over all the fluid particles to extrapolate properties from fluid,
whereas the Colagrossi method requires the creation and deletion of
particles in every timestep. In the case of the Marrone method, one needs
to solve an additional $4 \times 4$ matrix for each particle in the solid
boundary, excluding the extrapolation step.

\begin{figure}[htbp]
  \centering
  \includegraphics[width=0.8\linewidth]{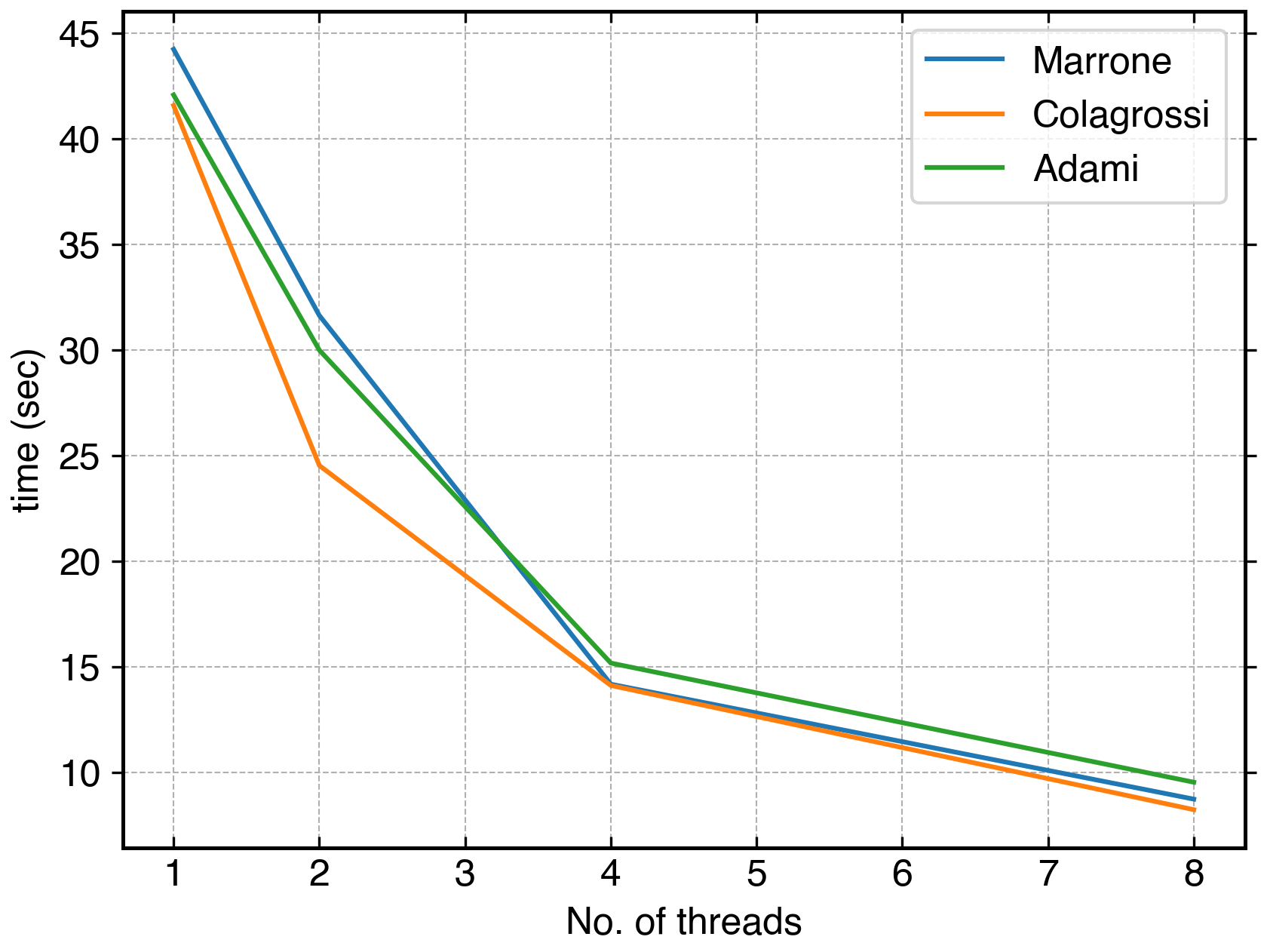}
\caption{The time taken with the increase in the number of threads for
different solid boundary conditions.}
  \label{fig:perf}
\end{figure}

In \cref{fig:perf}, we plot the time taken versus the no of parallel
computing threads for 100 timesteps for a $100 \times 100$ domain. Clearly,
the time taken are very close despite the fact that different amount of
computations are required. This is due to a very low number of solid particle
compared to the fluid particles. In the present case, the fluid particles
are $10000$, whereas the solid particles on which boundary condition is
implemented are $600$. Therefore, we demonstrated that a higher order
boundary implementation does not affect the performance of the solver.
Furthermore, these results can be easily extrapolated to open-boundary
implementations. In the next section, we propose an algorithm to obtain a
convergent solver for a problem containing inlet, outlet, and solid
boundaries.

\section{Complete second order convergent SPH scheme}
\label{sec:discussion}

In the previous section, we have shown that some of the boundary condition
implementations are convergent using the MMS. Recently,
\citet{negiHowTrainYour2021a}, and \citet{negi2021numerical} have used
verification methods to procedurally obtain a second order accurate WCSPH
scheme without boundary as discussed in \cref{sec:sph}. In this section, we
extend their method to propose a second-order convergent scheme with
second-order convergent boundary implementations. For brevity, we use short
names to represent an equation in the algorithm. For example,
\textit{EvaluateVelocityOnGhost(dest, sources)} can be written
formally as shown in algorithm \ref{alg:eq}.

\begin{algorithm}[htbp]
\caption{Psuedo-code of equation \textit{EvaluateVelocityOnGhost}.}\label{alg:eq}
\For(){i in dest}{
  $u_i = 0$\;
  \For(){j in sources}{
    $u_i  = u_i + u_j \omega_j W_{ij}$\;
}
}
\end{algorithm}
In algorithm \ref{alg:eq}, $i$ is the loop index and $j$ is over all the
neighbors of $i$\textsuperscript{th} element. We note that all the
extrapolation/equation require a corrected kernel/gradient and is implied.

Considering a fluid domain, the inlet continuously feeding particles to the
fluid, and the outlet continuously consuming particles from the fluid and an
arbitrary shaped solid body. In algorithm \ref{alg:soc}, we show the algorithm for a
convergent time-accurate WCSPH scheme to simulate the flow for a given
initial condition. We denote all the fluid particles by $\mathscr{F}$, all
solid particles by $\mathscr{S}$, and all inlet/outlet particles by
$\mathscr{IO}$. All the virtual particles required for solid particles are
denoted by $\mathscr{M(S)}$.

\begin{algorithm}[htbp]
\caption{A Second order convergent scheme.}\label{alg:soc}
\While{$t < t_{final}$}{
\For(){i in $\mathscr{F}$}{
  EvaluatePressure(dest=i, sources=$\phi$)\;
}
\For(){i in $\mathscr{IO}$}{
EvaluateVelocityOnInletOutlet(dest=i, sources=$\mathscr{F}$)\;
EvaluatePressureOnInletOutlet(dest=i, sources=$\mathscr{F}$)\;
}
\For(){i in $M(\mathscr{S})$}{
EvaluateVelocityOnGhost(dest=i, sources=$\mathscr{F} \cup \mathscr{IO}$)\;
EvaluatePressureOnGhost(dest=i, sources=$\mathscr{F} \cup \mathscr{IO}$)\;
}
\For(){i in $\mathscr{S}$}{
EvaluateSlipVelocityOnSolidFromGhost(dest=i, sources=$\phi$)\;
EvaluateNoSlipVelocityOnSolidFromGhost(dest=i, sources=$\phi$)\;
EvaluatePressureOnSolidFromGhost(dest=i, sources=$\phi$)\;
}
\For(){i in $\mathscr{F} \cup \mathscr{S}$}{
ComputeVelocityGradient(dest=i, sources=$\mathscr{F} \cup \mathscr{IO}$)\;
\CommentSty{\# use extrapolated no slip velocity for solid}\\
ComputeVelocityGradientSolid(dest=i, sources=$\mathscr{S}$)\;
}
\For(){i in $\mathscr{F}$}{
ContinuityEquation(dest=i, sources=$\mathscr{F} \cup \mathscr{IO}$)\;
\CommentSty{\# use extrapolated slip velocity for solid}\\
ContinuityEquationSolid(dest=i, sources=$\mathscr{S}$)\;
PressureForces(dest=i, sources=$\mathscr{F} \cup \mathscr{S} \cup \mathscr{IO}$)\;
ComputeViscousForces(dest=i, sources=$\mathscr{F} \cup \mathscr{IO}$)\;
}
\For(){i in $\mathscr{F} \cup \mathscr{IO}$}{
  Integrate(dest=i, sources=$\phi$)
}

}
\end{algorithm}

The algorithm starts with the evaluation of pressure from the equation of
state given by \cref{eq:eos} in \textit{EvaluatePressure}. In the next
step, we evaluate pressure and velocity on inlet and outlet domain using
the hybrid method~\cite{negi_nrbc_2020} in
\textit{EvaluatePressureOnInletOutlet}, and
\textit{EvaluateVelocityOnInletOutlet}, respectively. We note that the
inlet and outlet properties are updated, and then these are used as the
source for solid properties in case of overlaps. We use the method by
\citet{marrone-deltasph:cmame:2011} to evaluate properties on solid
particles. We first evaluate first order accurate values on virtual
particles in \textit{EvaluateVelocityOnGhost}, and
\textit{EvaluatePressureOnGhost} and then use these values to obtain
pressure, slip, and no-slip velocities in
\textit{EvaluatePressureOnSolidFromGhost},
\textit{EvaluateSlipVelocityOnSolidFromGhost},
\textit{EvaluateNoSlipVelocityOnSolidFromGhost}, respectively. In order to
obtain second-order accurate viscous operator, we require velocity gradient
on each particle, which is computed in \textit{ComputeVelocityGradient},
and \textit{ComputeVelocityGradientSolid}. We note that we consider no-slip
extrapolated velocity for gradient computation. Finally, we evaluate
accelerations due to various forces in \textit{ContinuityEquation,
ContinuityEquationSolid, PressureForces, ComputeViscousForces}. We note
that, in the evaluation of the continuity equation, we use extrapolated
slip velocity on solids~\cite{muta2020efficient}. We use the computed and
extrapolated properties to integrate the particle properties and position.
We use IPST to make the particles more uniform and update the properties
using a first-order accurate correction. The IPST can be performed after
every few timesteps. In our simulations, we perform shifting after every 10
iterations.

\begin{figure}[htbp]
  \centering
  \includegraphics[width=0.8\textwidth]{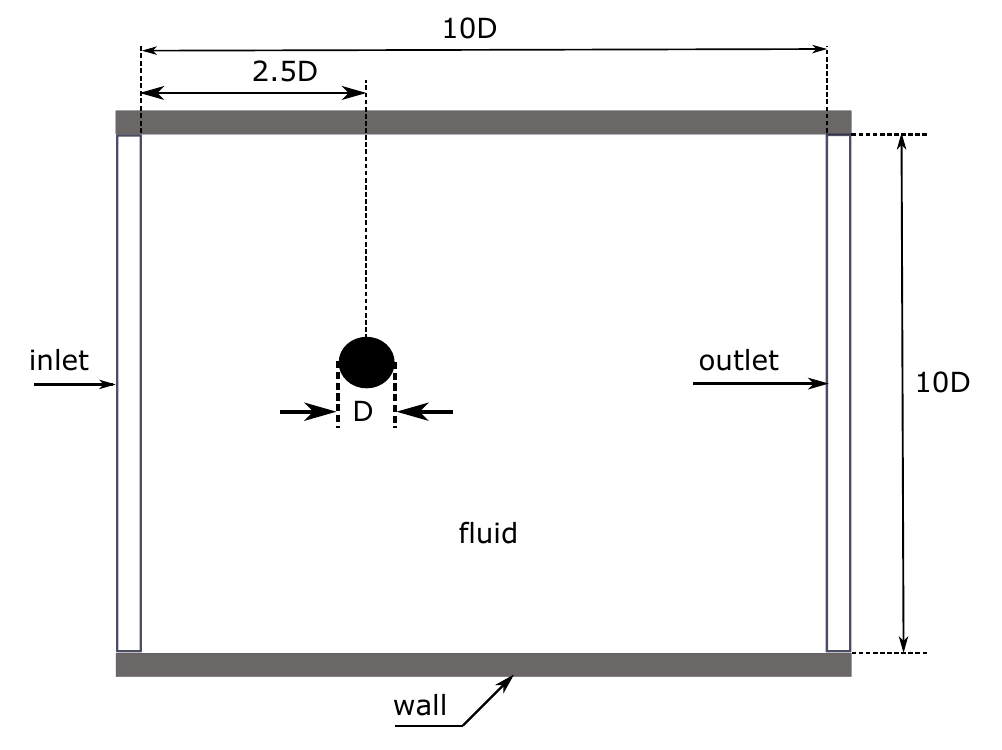}
  \caption{Description of the domain of the flow fast cylinder problem.}
  \label{fig:fpc_domain}
\end{figure}

In order to show the accuracy achieved by the proposed algorithm, we solve
the flow past a circular cylinder. We consider the domain shown in
\cref{fig:fpc_domain}. We consider inflow velocity $U = 1m/s$, Reynolds
number $Re=200$, and a cylinder of diameter $D=2m$. We set the dynamic
viscosity $\nu = U D/Re$. We discretize the domain with $\Delta x = D/40$.
The total particles in the domain are approximately $0.18M$. We simulate
the problem using the algorithm \ref{alg:soc} with artificial speed of
sound $c_o=10m/s$ for $200 sec$. We set the initial pressure $p_o = \rho
c_o^2$, density $\rho_o=1.0$, and velocity $\ten{u}_o = U\hat{i}$. We add
an additional density damping proposed by \citet{antuono-deltasph:cpc:2010}
to the continuity equation with $\delta=0.0625$, to reduce high-frequency
pressure oscillations (see the variations of SOC schemes proposed in sec.
II.F of \cite{negi2021numerical}).

\begin{figure}[htbp]
  \centering
  \includegraphics[width=0.8\textwidth]{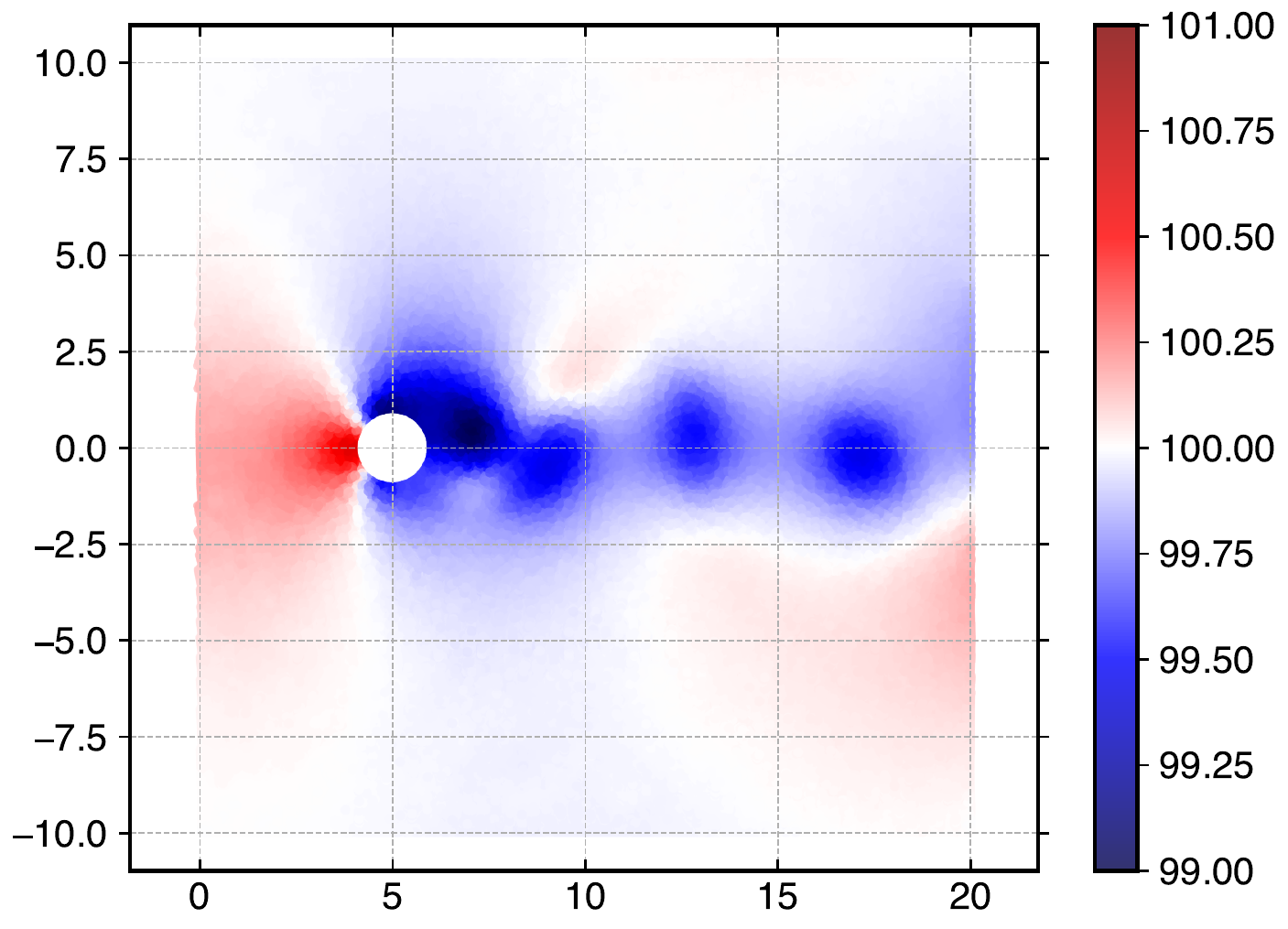}
  \caption{Average pressure at $t=100sec$.}
  \label{fig:fpc_pavg}
\end{figure}

\begin{figure}[htbp]
  \centering
  \includegraphics[width=0.8\textwidth]{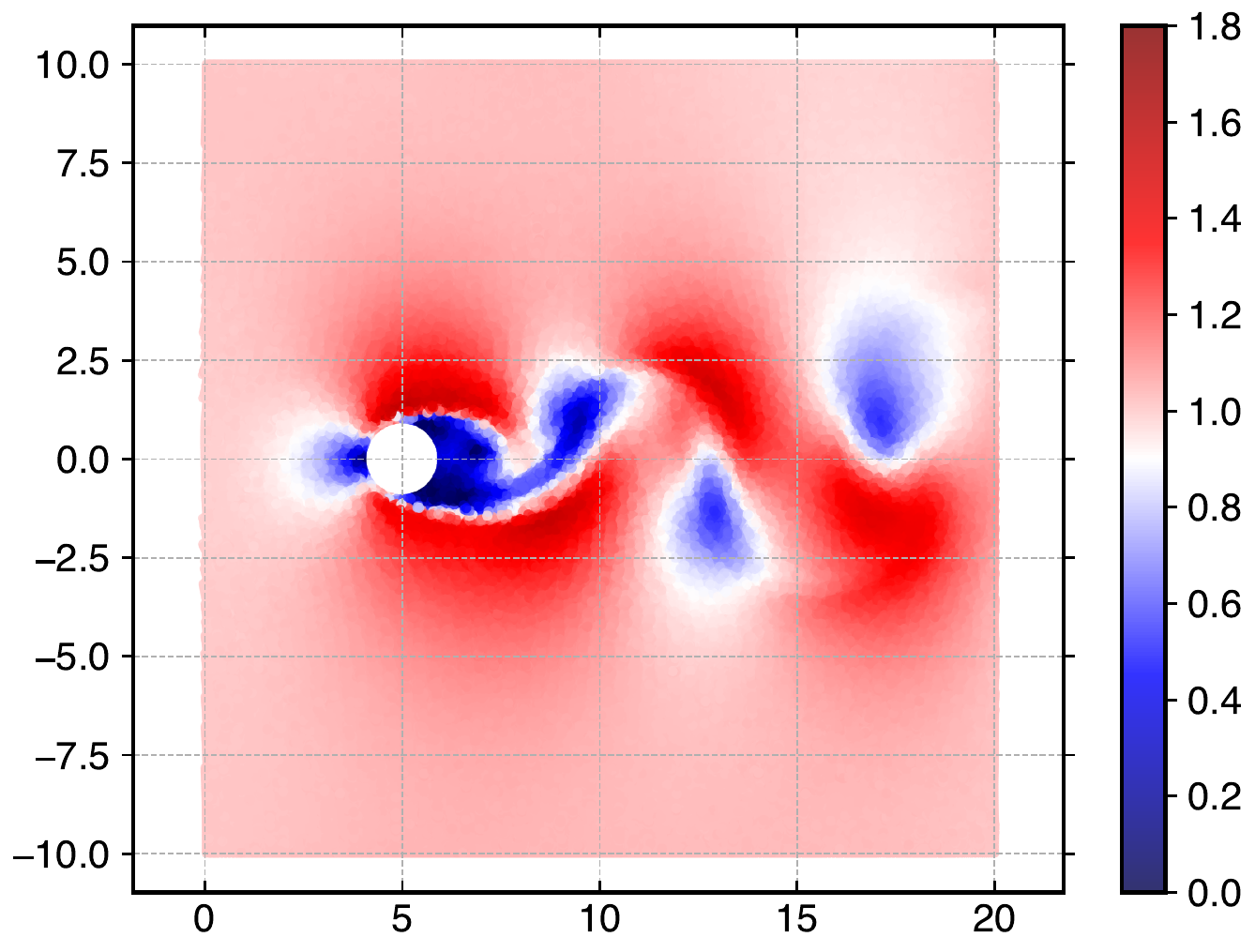}
  \caption{Velocity magnitude at $t=100sec$.}
  \label{fig:fpc_vmag}
\end{figure}

In \cref{fig:fpc_pavg} and \cref{fig:fpc_vmag}, we plot the average
pressure and velocity magnitude at $t=100sec$, respectively. We compute the
average pressure $(p_{avg})_i=\sum p_j/ N$, where the sum is taken over all
the $N$ neighbors of the $i$\textsuperscript{th} particle. Cleary, the
solution is free from any high frequency pressure oscillations. Futhermore,
the pressure in the domain remains in the vicinity of the reference
pressure $\rho c_o^2 = 100 Pa$ for the entire simulation. We also compute
the coefficient of lift $c_l$ and drag $c_d$ for the cylinder using the
method proposed in \citet{negi_nrbc_2020}.

\begin{figure}[htbp]
  \centering
  \includegraphics[width=0.8\textwidth]{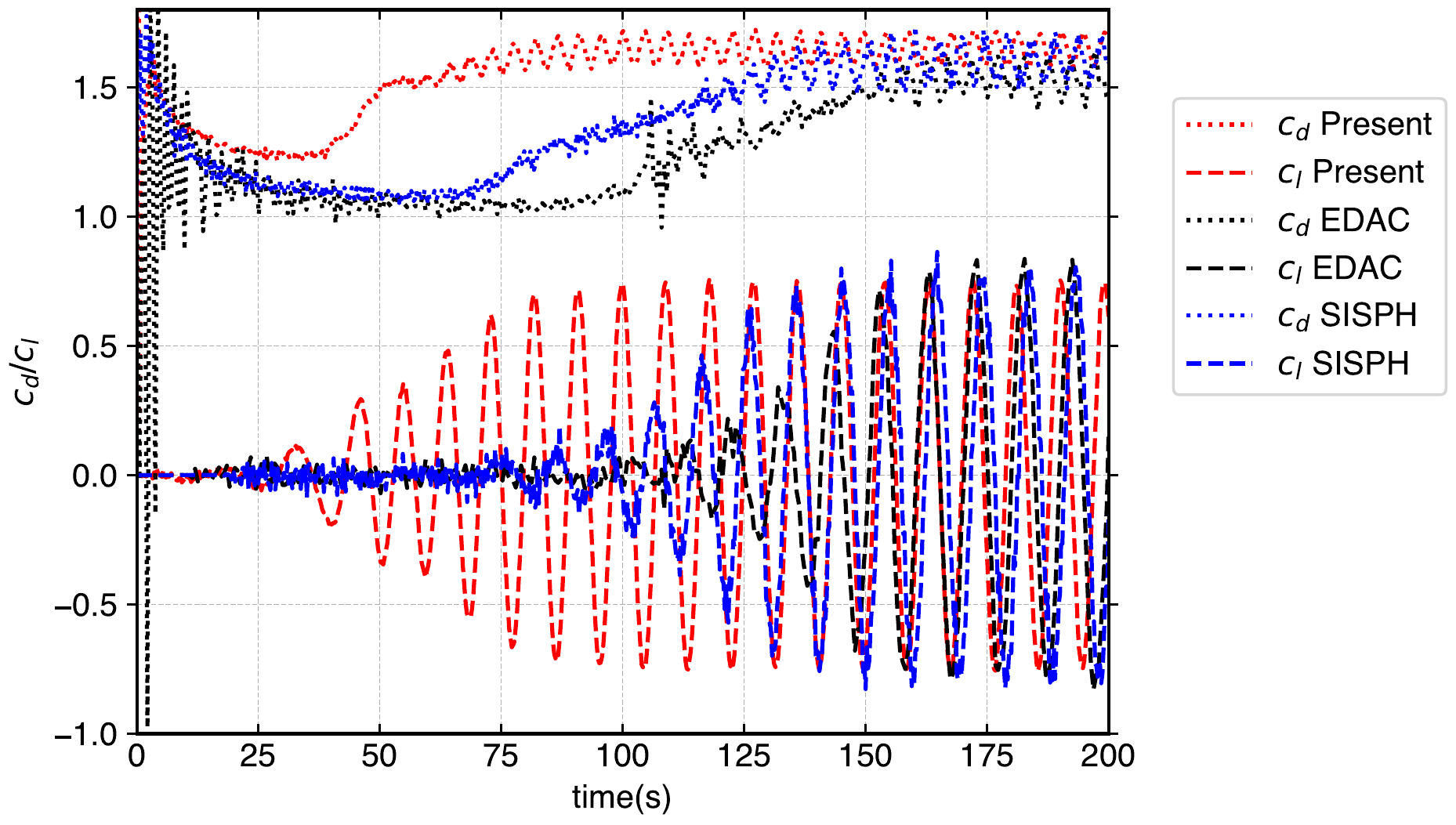}
  \caption{$c_d$ and $c_l$ variation for the flow past a cylinder with time.}
  \label{fig:fpc_cdcl}
\end{figure}

In \cref{fig:fpc_cdcl}, we plot the variation of $c_d$ and $c_l$ with time
for the present method with results of \citet{negi_nrbc_2020} referred to
as "EDAC" and \citet{muta2020efficient} referred to as "SISPH". We obtain
the mean $c_d$ value of $1.65$ and $c_l$ value of $0.74$, after the
shedding is established. These values are closer to SISPH values, where a
pressure Poisson equation is solved to obtain pressure. Futhermore, both
the drag and lift coefficient values are free from disturbances.

\section{Conclusions}
\label{sec:conclusions}

The convergence of the boundary condition implementations is a grand
challenge\cite{vacondio_grand_2020} in SPH. In order to obtain a convergent
boundary implementation, we require a scheme that must be convergent and
should have the order of error lower than that in the boundary implementation.
In this paper, we use the second order convergent scheme proposed by
\citet{negi2021numerical} to verify various boundary condition
implementations. We use the MMS and construct MS for specific boundary
conditions and domains. For a solid boundary, we verify methods for Neumann
pressure, slip, and no-slip boundary conditions. In order to cover all the
aspects of arbitrary geometries, we test the convergence on a straight,
convex, and concave boundary. In the case of open boundaries, we consider a
square domain with inlet and outlet regions simulating a wind-tunnel. We
manufactured solutions for inlets and outlets for both Neumann pressure and
velocity. Additionally, we manufactured solutions depicting waves of pressure
and velocity passing through inlet/outlets.

We show that the method proposed by \citet{marrone-deltasph:cmame:2011} is
convergent for all kinds of the domain and boundary conditions on a solid
boundary. Some other method like \citet{Adami2012,colagross2003a} are
second-order convergent in inviscid flow and packed domains. Almost all
boundary implementations are second-order on a straight boundary. In the
case of open boundaries, the mirror and simple-mirror work well in the
absence of a wave traveling through the boundary. The hybrid and do-nothing
boundaries are bounded by the $O(M^2)$, where $M$ is the Mach number of the
flow. However, in the case of a wave traveling through the domain, the
mirror and simple-mirror method diverges, and hybrid and do-nothing methods
converge with second-order accuracy. Finally, we discuss an algorithm to
apply these boundary conditions in order to get a convergent solver. We use
the method proposed by \citet{marrone-deltasph:cmame:2011} for solids and
\citet{negi_nrbc_2020} for inlet and outlet boundaries. We demonstrate the
accuracy of the proposed algorithm by solving the flow past a circular
cylinder. We achieved the accuracy close to the results obtained using
incompressible SPH solvers.

The manufactured solutions created in this paper can be used for any meshless
solver for the specified domains. In this paper, we carried out appropriately
834 simulations, which take around 70 hours, which demonstrates the efficiency
of the MMS. In the future, we would like to use the MMS to obtain a
second-order convergent adaptive solver. Since the second-order convergence
comes at the cost of performing kernel correction, the adaptivity in space and
time will compensate for this and will result in a convergent and fast WCSPH
solver.

We believe that with the identification of convergent boundary conditions
along with second order convergent SPH schemes, it should be possible to build
accurate and general purpose SPH-based solvers for the simulation of a wide
variety of incompressible and weakly compressible fluid flow simulations.
The recent advancements in adaptive resolution in SPH, will also facilitate
the efficient simulation of such problems. In the future, we would like to use
MMS to obtain a second-order convergent free surface boundary condition
which are used widely in fluid and structure applications.

\appendix

\section{Manufactured solutions}
\label{apn:mms}

\subsection{Neumann pressure boundary}
\label{apn:new_pres}

In this boundary condition, we ensure that $\nabla p \cdot \ten{n} = 0$,
where $\ten{n}$ is normal to the boundary surface. For the straight
domain the normal $\ten{n} = \ten{j}$, therefore, we can construct the MS
given by,
\begin{equation}
  \begin{split}
  u(x, y) &= \left(y - 1\right) \sin{\left(2 \pi x \right)} \cos{\left(2 \pi y \right)}\\
  v(x, y) &= - \left(y - 1\right) \sin{\left(2 \pi y \right)} \cos{\left(2 \pi x \right)}\\
  p(x, y) &= x^{2} + \cos{\left(4 \pi x \right)}.
  \end{split}
  \label{eq:mms_pres_num_d1}
\end{equation}

\begin{figure}[htbp!]
  \centering
  \includegraphics[width=\linewidth]{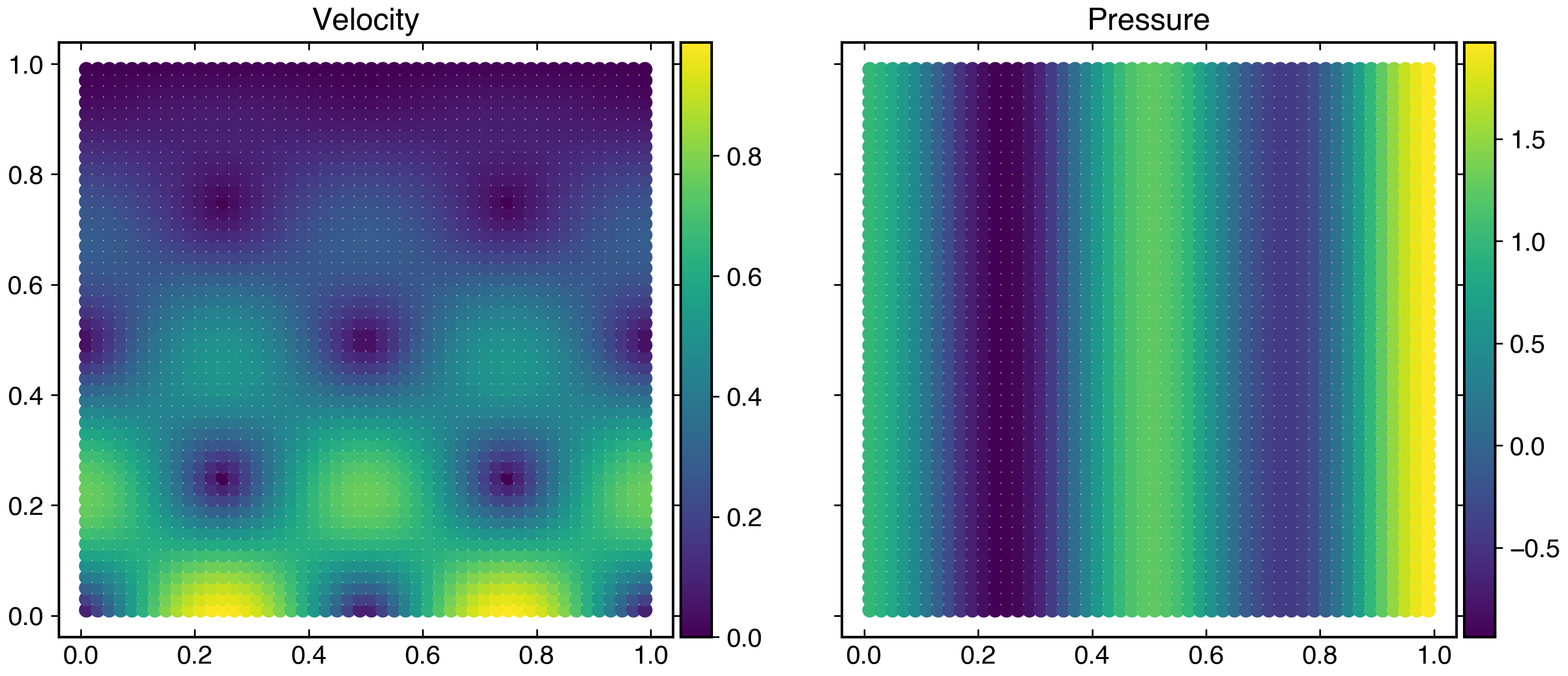}
\caption{Velocity and pressure contours on the straight domain in
\cref{fig:domains} of the MS in \cref{eq:mms_pres_num_d1}.}
  \label{fig:mms_pres_num_d1}
\end{figure}

In \cref{fig:mms_pres_num_d1}, we show the contour plot of the above MS in
straight domain.

In the case of the convex domain, the normal to the surface is given by
$\ten{n} = (x-0.5)\ten{i} + (y-0.5)\ten{j}$, therefore we can construct a
MS is given by

\begin{equation}
  \begin{split}
    u(x, y) &= \left(y - 1\right) \sin{\left(2 \pi x \right)} \cos{\left(2 \pi y \right)}\\
    v(x, y) &= - \left(y - 1\right) \sin{\left(2 \pi y \right)} \cos{\left(2 \pi x \right)}\\
    p(x, y) & = \tan^{-1} \left( \frac{(y-0.5)^2}{(x-0.5)^2} \right)
  \end{split}
  \label{eq:mms_pres_num_d5}
\end{equation}

\begin{figure}[htbp!]
  \centering
  \includegraphics[width=\linewidth]{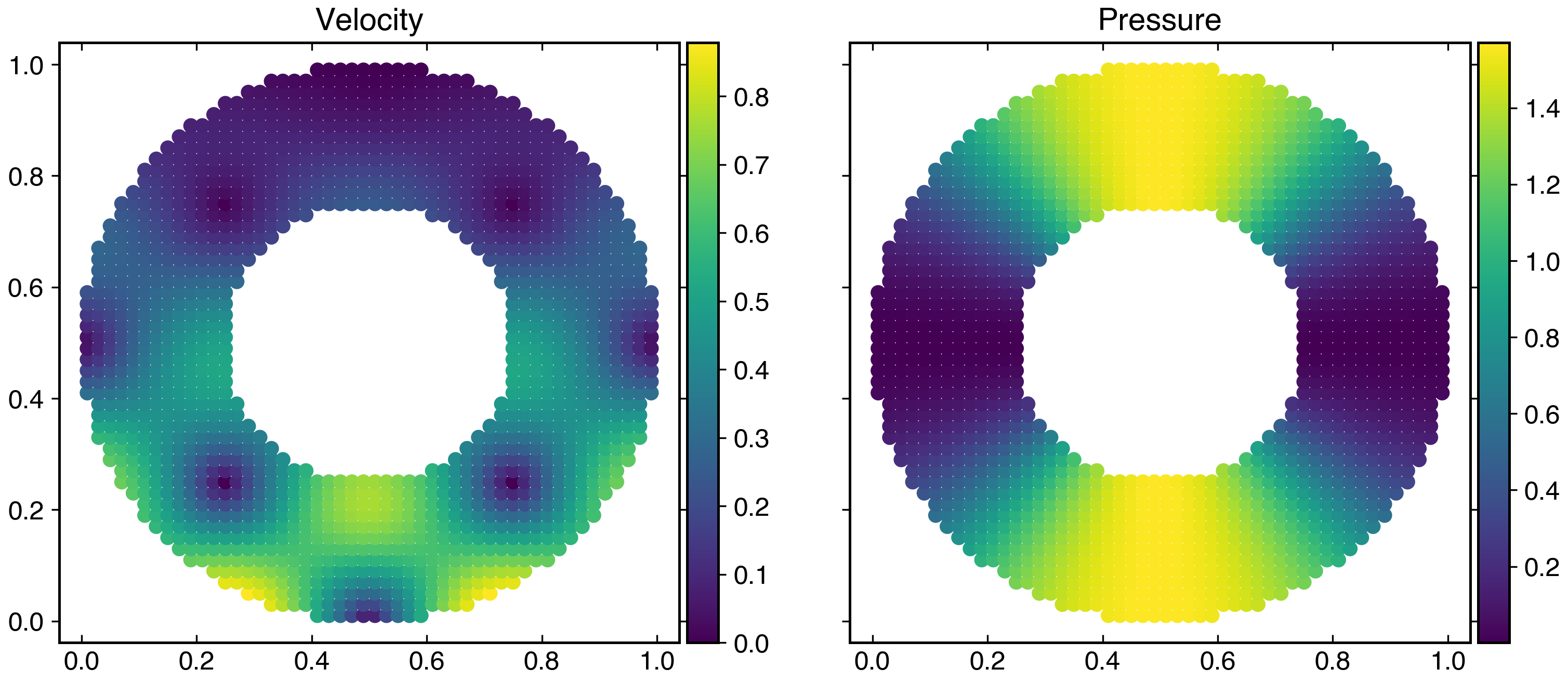}
\caption{Velocity and pressure contours on the convex/concave domain of the MS in
\cref{eq:mms_pres_num_d5}.}
  \label{fig:mms_pres_num_d5}
\end{figure}

In \cref{fig:mms_pres_num_d5}, we show the contour plot of the above MS in
the convex domain.

Since, for the concave domain, the normal remains the same, we can use the same
MS since it satisfies $\nabla p \cdot \ten{n}=0$ at the surface of
interest. We note that for the packed version of the domain in
\cref{fig:domains_pack}, we can use the same MS as used in the unpacked version
as the surface of interest is exactly the same.

\subsection{Slip boundary condition}
\label{apn:slip}

For the slip boundary condition, we ensure that $\ten{u} \cdot \ten{n}= 0$
at the boundary surface. For the straight domain, we construct the MS given by

\begin{equation}
  \begin{split}
    u(x, y) &= \left(y - 1\right) \sin{\left(2 \pi x \right)} \cos{\left(2 \pi y \right)}
      + 1\\
    v(x, y) &= \left(\left(y - 1\right)^{2}\right) \sin{\left(2 \pi y \right)}\\
    p(x, y) &= \cos{\left(4 \pi x \right)} + \cos{\left(4 \pi y \right)}\\
  \end{split}
  \label{eq:mms_slip_d1}
\end{equation}

\begin{figure}[htbp!]
  \centering
  \includegraphics[width=\linewidth]{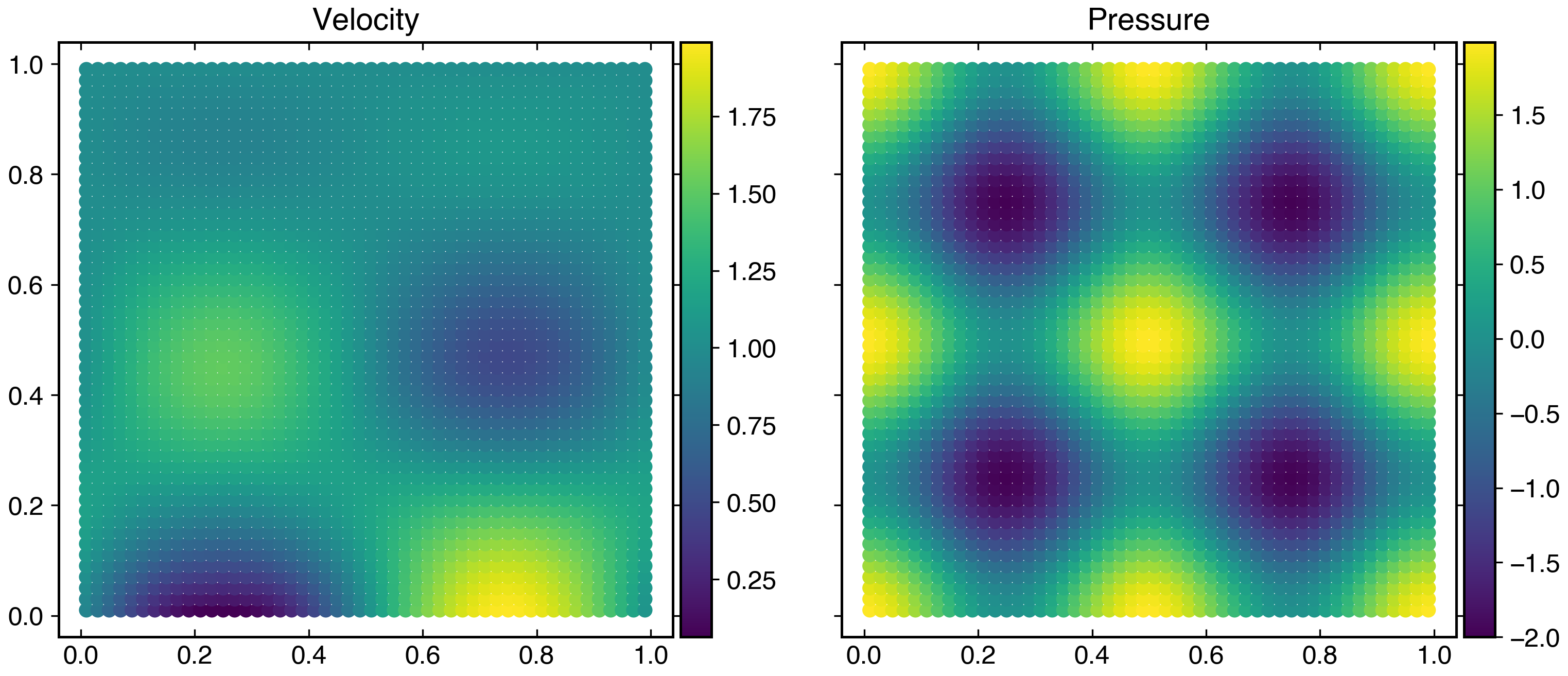}
\caption{Velocity and pressure contours on the straight domain of the
MS in \cref{eq:mms_slip_d1}.}
  \label{fig:mms_slip_d1}
\end{figure}

In \cref{fig:mms_slip_d1}, we plot the velocity and pressure contour
generated by the MS in \cref{eq:mms_slip_d1}.

\begin{figure}[htbp!]
  \centering
  \includegraphics[width=\linewidth]{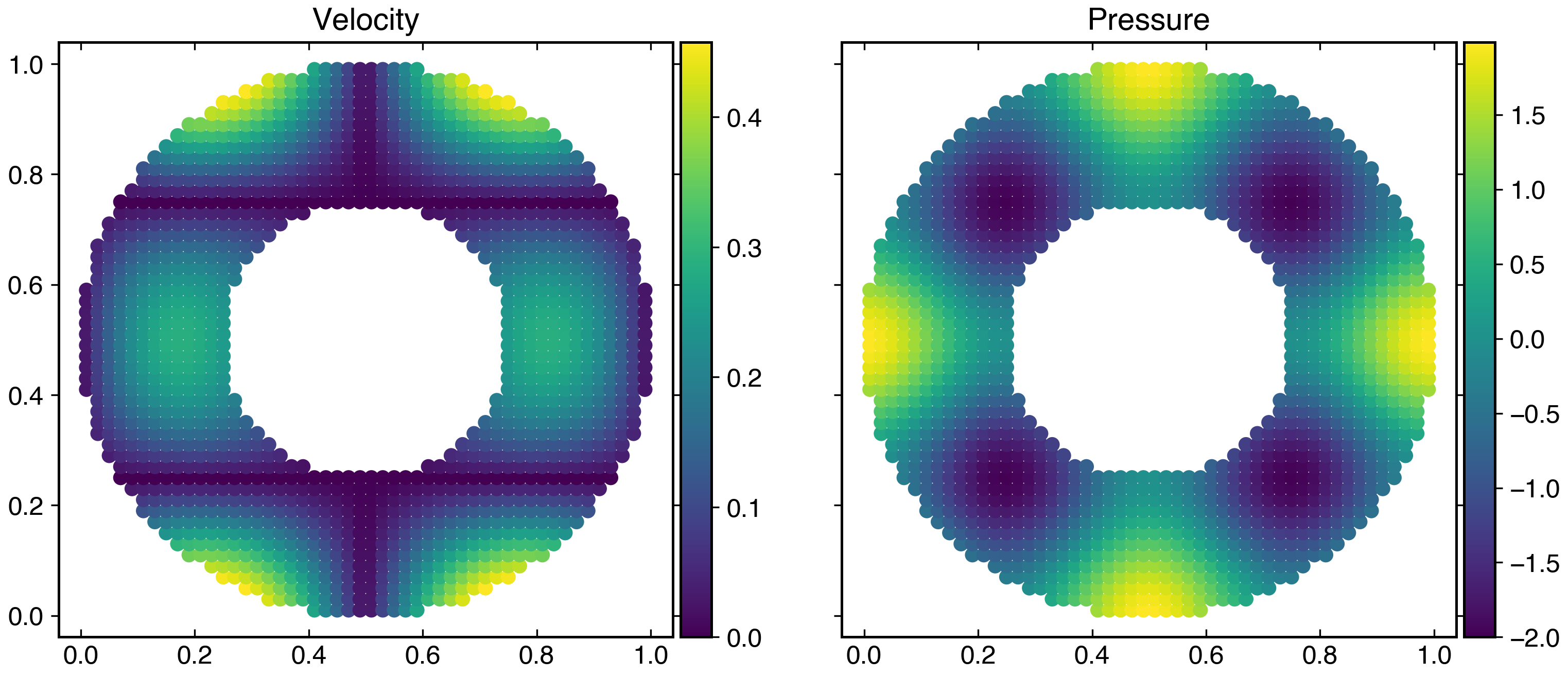}
\caption{Velocity and pressure contours on the convex/concave domain of the
MS in \cref{eq:mms_slip_d5}.}
  \label{fig:mms_slip_d5}
\end{figure}

For the convex domain the normal $\ten{n} = (x-0.5)\ten{i} +
(y-0.5)\ten{j}$, therefore we construct the MS given by
\begin{equation}
  \begin{split}
    u(x, y) &= \left(y - 0.5\right) \sin{\left(2 \pi x \right)} \cos{\left(2 \pi y \right)}\\
    v(x, y) &= - \left(x - 0.5\right) \sin{\left(2 \pi x \right)} \cos{\left(2 \pi y \right)}\\
    p(x, y) &= \cos{\left(4 \pi x \right)} + \cos{\left(4 \pi y \right)}\\
  \end{split}
  \label{eq:mms_slip_d5}
\end{equation}
such that $\ten{u} \cdot \ten{n}=0$. In \cref{fig:mms_slip_d5}, we plot the
velocity and pressure contour generated from the MS in
\cref{eq:mms_slip_d1}. We note that since for the concave as well as packed
domains, the normal remains the same therefore we can use the same MS in
\cref{eq:mms_slip_d5} for all these domains.

\subsection{No-slip boundary condition}
\label{apn:noslip}

For a no-slip boundary condition, we ensure that the velocity at the
surface is zero for a stationary wall. For the straight domain, we
construct the MS given by

\begin{equation}
  \begin{split}
    u(x, y, t) &= \left(1 - y\right)^{2} e^{- 10 t} \sin{\left(2 \pi x \right)} \cos{\left(2
      \pi y \right)}\\
    v(x, y, t) &=- \left(1 - y\right)^{2} e^{- 10 t} \sin{\left(2 \pi y \right)} \cos{\left(2
      \pi x \right)}\\
    p(x, y, t) &= \left(\cos{\left(4 \pi x \right)} + \cos{\left(4 \pi y \right)}\right) e^{-
      10 t}
  \end{split}
  \label{eq:mms_noslip_d1}
\end{equation}

\begin{figure}[htbp!]
  \centering
  \includegraphics[width=\linewidth]{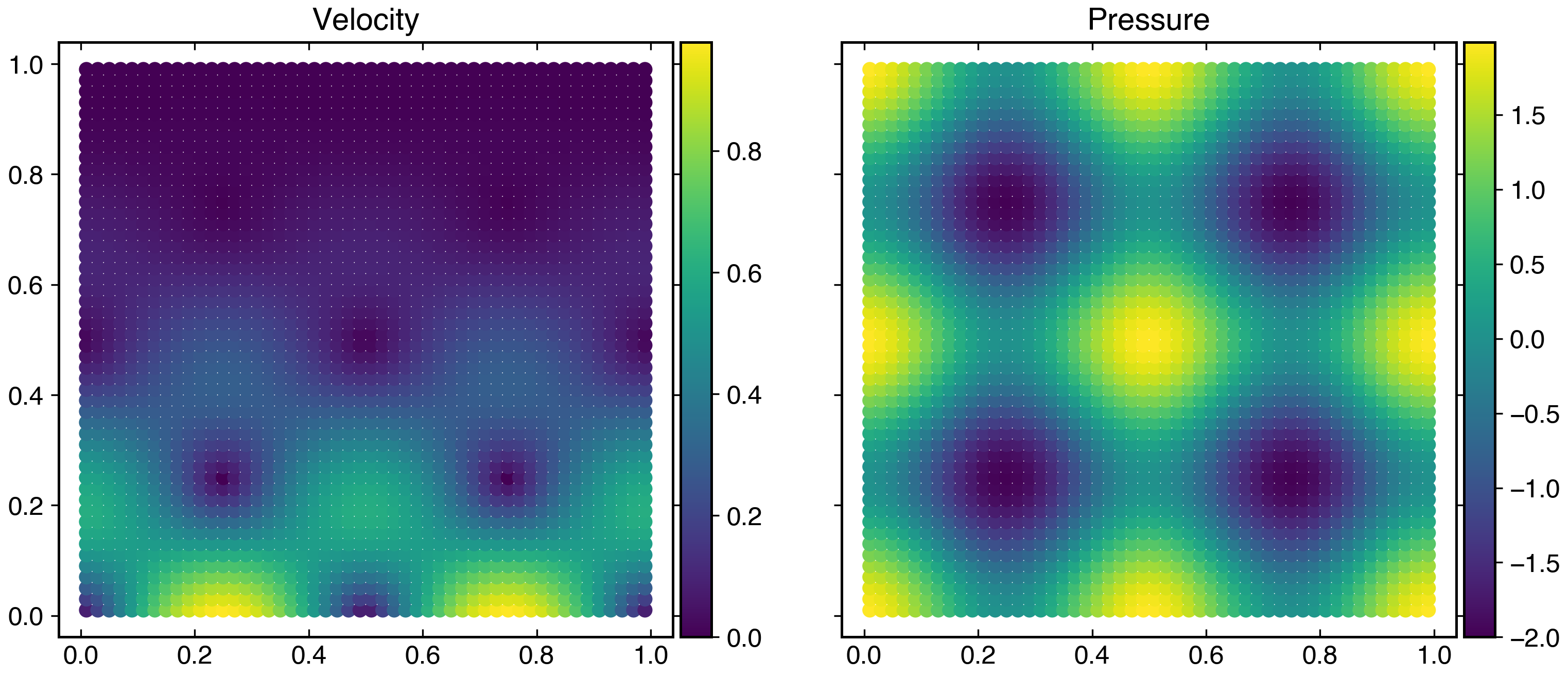}
\caption{Velocity and pressure contours on the straight domain of the MS in
\cref{eq:mms_noslip_d1}.}
  \label{fig::mms_noslip_d1}
\end{figure}

In \cref{fig::mms_noslip_d1}, we plot the contour plot for the velocity and
pressure generated by the MS in \cref{eq:mms_noslip_d1}.

In order to construct an MS for the convex domain in \cref{fig:domains}, we
construct the MS such that the velocity is zero on the inner surface of the
domain, given by

\begin{equation}
  \begin{split}
    u(x, y, t) &= \left(- \left(x - 0.5\right)^{2} - \left(y - 0.5\right)^{2} + 0.0625\right)\\
    &e^{- 10 t} \sin{\left(\pi \left(2 \left(x - 0.5\right)^{2} + 2 \left(y - 0.5\right)^{2}\right)
    \right)}\\
    v(x, y, t) &=- \left(- \left(x - 0.5\right)^{2} - \left(y - 0.5\right)^{2} + 0.0625\right)\\
    &e^{- 10 t} \cos{\left(\pi \left(2 \left(x - 0.5\right)^{2} + 2 \left(y - 0.5\right)^{2}\right)
    \right)}\\
    p(x, y, t) &= \left(\cos{\left(4 \pi x \right)} + \cos{\left(4 \pi y \right)}\right) e^{-
    10 t}
  \end{split}
  \label{eq:mms_noslip_d5}
\end{equation}

\begin{figure}[htbp!]
  \centering
  \includegraphics[width=\linewidth]{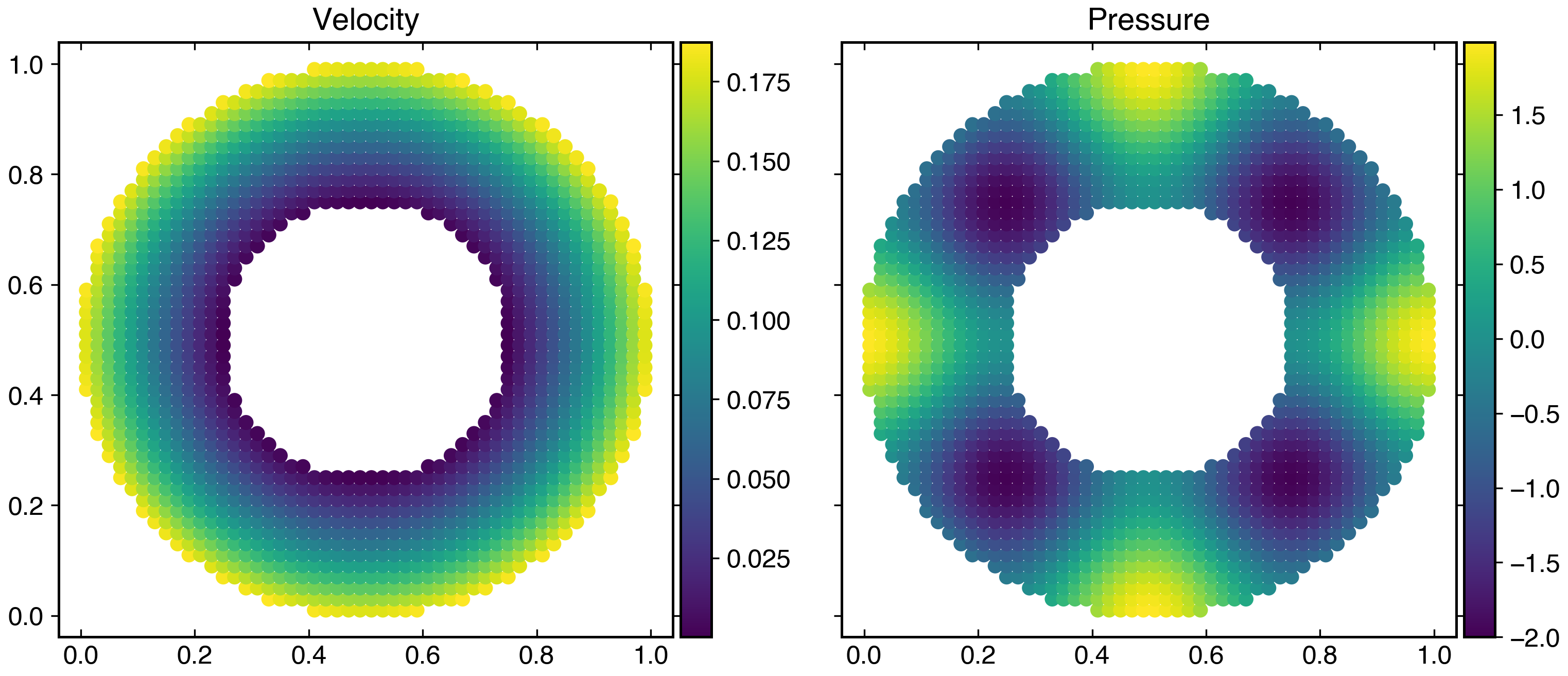}
\caption{Velocity and pressure contours on the convex domain of the MS in
\cref{eq:mms_noslip_d5}.}
  \label{fig::mms_noslip_d5}
\end{figure}

In \cref{fig::mms_noslip_d5}, we plot the velocity and pressure contour
generated from the MS in \cref{eq:mms_noslip_d5}.

In order to construct the MS for the concave domain in \cref{fig:domains}, we
make the velocity zero on the outer surfaces given by

\begin{equation}
  \begin{split}
    u(x, y, t) &= \left(- \left(x - 0.5\right)^{2} - \left(y - 0.5\right)^{2} + 0.25\right) \\
    &e^{- 10 t} \sin{\left(\pi \left(2 \left(x - 0.5\right)^{2} + 2 \left(y - 0.5\right)^{2}\right)
    \right)}\\
    v(x, y, t) &=- \left(- \left(x - 0.5\right)^{2} - \left(y - 0.5\right)^{2} + 0.25\right)\\
    &e^{- 10 t} \cos{\left(\pi \left(2 \left(x - 0.5\right)^{2} + 2 \left(y - 0.5\right)^{2}\right)
    \right)}\\
    p(x, y, t) &= \left(\cos{\left(4 \pi x \right)} + \cos{\left(4 \pi y \right)}\right) e^{-
    10 t}
  \end{split}
  \label{eq:mms_noslip_d6}
\end{equation}

\begin{figure}[htbp!]
  \centering
  \includegraphics[width=\linewidth]{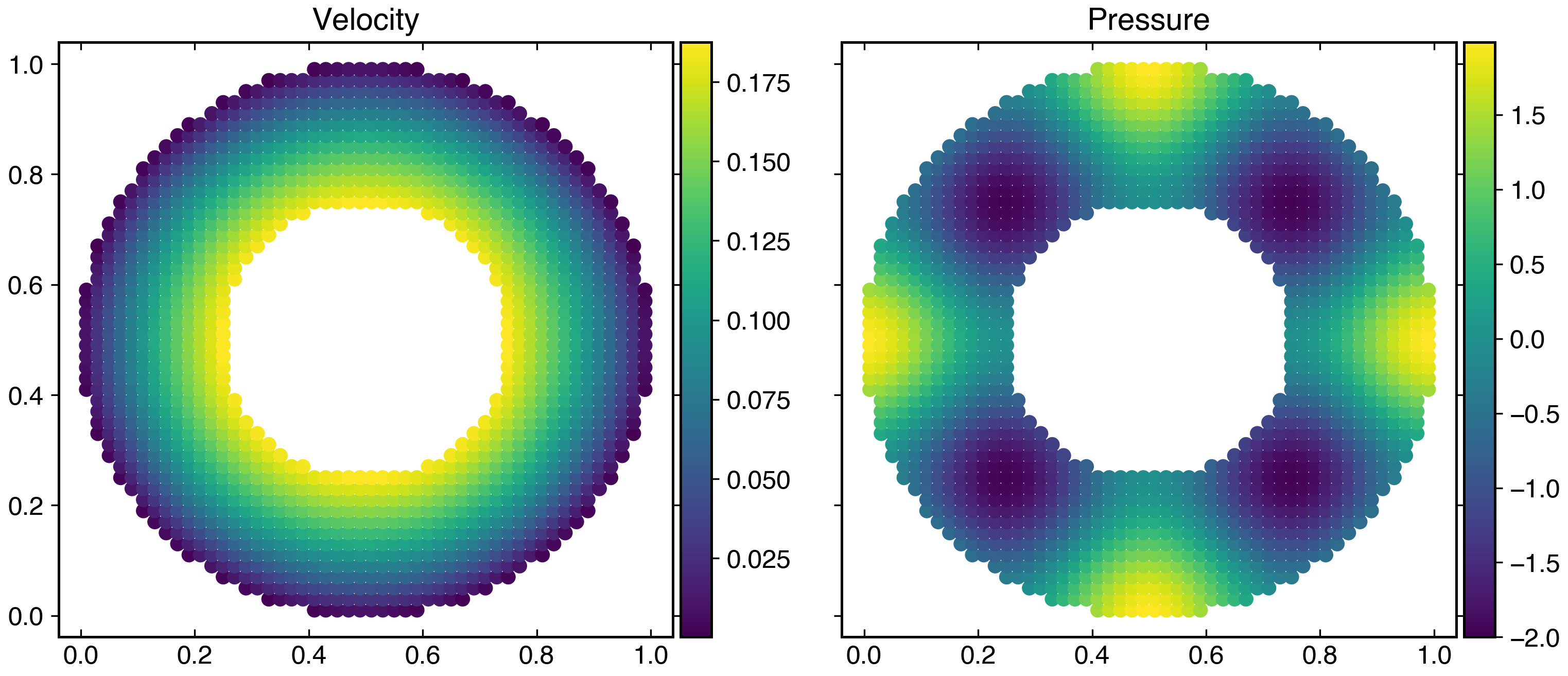}
\caption{Velocity and pressure contours on the concave domain of the MS in
\cref{eq:mms_noslip_d6}.}
  \label{fig::mms_noslip_d6}
\end{figure}

In \cref{fig::mms_noslip_d6}, we plot the contour for velocity and pressure
generated from the \cref{eq:mms_noslip_d6} for the concave domain. We note
that the MS described remains the same for the corresponding packed version
of the domains.

\subsection{Inlet and outlet velocity boundary}
\label{apn:io_vel}

At the inlet, we make sure that $\nabla \ten{u} \cdot \ten{n} = 0$.`' Since
the inlet is usually straight. We consider one type of inlet with constant
normal $\ten{n} = -\ten{i}$ similarly outlet with normal $\ten{n}
=\ten{i}$. We use the MS given by
\begin{equation}
  \begin{split}
    u(x, y, t) &= y \left(y - 1\right) e^{- 10 t} \cos{\left(2 \pi y \right)} + 1\\
    v(x, y, t) &= - x^{2} \left(x - 1\right)^{2} e^{- 10 t} \sin{\left(2 \pi y \right)}\\
    p(x, y, t) &= \left(\cos{\left(4 \pi x \right)} + \cos{\left(4 \pi y \right)}\right) e^{-
      10 t}
  \end{split}
  \label{eq:mms_io_vel}
\end{equation}

\begin{figure}[htbp!]
  \centering
  \includegraphics[width=\linewidth]{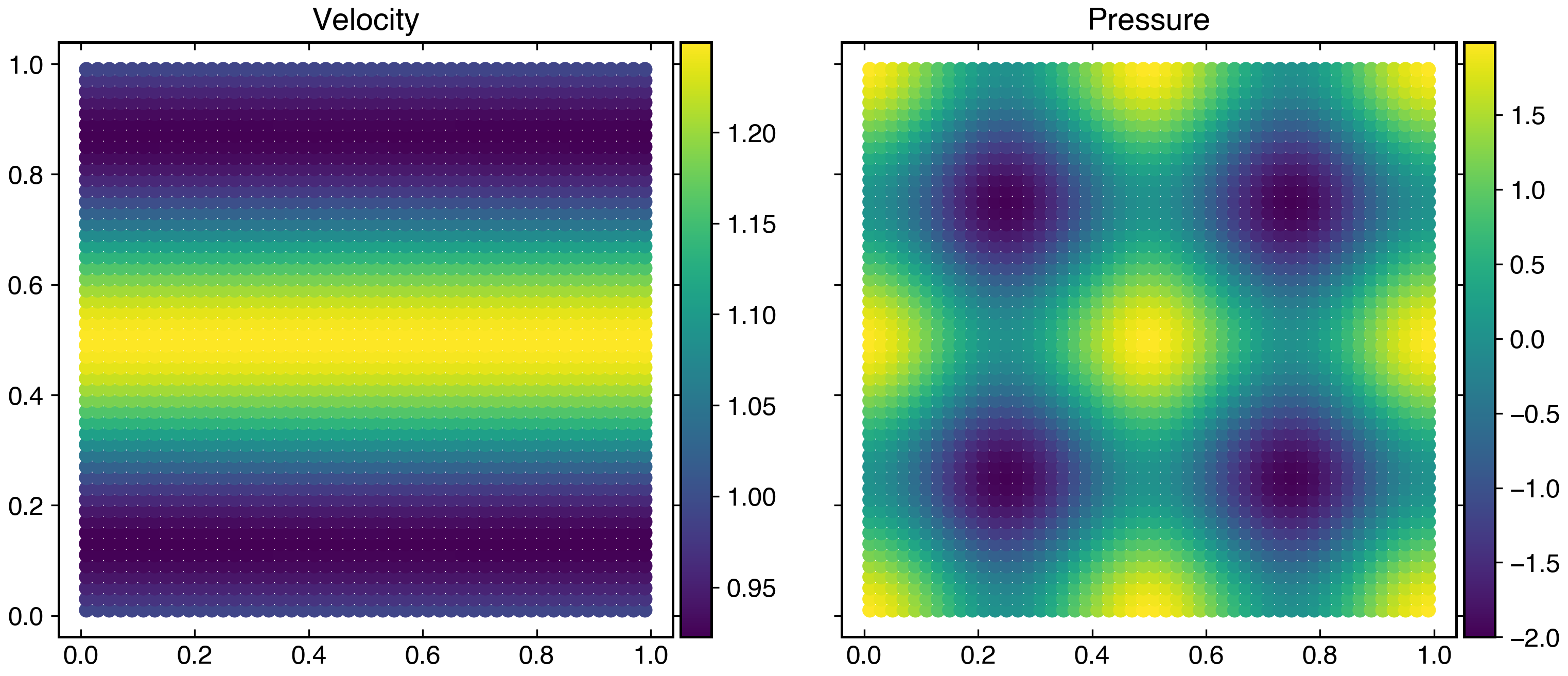}
\caption{Velocity and pressure contours on the domain in \cref{fig:io} of the MS in
\cref{eq:mms_io_vel}.}
  \label{fig::mms_io_vel}
\end{figure}

In the \cref{fig::mms_io_vel}, we plot the velocity and pressure contour
for the MS in \cref{eq:mms_io_vel}.

Additionally, we also simulate the wave passing through the inlet and
outlet. We also must satisfy the boundary condition. For the inlet, we
construct the MS given by

\begin{equation}
  \begin{split}
    u(x, y, t) &= x^{2} y \left(y - 1\right) e^{- 200 \left(x-0.1-40 t \right)^{2}}
    \cos{\left(2 \pi y \right)} + 1\\
    v(x, y, t) &= 0.0\\
    p(x, y, t) &= \cos{\left(4 \pi x \right)} + \cos{\left(4 \pi y \right)}
  \end{split}
  \label{eq:mms_in_vel}
\end{equation}

\begin{figure}[htbp!]
  \centering
  \includegraphics[width=\linewidth]{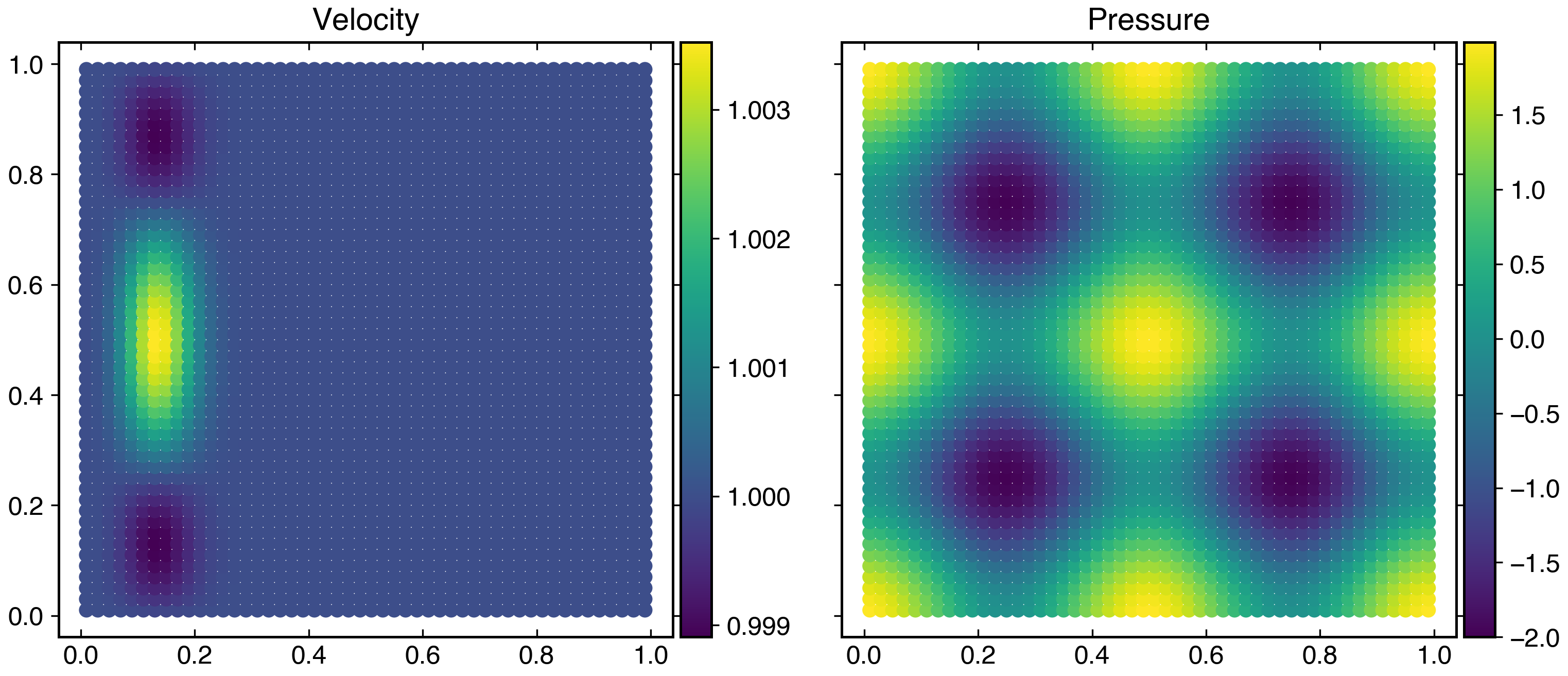}
\caption{Velocity and pressure contours on the domain in \cref{fig:io} of the MS in
\cref{eq:mms_io_vel}.}
  \label{fig::mms_in_vel}
\end{figure}

In \cref{fig::mms_in_vel}, we plot the velocity and pressure contour for
the MS in \cref{eq:mms_in_vel}.

We construct the wave of velocity passing through the outlet given by

\begin{equation}
  \begin{split}
    u(x, y, t) &= (x-1)^{2} y \left(y - 1\right) e^{- 200 \left(x-0.9+40 t \right)^{2}}
    \cos{\left(2 \pi y \right)} + 1\\
    v(x, y, t) &= 0.0\\
    p(x, y, t) &= \cos{\left(4 \pi x \right)} + \cos{\left(4 \pi y \right)}
  \end{split}
  \label{eq:mms_out_vel}
\end{equation}

\begin{figure}[htbp!]
  \centering
  \includegraphics[width=\linewidth]{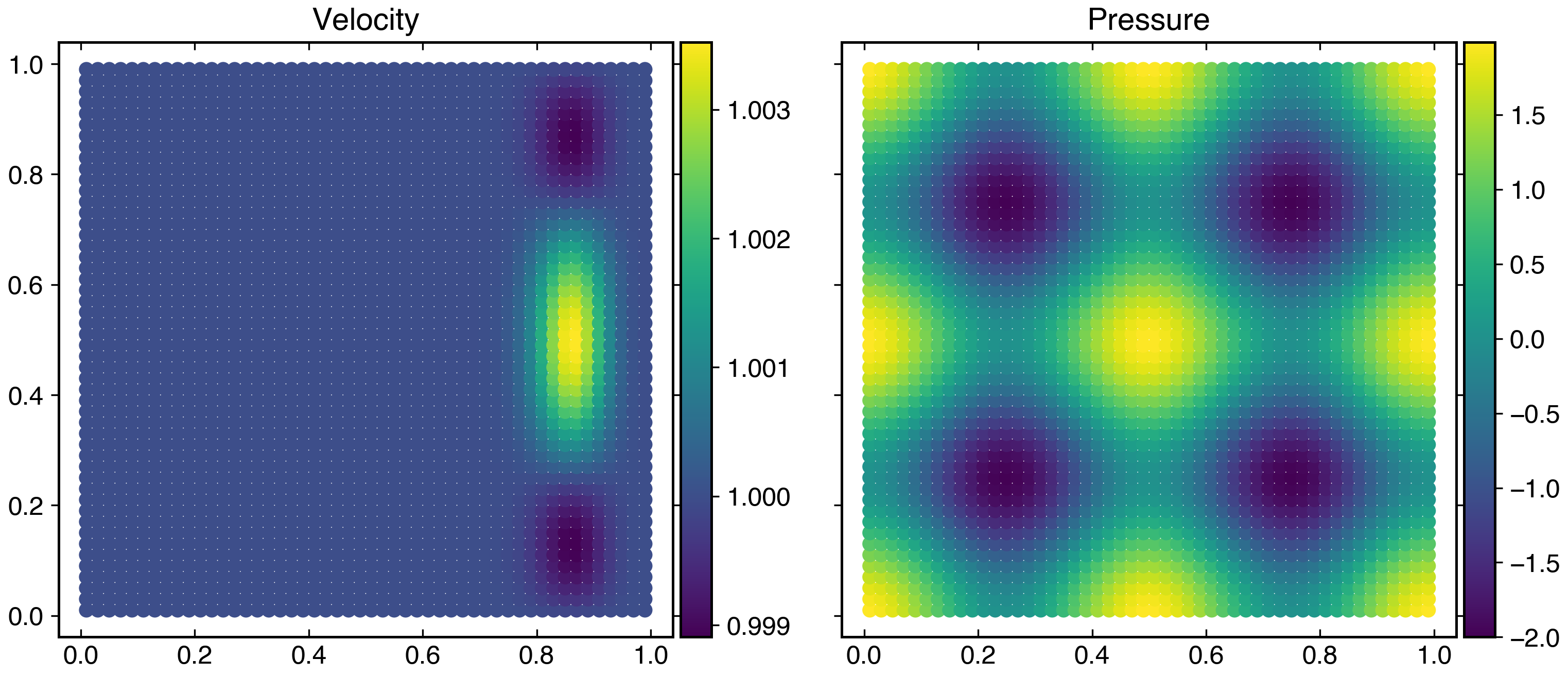}
\caption{Velocity and pressure contours on the domain in \cref{fig:io} of the MS in
\cref{eq:mms_io_vel}.}
  \label{fig::mms_out_vel}
\end{figure}

In \cref{fig::mms_out_vel}, we plot the velocity and pressure generated by
the MS in \cref{eq:mms_out_vel}.

\subsection{Inlet and outlet pressure boundary}
\label{apn:io_pres}

At the inlet, for pressure ,we make sure that $\nabla p \cdot \ten{n} =
0$. For inlet as well as the outlet, we use the MS given by
\begin{equation}
  \begin{split}
    u(x, y, t) &= y \left(y - 1\right) e^{- 10 t} \cos{\left(2 \pi y \right)} + 1\\
    v(x, y, t) &= - x \left(x - 1\right) e^{- 10 t} \sin{\left(2 \pi y \right)}\\
    p(x, y, t) &= y \left(y - 1\right) e^{- 10 t} \cos{\left(2 \pi y \right)}
  \end{split}
  \label{eq:mms_io_pres}
\end{equation}

\begin{figure}[htbp!]
  \centering
  \includegraphics[width=\linewidth]{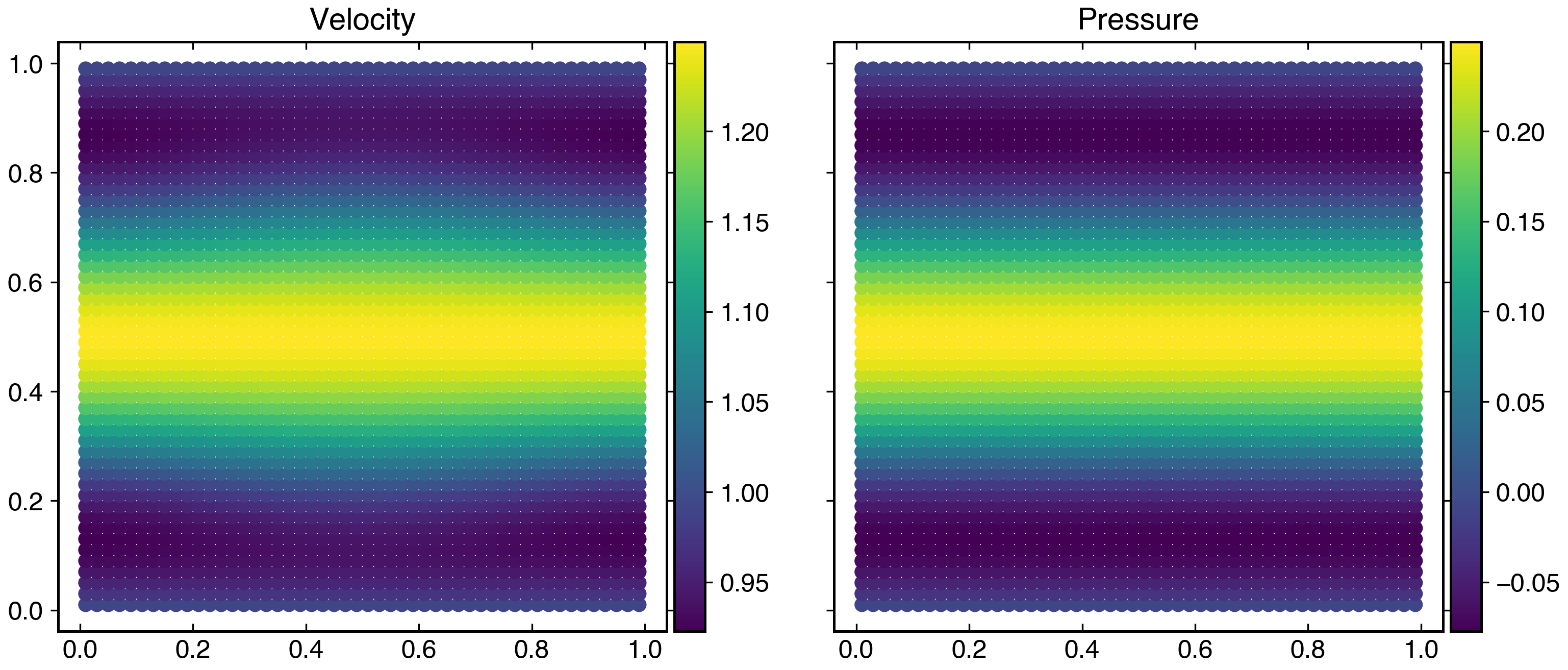}
\caption{Velocity and pressure contours on the domain in \cref{fig:io} of the MS in
\cref{eq:mms_io_pres}.}
  \label{fig::mms_io_pres}
\end{figure}

In \cref{fig::mms_io_pres}, we plot the velocity and pressure contour
generated from the MS in \cref{eq:mms_io_pres}.

In order to simulate a pressure wave passing through both inlet and outlet,
we construct MSs with pressure moving with the artificial speed of sound. For
the inlet, we construct the MS given by

\begin{equation}
  \begin{split}
    u(x, y, t) &= y \left(y - 1\right) \cos{\left(2 \pi y \right)} + 1\\
    v(x, y, t) &= 0.0\\
    p(x, y, t) &= x^{2} e^{- 200 \left(x-0.1-40 t \right)^2} \cos{\left(2
    \pi y \right)}
  \end{split}
  \label{eq:mms_in_pres}
\end{equation}

\begin{figure}[htbp!]
  \centering
  \includegraphics[width=\linewidth]{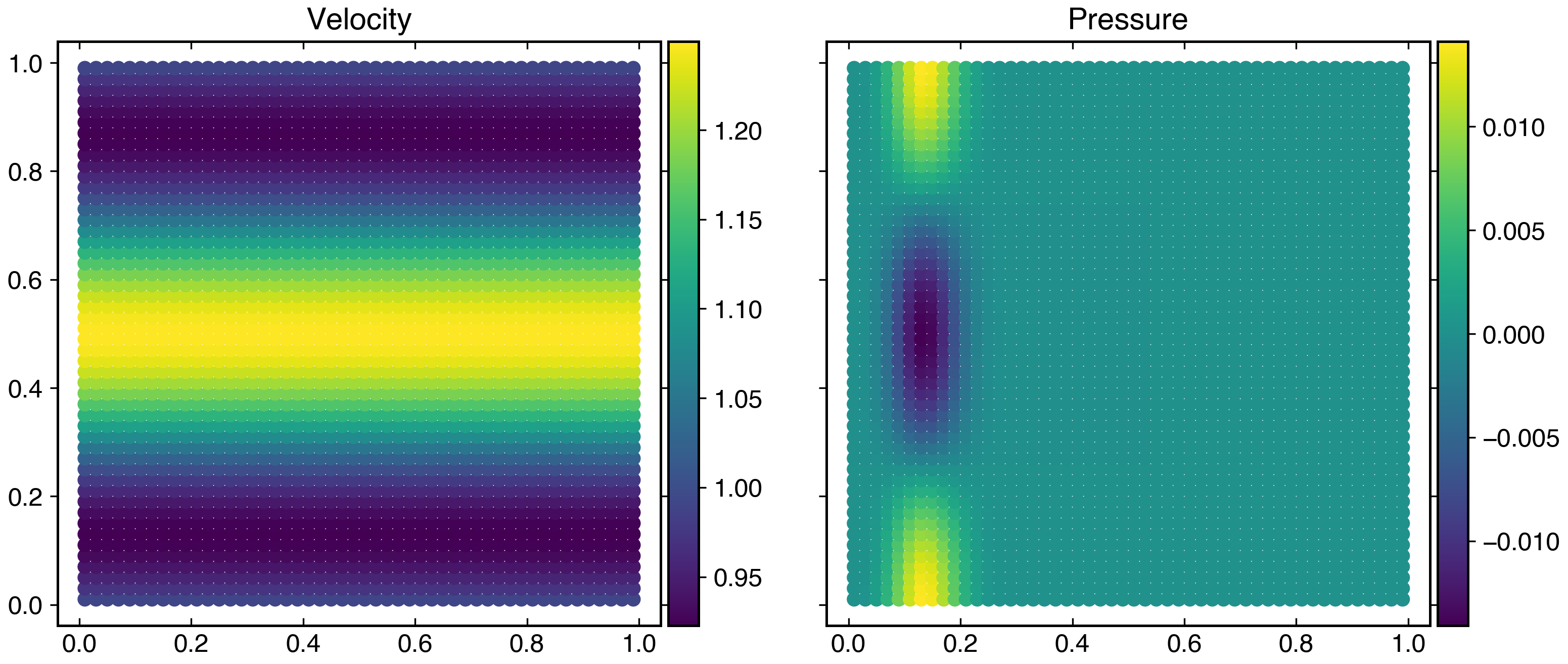}
\caption{Velocity and pressure contours on the domain in \cref{fig:io} of the MS in
\cref{eq:mms_io_pres}.}
  \label{fig::mms_in_pres}
\end{figure}

In \cref{fig::mms_in_pres}, we plot the velocity and pressure generated
from the MS in \cref{eq:mms_in_pres}.

In the case of the outlet, we construct the MS given by

\begin{equation}
  \begin{split}
    u(x, y, t) &= y \left(y - 1\right) \cos{\left(2 \pi y \right)} + 1\\
    v(x, y, t) &= 0.0\\
    p(x, y, t) &= (x-1)^{2} e^{- 200 \left(x-0.9+40 t \right)^2} \cos{\left(2
    \pi y \right)}
  \end{split}
  \label{eq:mms_out_pres}
\end{equation}

\begin{figure}[htbp!]
  \centering
  \includegraphics[width=\linewidth]{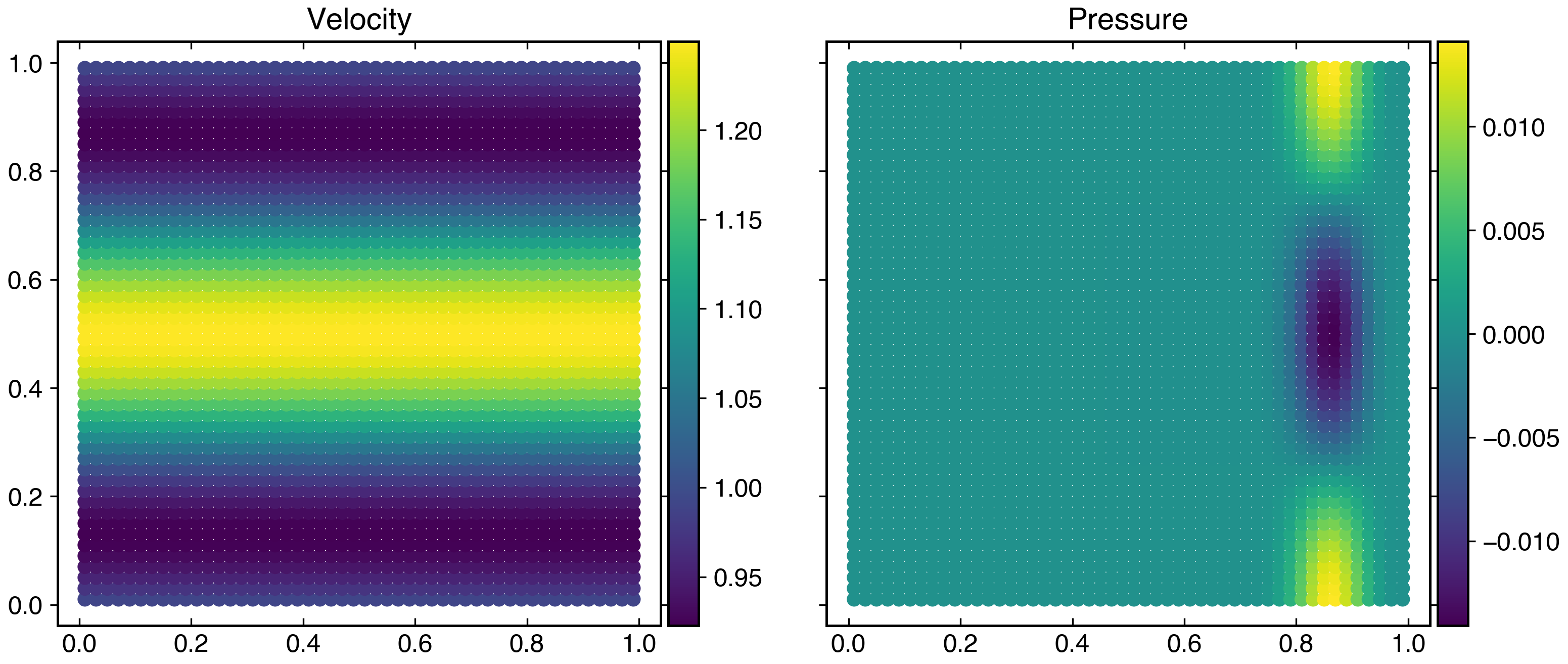}
\caption{Velocity and pressure contours on the domain in \cref{fig:io} of the MS in
\cref{eq:mms_io_pres}.}
  \label{fig::mms_out_pres}
\end{figure}

In \cref{fig::mms_out_pres}, we plot the velocity and pressure due to MS in
\cref{eq:mms_out_pres}.

\bibliographystyle{model6-num-names}
\bibliography{references}

\end{document}